\documentclass[12pt,a4paper,leqno]{amsart}
\usepackage{amsfonts,amssymb,amsmath,amstext,amscd,latexsym,euscript,amsthm}
\advance\textwidth by 3.2cm 
\advance\oddsidemargin by -1.6cm \advance\evensidemargin by -1.6cm
\input{colordvi}
\newtheoremstyle{plain}{2pt}{2pt}{\normalfont\sf}{\parindent}{\normalfont\bf}{.\,\,}{.0em}{}
\newtheoremstyle{definition}{2pt}{2pt}{\normalfont\sl}{\parindent}{\normalfont\bf}{.\,\,}{.0em}{}
\newtheoremstyle{proof}{2pt}{2pt}{\normalfont}{\parindent}{\normalfont\scshape}{.}{.0em}{}
\theoremstyle{definition}
\newtheorem{definition}{Definition}[section]
\newtheorem{example}{Example}[section]
\newtheorem{remark}{Remark}[section]
\newtheorem*{definition*}{Definition}
\newtheorem*{example*}{Example}
\newtheorem*{remark*}{Remark}
\theoremstyle{proof}

\theoremstyle{plain}
\newtheorem{theorem}{Theorem}[section]
\newtheorem{lemma}{Lemma}[section]
\newtheorem{proposition}{Proposition}[section]
\newtheorem{corollary}{Corollary}[section]

\newtheorem*{theorem*}{Theorem}
\newtheorem*{lemma*}{Lemma}
\newtheorem*{proposition*}{Proposition}
\newtheorem*{corollary*}{Corollary}
\newtheorem*{statement*}{Statement}
\newtheorem*{problem*}{Problem}
\newtheorem*{hypothesis*}{Hypothesis}
\numberwithin{equation}{section}
 \DeclareMathOperator{\span1}{span}
 \DeclareMathOperator{\supp}{supp}

\newcommand{\ex}{\exists}
\newcommand{\mmb}[1]{\mathbb{#1}}
\newcommand{\mmf}[1]{\mathbf{#1}}
\newcommand{\mmm}[1]{\mathrm{#1}}
\newcommand{\mcl}[1]{\mathcal{#1}}

\newcounter{num}

\newcounter{rlab}

\newcommand{\ov}[1]{\overline{#1}}

\def\Nat{{\mathbb N}}
\begin{document}

\mbox{ }

\mbox{ }

\mbox{ }

\mbox{ }

\title{Chernoff and Trotter-Kato theorems for locally convex spaces.}
\author{Neklyudov A. Y.}

 \email{neklyudov.alex@gmail.com}


\begin{abstract} We develop new approach for studying the abstract Cauchy
problem
 $\dot{x}={\mmm{A}}x$, $x(0)=x_0\in\mcl{D}(\mmm{A})$, for linear operator  $\mmm{A}$ defined on locally convex space $\mmf{X}$.
 This approach was firstly introduced in the paper
 "Chernoff and Trotter type product formulas"\, to study the problem for Banach spaces.
 In this paper
 we not only generalize the results of the previous paper to more general topological spaces
 but also get new results for Banach spaces.
 In particular, we prove the "local"\, extension of Chernoff-Trotter-Kato type theorems.
 Applying this result, we prove  Chernoff, Lie-Trotter and Trotter-Kato theorems for locally convex
spaces. Also we find necessary and sufficient conditions for the
validity of the Chernoff and Trotter product formulas.
\end{abstract}\maketitle
\maketitle
\tableofcontents

\section{Introduction}
The paper is devoted to extensions of Chernoff, Lie-Trotter and
Trotter-Kato theorems for locally convex spaces. In the paper we
have developed the approach which was firstly introduced in
\cite{Nekl2} for studying the abstract Cauchy problem
$\dot{x}={\mmm{A}}x$, $x(0)=x_0\in\mcl{D}(\mmm{A})$ in Banach
spaces. As a result, we not only generalize the results of that
paper to more general topological spaces,
 but also get new results for Banach spaces. In particular, we find new proofs of Chernoff, Lie-Trotter and Trotter-Kato theorems~(\cite{Engel}) and,
 moreover, prove the local version of these theorems.

Chernoff, Lie-Trotter and Trotter-Kato theorems are classical
results in the theory of $C_0$-semigroups in Banach spaces
(\cite{Dav}, \cite{Engel}). Initially the "resolvent"\, technique
was used to prove these theorems. So, the development of this
technique for the theory of semigroups in locally convex spaces was
one of the central problems in works
~\cite{Grosu},~\cite{Komatsu},~\cite{Komura},~\cite{Ouchi},~\cite{Schwartz},~\cite{Yosida}
and others. Let us draw attention to~\cite{Albanese0}. There,
applying the "resolvent"\, technique developed in~\cite{Ouchi}, the
extensions of Trotter-Kato theorem were proved for locally
equicontinuous semigroups defined on sequentially complete locally
convex spaces (Theorems 15, 17). However, unlike our variant of the
extension of Trotter-Kato theorem (Theorem~\ref{t-k}),
 these results are not applicable to arbitrary equicontinuous
semigroups. Indeed,
it suffices to construct a generator $\mmm{Z}$ of the
$lC_0$-semigroups $\mmm{T}$, such that there doesn't exist
 $\lambda\in\mmb{C}$, for which the image of the operator
 $\lambda\mmm{I}-\mmm{Z}$ is dense in $\mmm{X}$ (see the example of such semigroup in Example~\ref{exa}).

Chernoff theorem \cite{Chernoff} is the extension of Lie-Trotter
theorem \cite{Trotter} and can be considered as one of the central
theorems in the theory of semigroups. Chernoff theorem is applied to
find representations of solutions for many equations in mathematical
physics such as Shr\"{o}dinger equation and heat equation. In
particular, E. Nelson \cite{Nelson} applied Lie-Trotter formula for
representing the solution of Shr\"{o}dinger equation by Feynman path
integrals in phase space. Chernoff theorem was used in \cite{Witt}
  to find the solution of heat equation on a compact Riemannian manifold  and it was also used in \cite{Hamil}
  to represent the solution of Shr\"{o}dinger equation by Feynman path integrals.

 In the present paper we obtain Chernoff type formula for the local solution of the abstract Cauchy problem
\begin{eqnarray}
 \dot{x} =\mmm{A}x,\,\,
x(0)=x_0\in\mcl{D}(\mmm{A}),\label{new:1}
\end{eqnarray} for any  linear operator $\mmm{A}\in\mcl{Z}_{\mmf{X}}$ on a locally convex space $\mmf{X}$ (see Definitions~\ref{deft:20},~\ref{deft:21} and Theorem~\ref{chern1}). 
Applying this result,
we find necessary and sufficient conditions for the existence of the
local solution of the abstract Cauchy problem~\ref{new:1} (see
Theorems~\ref{vvmain:10a},~\ref{vvt:54} and Corollary~\ref{t:55}).
Moreover, we prove the extensions of Chernoff, Lie-Trotter and
Trotter-Kato theorems for sequentially complete locally convex
spaces (see Theorems~\ref{chern}-\ref{t:22}). As a consequence of
these results, we find necessary and sufficient conditions that the
closure of a linear operator $\mmm{A}$ is the generator of the
locally equicontinuous semigroup (see
Theorems~\ref{vt:4}-\ref{vt:54},
Corollaries~\ref{vmain:10}-\ref{vt:55} and Remark~\ref{ref:9.1}).
 In addition, we show that if there exists a local solution of the abstract Cauchy problem
\begin{eqnarray}\dot{x}=\mmm{Z}x,\,\, x(0)=x_0\in\mcl{D}(\mmm{Z}),
\label{vvv:1}\end{eqnarray}
where $\mmm{Z}$ is a densely defined dissipative linear operator in a
Hilbert space  $\mmf{X}$, then the solution can be represented by
Chernoff type formula (Corollary~\ref{chern1hilbert}).  In this case
Chernoff theorem is not applicable because the closure of $\mmm{Z}$
may not be the generator of the $C_0$-semigroup. By
Theorem~\ref{chern1}, we also easily get representations of the
local solution for Shr\"{o}dinger equation by Feynman path
integrals (\cite{Hamil}). In conclusion, note that the proof of
Theorem~\ref{chern1} can be used to find an estimate of the rate of convergence of the
approximation to the local solution
of the abstract Cauchy problem \eqref{new:1}.

\textbf{Acknowledgments}. The author expresses the gratitude to his
PhD advisor Professor O.~G.~Smolyanov for his attention to the
work and helpful discussions. E.~T.~Shavgulidze, N.~N.~Shamarov,
S.~V.~Kozyrev were among the many colleagues who read the paper and
helped to improve it by their comments. Also I would like to thank Z.~Brzezniak for his attention to the
work. Finally, I'm grateful to
M.~Y.~Neklyudov for reading the final version of the paper and
making many valuable remarks.

\section{Preliminaries.}\label{1}
Let $\mmf{E}$ be a locally convex space on the field $\mmb{T}\in\{\mmb{R},\mmb{C}\}$ with topology generated by semi-norms
$\|\cdot\|_{\alpha},\,\,\alpha\in\Omega;$ $\mcl{L}(\mmf{E})$
be the space of all bounded linear operators in  $\mmf{E}$  with the topology of pointwise convergence (strong operator topology); $\mmm{I}$ be the identity operator in $\mmf{E}$.
 For any function $\mmm{C}$ we denote $\mcl{D}(\mmm{C})$ the domain of
$\mmm{C}$. Linear
 subspaces of
$\mmf{E}$ will be considered as the locally convex spaces with topology generated by the natural
semi-norms (i.e. the semi-norms inherited from the locally convex
space $\mmf{E}$). For any set $\mmf{A}\subset\mmf{E}$ let the
function
$\|\cdot\|_{\mmf{A}}:\mmf{E}\mapsto[0,\,\infty)\cup\{\infty\}$
be
defined by the equality
$$\|x\|_{\mmf{A}}=\begin{cases}
\sup{\{t>0|x/t\notin\mmf{A}\}},&\text{if
$\{t>0|x/t\notin\mmf{A}\}\neq{\emptyset}$;}\\
0,&\text{otherwise,}
\end{cases}$$ where $x\in\mmf{E}$; $\mmf{A}^{\circ}$ be the set consisting of
all linear continuous functionals $f$
on $\mmf{E}$ such that
$\sup\limits_{x\in\mmf{A}}{|f(x)|}\leq 1$. For any set
$\mmf{B}$ denote $2^{\mmf{B}}$  the set consisting of all subsets of  $\mmf{B}$.
 For any
function
 $\mmm{S}:[0, \,\infty) \mapsto\mcl{L}(\mmf{E})$ we denote  $\mmm{S}^*$
 the function from $[0, \infty)$ to $\mcl{L}(\mmf{E}^*)$ such that
 $\mmm{S}^*(s)=(\mmm{S}(s))^*$ for each $s\geq 0$.
Recall that the dual space $\mmf{E}^{*}$ is the locally convex space of all linear
continuous functionals on  $\mmf{E}$ with the topology generated by
the semi-norms  of the form
$\|\cdot\|_{\mmf{A}^{\circ}}$, where $\mmf{A}$ is a bounded subset of  $\mmf{E}$.

 Further the expressions $(x,\phi)$ and $\|x\|_{\phi}$, where $x\in\mmf{E}$ and $\phi\in\mmf{E}^*$, is defined by the equalities  $(x,\phi)=\phi(x)$
 and
$\|x\|_{\phi}=|(x,\phi)|$. It is not difficult to show that
 $\|f\|_{\mmf{A}^{\circ}}=\sup\limits_{x\in\mmf{A}}{|(x,f)|}$ for
 any $f\in\mmf{E}^*$.

\begin{definition}\label{deft:1}
Linear operator $\mmm{A}:\mcl{D}(\mmm{A})\to\mmf{E}$ is closed if
its graph $\{(x,\mmm{A}x)|x\in\mcl{D}(\mmm{A})\}$ is closed in the
product space $\mmf{E}\times\mmf{E}$.
\end{definition}
\begin{definition}\label{deft:2}
Linear operator $\mmm{A}:\mcl{D}(\mmm{A})\to\mmf{E}$ is closable if
the closure of its graph $\{(x,\mmm{A}x)|x\in\mcl{D}(\mmm{A})\}$ is
the graph of the linear operator.
\end{definition}

\begin{definition}\label{deft:2:1}
The set $\mmf{A}\subset\mmf{E}$ is sequentially dense in $\mmf{E}$,
if for each $x\in\mmf{E}$ there exists a sequence
$\{x_k\}_{k\in\mmb{N}}$, $x_k\in\mmf{A}$, such that there exists the
limit $\lim_{k\to\infty}{x_k}=x$.
\end{definition}

\begin{definition}\label{deft:2:2}
The set $\mmf{A}\subset\mmf{E}^{*}$ is *-dense in $\mmf{E}^{*}$, if
for each  $\phi\in\mmf{E}^{*}$ there exists a sequence
 $\{\phi_k\}_{k\in\mmb{N}}$, $\phi_k\in\mmf{A}$,
such that for each $x\in\tilde{\mmf{E}}$ there exists the limit
$\lim_{k\to\infty}{(x,\phi_k)}=(x,\phi)$, where $\tilde{\mmf{E}}$ is
the completion of the space $\mmf{E}$.
\end{definition}

\begin{definition}\label{deft:3}
Let $\mmf{K}$ be a non-empty set. A filter base $\mathfrak{F}\subset
2^{K}$ is a system of subsets of  $\mmf{K}$ satisfying the following
condition:

\noindent i) The intersection of any two sets of  $\mathfrak{F}$
contains a set of $\mathfrak{F}.$

\noindent ii) $\mathfrak{F}$ is non-empty and the empty set is not
in $\mathfrak{F}$.
\end{definition}

\begin{definition}\label{deft:4}
Let $\mmm{G}$ be a function from non-empty set $\mmf{K}$ to
$\mmf{E}$. Let $\mathfrak{F}\subset 2^{K}$ be a filter base. The
point $x\in\mmf{X}$ is said to be a limit of the function $\mmm{G}$ with
respect to the filter base $\mathfrak{F}$,
$$\lim\limits_{\mathfrak{F}}{\mmm{G}}=x,$$ if for any $\alpha\in\Omega$
and $\epsilon>0$ there exists
$\mmf{L}_{\alpha,\epsilon}\in\mathfrak{F}$ such that
$\|\mmm{G}(y)-x\|_{\alpha}\leq\epsilon$ for any $y\in
\mmf{L}_{\alpha,\epsilon}$.
\end{definition}

\begin{definition}\label{deft:5}
Let $\mmm{G}$ be a function from non-empty set $\mmf{K}$
to $\mmb{R}$. Let $\mathfrak{F}\subset 2^{K}$ be a filter base. The
point $x\in\mmb{R}$ is said to be the upper limit of the function
 $\mmm{G}$ with respect to the filter base $\mathfrak{F}$
$$\limsup_{\mathfrak{F}}{\mmm{G}}=x$$ if for any $\mmf{L}\in\mathfrak{F}$
and $\epsilon>0$ there exists $\mmf{L}_{\epsilon}\in\mathfrak{F}$
such that $\mmf{L}_{\epsilon}\subset\mmf{L}$ and
$$0\leq\sup\limits_{y\in\mmf{L}_{\epsilon}}{(\mmm{G}(y)-x)}<\epsilon.$$
\end{definition}

\begin{definition}\label{deft:6}
A family of functions
$\mmm{G}_{\beta}\in\mcl{L}(\mmf{E}),\,\,\beta\in\Gamma,$ is
equicontinuous if for any   $\alpha\in\Omega$ there exists
$\{{\alpha}_i\}_{i=1}^{n},\,\,{\alpha}_i\in\Omega,$ and
$\{b_i\}_{i=1}^{n},\,\,b_i>0,$ such that we have
 $$\sup\limits_{\beta\in\Gamma}
 {\|\mmm{G}_{\beta}x\|_{\alpha}}\leq\sum_{i=1}^{n}{b_i\|x\|_{{\alpha}_i}}$$
for each $x\in\mmf{E}$.
\end{definition}

\begin{definition}\label{deft:7}
For any function $\mmm{F}:[0, \infty) \mapsto\mcl{L}(\mmf{E})$ and
 $s\geq 0$ let
$\mmf{B}_{s}^{\mmm{F}}\subset\mcl{L}(\mmf{E})$ be the family of
functions
  $\{\mmm{F}^m(\frac{d}{n})|md/n\leq
s,\,\, n,\,m\in\mmb{N},\,\, d\geq 0\}.$
\end{definition}
\begin{definition}\label{deft:8}
 The set $\EuScript{E}_{\mmf{E}}$ is the set of all functions
$\mathrm{F}:[0,\,\infty)\mapsto\mathcal{L}(\mathbf{E})$ such that
$\mathrm{F}(0)=\mathrm{I}$ and the set $\mmf{B}_{s}^{\mmm{F}}$ is
equicontinuous for each  $s\geq 0$.
\end{definition}
\begin{definition}\label{deft:9} For any $\mmm{F}\in\EuScript{E}_{\mmf{E}},
\,\,\alpha\in\Omega,\,\,s\geq
 0$  let the function
$\|\cdot\|_{\alpha}^{\mmm{F},s}:\mmf{E}\mapsto[0,\,\infty)$ be
defined by the equality
 $\|x\|_{\alpha}^{\mmm{F},s}=\sup\limits_{g\in\mmf{B}_{s}^{\mmm{F}}}
 {\|g(x)\|_{\alpha}},\,\,
 x\in\mmf{E}$.
 \end{definition}

 \begin{definition}\label{deft:10}
  By Definition \ref{deft:9} and induction over n,
  for any $\{\mmm{F}_i\}_{i=1}^{n},\,\,\mmm{F}_i\in\EuScript{E}_{\mmf{E}},\,\,\alpha\in\Omega,\,\,\{s_i\}_{i=1}^{n},\,\,s_i\geq
 0,$ let the function
$$\|\cdot\|_{\alpha}^{\mmm{F}_1,s_1;\mmm{F}_2,s_2;\ldots;
\mmm{F}_n,s_n}:\mmf{E}\mapsto[0,\,\infty)$$  be defined by the
equality
 $$\|x\|_{\alpha}^{\mmm{F}_1,s_1;\mmm{F}_2,s_2;
 \ldots;\mmm{F}_n,s_n}=\sup\limits_{g\in\mmf{B}_{s_n}^{\mmm{F}_n}}
 {{\|g(x)\|_{\alpha}^{\mmm{F}_1,s_1;\mmm{F}_2,s_2;\ldots;\mmm{F}_{n-1},s_{n-1}}}},\,
 x\in\mmf{E},$$ for each natural $n\geq 2$.
 \end{definition}
 \begin{remark}\label{r:1}
 In the case of  $\mathrm{F}\in\EuScript{E}_{\mmf{E}}$, the set of
 semi-norms
 $\{\|\cdot\|_{\alpha}^{\mmm{F},s}|\alpha\in\Omega,\,s\geq 0\}$
 generates the topology of the space $\mmf{E}$. It follows from the definition of equicontinuity.
In the same way if
 $\mmm{F}_1,\,\mmm{F}_2,\ldots,\mmm{F}_n\in\EuScript{E}_{\mmf{E}}$,
 then the set of semi-norms
 $\{\|\cdot\|_{\alpha}^{\mmm{F}_1,s_1;\mmm{F}_2,s_2;\ldots;\mmm{F}_n,s_n}|\alpha\in\Omega,\,s_1,s_2,\ldots ,s_n\geq 0\}$
 generates the topology of the space $\mmf{E}$.
 \end{remark}

 \begin{definition}\label{deft:11} For any function $\mmm{F}:[0,\, \infty) \mapsto\mcl{L}(\mmf{E})$
let $\mmf{E}_{\mmm{F}}$ be the set consisting of all $f\in\mmf{E}$
such that there exists a sequence $\{f_s\}_{s\geq 0}$ that there
exist the limits
  $\lim\limits_{h\to
0}{f_h}=f$ and $\lim\limits_{h\to
0}{h^{-1}(\mmm{F}(h)-\mmm{F}(0))f_h}$.
\end{definition}

\begin{definition}\label{deft:12} By the (strong) many-valued effective derivative at the point $0$
of a function $\mmm{F}:[0,\, \infty) \mapsto\mcl{L}(\mmf{E})$ we mean the
linear map $\mmm{F}'_{mef}(0):\mmf{E}_{\mmm{F}}\mapsto 2^{\mmf{E}}$
such that for any $f\in\mmf{E}_{\mmm{F}}$ $\mmm{F}'_{mef}(0)f$
consists of all $g\in\mmf{E}$ for which there exists a sequence
 $\{f_s\}_{s\geq 0}$ such that there exist the limits
  $\lim\limits_{h\to
0}{f_h}=f$ and  $\lim\limits_{h\to
0}{h^{-1}(\mmm{F}(h)-\mmm{F}(0))f_h}=g$.
\end{definition}
\begin{definition}\label{deft:13} We say that $\mmm{F}'_{mef}(0)$
is single-valued if $\mmm{F}'_{mef}(0)f$ consists of one element for
each  $f\in\mmf{E}_{\mmm{F}}$.
\end{definition}
\begin{definition}\label{deft:14}
 Let $\mmm{F}$ be a function from  $[0,\, \infty)$ to $\mcl{L}(\mmf{E})$ and $\mmm{F}'_{mef}(0)$ be single-valued.
  We call the (strong) effective derivative at the point  $0$
of the function $\mmm{F}$ the linear map
$\mmm{F}'_{ef}(0):\mmf{E}_{\mmm{F}}\mapsto {\mmf{E}}$ such that for
any $f\in\mmf{E}_{\mmm{F}}$
$\mmm{F}'_{ef}(0)f\in\mmm{F}'_{mef}(0)f.$ So, we say that there
exists  $\mmm{F}'_{ef}(0)$ if $\mmm{F}'_{mef}(0)$ is single-valued.
\end{definition}

\begin{definition}\label{deft:15} We call the (strong)  derivative at the point $0$ of a function
$\mmm{F}:[0,\, \infty) \mapsto\mcl{L}(\mmf{E})$ the linear map
$\mmm{F}'(0):\mcl{D}(\mmm{F}'(0))\mapsto\mmf{E}$ defined by the
equality  $\mmm{F}'(0)\psi=\lim\limits_{h\to
0}{h^{-1}(\mmm{F}(h)\psi-\mmm{F}(0)\psi)},
\,\psi\in\mcl{D}(\mmm{F}'(0)),$ where $\mcl{D}(\mmm{F}'(0))$ is the
space of all $\psi\in\mmf{E}$ for which the above limit exists.
\end{definition}

\begin{definition}\label{deft:16}
A function $\mmm{T}:[0,\,\infty)\to\mcl{L}(\mmf{E})$ is called
a $lC_0$-semigroup if the following conditions hold:

\noindent 1) $\mmm{T}(0)=\mmm{I}$,
$\mmm{T}(l+m)=\mmm{T}(l)\mmm{T}(m)$ for each $l,\,m\in[0,\,\infty)$.

\noindent 2) The function $\mmm{T}$ is continuous.

\noindent 3) $\mmm{T}\in\EuScript{E}_{\mmf{E}}.$
\end{definition}

\begin{definition}\label{deft:17}
A linear operator $\mmm{Z}$ is called the generator of
the $lC_0$-semigroup $\mmm{T}$ if $\mmm{Z}$ is the (strong) derivative
at the point $0$ of the function $\mmm{T}$.
\end{definition}
It is well-known fact that there exists one-to-one correspondence
between $lC_0$-semigroups and its generators.

\begin{definition}\label{deft:18}
A set $\mcl{D}$ is called a core of a generator $\mmm{Z}$ of
the $lC_0$-semigroup $\mmm{T}$ if the closure of the subset
$\{(x,\mmm{Z}x)|x\in\mcl{D}\}$ of $\mmf{E}\times\mmf{E}$ coincides
with the graph of the generator $\mmm{Z}$.
\end{definition}

\begin{definition}\label{deft:19}
 The set $\EuScript{F}_{\mmf{E}}$ is the set of all functions
$\mathrm{F}:[0,\,\infty)\mapsto\mathcal{L}(\mathbf{E})$ such that
$\mathrm{F}\in\EuScript{E}_{\mmf{E}}$ and
$\mcl{D}(\mmm{F}'_{mef}(0))$ is dense in $\mmf{E}$.
\end{definition}

\begin{definition}\label{deft:20} The set $\mcl{Z}_{\mmf{E}}$ is the set
 of all densely defined linear operators $\mmm{Z}$ in  $\mathbf{E}$
 for which
 there exist
 $\mathrm{F}\in\EuScript{F}_{\mmf{E}}$ such that
$\mcl{D}(\mmm{Z})\subset\mcl{D}(\mmm{F}'_{mef}(0))$ and
$\mmm{Z}f\in\mmm{F}'_{mef}(0)f$ for each $f\in\mcl{D}(\mmm{Z})$.
\end{definition}

We show in Theorem \ref{cor:1} that the set $\mcl{Z}_{\mmf{E}}$
consists of closable operators. Moreover, if $\mmf{E}$ is a
Hilbert space, then, by Proposition \ref{prop:11}, the set of all dense
defined dissipative operators is contained in $\mcl{Z}_{\mmf{E}}$.

In the article we assume that  $\mmf{X}$ is a Hausdorff locally
convex space on the field $\mmb{T}\in\{\mmb{R},\mmb{C}\}$)
 with the topology generated by semi-norms
$\|\cdot\|_{\alpha},\,\,\alpha\in\Omega,$ and $\mcl{B}$ is either
a separable reflexive Banach space or Hilbert space.

\begin{remark}\label{remt:1} The Hausdorff  condition on $\mmf{X}$ is imposed
for simplicity. Actually this condition is not essential because
instead of $\mmf{X}$  we can consider its appropriate
factor space.
\end{remark}

\section{Existence of effective derivative $\mmm{S}'_{ef}(0)$ for $\mmm{S}\in\EuScript{F}_{\mmf{X}}$.}\label{1.2}
\begin{lemma}\label{t:5}
Let $\mathrm{F}\in\EuScript{F}_{\mmf{X}}$. Then there exists
$$\lim\limits_{t i/k \to
0}{({\mathrm{F}}^{i}(t/k)-\mathrm{I})g}=0,i,\,\, k \in
\mathbb{N},\,\,t> 0,$$ for any $g \in \mmf{X}$.
\end{lemma}

\begin{remark} The limit in Lemma~\ref{t:5} is considered as a limit
w. r. t. the filter base consisting of the sets of the form
$\{(i,k,t)|\,t\frac{i}{k}< \epsilon,\,\,i,\,\, k \in
\mathbb{N},\,\,t> 0 \}$, where $\epsilon>0$.
\end{remark}

\begin{proof}[Proof of Lemma \ref{t:5}] Fix $\epsilon_0>0$.
If $g \in \mcl{D}(\mmm{F}_{mef}'(0))$, then there exists a sequence
$\{g_s\}_{s\geq 0}$ such that there exist the limits
$\lim\limits_{s\to 0}{g_s}=g$, $\lim\limits_{s\to
0}{s^{-1}(\mmm{F}(s)-\mmm{F}(0))g_s}$. Put $h=\lim\limits_{s\to
0}{s^{-1}(\mmm{F}(s)-\mmm{F}(0))g_s}$. Then for any
$\alpha\in\Omega$ and $\epsilon\in(0,\,\epsilon_0)$ the following
chain of the inequalities is satisfied:
\begin{eqnarray}
\nonumber &
&{\textstyle\|({\mathrm{F}}^{i}(\frac{t}{k})-\mathrm{I})g\|_{\alpha}}\leq
{\textstyle{i
\frac{t}{k}}\bigl\|\bigl(\frac{\mathrm{F}(\frac{t}{k})^{i-1}
+\cdots+\mathrm{I}}{i}\bigr)
\frac{(\mathrm{F}(\frac{t}{k})g_{t/k}-g_{t/k})}{\frac{t}{k}}\bigr\|_{\alpha}}\\\nonumber
&+&\|{({\mathrm{F}}^{i}({\textstyle\frac{t}{k}})-\mathrm{I})(g-g_{t/k})}\|_{\alpha}
 \leq {\textstyle i
 \frac{t}{k}\bigl\|\frac{(\mathrm{F}(\frac{t}{k})g_{t/k}-g_{t/k})}
 {\frac{t}{k}}\bigr\|_{\alpha}^{\mmm{F},\epsilon_0}}\\&+&{\textstyle  \|{({\mathrm{F}}^{i}
 ({\textstyle\frac{t}{k}})-\mathrm{I})(g-g_{t/k})}\|_{\alpha}\leq
2\epsilon   \| h
 \|_{\alpha}^{\mmm{F},\epsilon_0}}
 +2\|{g-g_{t/k}}\|^{\mmm{F},\epsilon_0}_{\alpha}\nonumber
\end{eqnarray}
for $ti/k<\epsilon$. Therefore, we have $\lim\limits_{ ti/k \to
0}{({\mathrm{F}}^{i}({\textstyle\frac{t}{k}})-\mathrm{I})g}=0$.

If $g \notin \mcl{D}(\mmm{F}_{mef}'(0)),$ then for any $\varepsilon
\in(0,\,\epsilon_0)$ and $\alpha\in\Omega$ there exists
$g'\in\mcl{D}(\mmm{F}_{mef}'(0))$ such that
$\|g-g'\|_{\alpha}^{\mmm{F},\epsilon_0}<\varepsilon /2 $ and

$$
\begin{array}{c}
\limsup\limits_{ ti/k \to
0}{\|({\mathrm{F}}^{i}({\textstyle\frac{t}{k}})-\mathrm{I})g}\|_{\alpha}\leq
\limsup\limits_{ ti/k \to
0}{\|({\mathrm{F}}^{i}({\textstyle\frac{t}{k}})-\mathrm{I})(g-g')}\|_{\alpha}\\+
\limsup\limits_{ ti/k \to
0}{\|({\mathrm{F}}^{i}({\textstyle\frac{t}{k}})-\mathrm{I})g'}\|_{\alpha}\leq\varepsilon.
\end{array}
$$
Since $\varepsilon$ and $\alpha$ are arbitrary the Lemma is proved.
\end{proof}

\begin{lemma}\label{t:7} Let $\mathrm{F}\in\EuScript{F}_{\mmf{X}}$ and $s>0$.
Then for any $g \in \mmf{X}$ there exists
$$\lim\limits_{ \mathfrak{F}_s}{({\mathrm{F}}^{i}(t/k)-{\mathrm{F}}^{l}(t/k))g}=0, $$
where the filter base $\mathfrak{F}_s$ consists of the sets of the
form $\{(i,l,k,t)|\,0<|t\frac{i-l}{k}|<\epsilon,\,\,
t\frac{i}{k}\leq s,\,\,t\frac{l}{k}\leq s,\,i,\, k,\, l \in
\mathbb{N},\,t>0 \}$, where $\epsilon>0$.
\end{lemma}
\begin{proof}[Proof of Lemma \ref{t:7}]
Fix $s>0$ and $\alpha\in\Omega$. By Lemma~\ref{t:5},
\begin{eqnarray}\limsup\limits_{\mathfrak{F}_s}{\|
({\mathrm{F}}^{i}(t/k)-{\mathrm{F}}^{l}(t/k))g\|_{\alpha}}\nonumber\\\leq
\lim\limits_{\mathfrak{F}_s}{\|({\mathrm{F}}^{|i-l|}(t/k)-\mathrm{I})g\|_{\alpha}^{\mmm{F},s}}=
\lim\limits_{t i/k \to
0}{\|({\mathrm{F}}^{i}(t/k)-\mathrm{I})g\|_{\alpha}^{\mmm{F},s}}=0\nonumber
\end{eqnarray}
  \end{proof}

\begin{lemma}\label{t:8} Let $\mathrm{F}\in\EuScript{F}_{\mmf{X}}$ and $s>0$
 be fixed.
 Then for any sequence $\{g_s\}_{s\geq 0},\,\,g_s\in\mmf{X},$  such
 that there exist the limits
 $\lim\limits_{s\to 0}{g_s}=g$ and
$\lim\limits_{s\to
0}{s^{-1}(\mmm{F}(s)-\mmm{F}(0))g_s}=h,\,\,g,\,\,h\in\mmf{X},$ we
have
$$\lim\limits_{  \mathfrak{F}_s}
{\frac{({\mathrm{F}}^{i}(t/k)-{\mathrm{F}}^{l}(t/k))g_{t/k}}{(i-l)t/k}-
({\mathrm{F}}^{\min(i,l)}(t/k))h}=0,$$ where the filter base
$\mathfrak{F}_s$ consists of the sets of the form
$\{(i,l,k,t)|\,0<|t\frac{i-l}{k}|<\epsilon,\,\, t\frac{i}{k}\leq
s,\,\,t\frac{l}{k}\leq s,\,i,\, k,\, l \in \mathbb{N},\,t>0 \}$,
where $\epsilon>0$.
\end{lemma}

\begin{proof}[Proof of Lemma \ref{t:8}]  Fix $s>0$.
 For any $\alpha\in\Omega$ choose
${\alpha}_1,\,\alpha_2,...,\alpha_n\in\Omega$
 and ${\beta}_1,\,\beta_2,...,\beta_n>0$ such that
 $\|x\|_{\alpha}^{\mmm{F},s}\leq\beta_1\|x\|_{\alpha_1}+
 \beta_2\|x\|_{\alpha_2}+
 \cdots+ \beta_n\|x\|_{\alpha_n}$, for each $x\in\mmf{X}$.
 Then we have the following chain of the inequalities:

\begin{eqnarray}
& &\nonumber{\textstyle\limsup\limits_{
\mathfrak{F}_s}{{\|\frac{({\mathrm{F}}^{i}(\frac{t}{k})-
{\mathrm{F}}^{l}(\frac{t}{k}))g_{t/k}}{(i-l)\frac{t}{k}}-
({\mathrm{F}}^{\min(i,l)}({\textstyle\frac{t}{k}}))h}}\|_{\alpha}}
\\&\leq&\nonumber{\textstyle
\limsup\limits_{
\mathfrak{F}_s}{{\|\frac{({\mathrm{F}}^{|i-l|}(\frac{t}{k})
-\mathrm{I})g_{t/k}}{|i-l|{\frac{t}{k}}}-h\|_{\alpha}^{\mmm{F},s}}}}
\\
&\leq& \nonumber{\textstyle\limsup\limits_{ ti/k \to
0}{{\|\frac{({\mathrm{F}}^{i}(\frac{t}{k})-\mathrm{I})g_{t/k}}{i
{\frac{t}{k}}}-h\|_{\alpha}^{\mmm{F},s}}}}
\\ &\leq& \nonumber{\textstyle\limsup\limits_{ ti/k \to
0}{{{\bigl\|\bigl(\frac{\mathrm{F}(\frac{t}{k})^{i-1}
+\cdots+\mathrm{I}}{i}\bigr)
(\frac{(\mathrm{F}(\frac{t}{k})g_{t/k}-g_{t/k})}{\frac{t}{k}}-h)}\bigr\|
_{\alpha}^{\mmm{F},s}}}}
\\&+& \nonumber{\textstyle\limsup\limits_{t i/k \to
0}{{{\bigl\|\bigl(\frac{\mathrm{F}(\frac{t}{k})^{i-1}
+\cdots+\mathrm{I}}{i}\bigr)h -h}\bigr\|_{\alpha}^{\mmm{F},s}}}}
\\&\leq& \nonumber
 {\textstyle\limsup\limits_{ ti/k \to 0}{{\sum_{l=1}^{n}{{\beta_l\bigl\|
\frac{(\mathrm{F}(\frac{t}{k})g_{t/k}-g_{t/k})}{\frac{t}{k}}-h}\bigr\|
_{\alpha_l}^{\mmm{F},s}}}}}
\\&+&
\nonumber{\textstyle\limsup\limits_{ ti/k \to
0}{{{\bigl\|\bigl(\frac{(\mathrm{F}(\frac{t}{k})^{i-1}-\mathrm{I})
+\cdots+(\mathrm{I}-\mathrm{I})}{i}\bigr)h}\bigr\|_{\alpha}^{\mmm{F},s}}}=0,}
\end{eqnarray}
where the last limit is equal to $0$ by Lemma~\ref{t:5}.\end{proof}

\begin{proposition}\label{t:1}
Assume that  $\mathrm{F}\in\EuScript{F}_{\mmf{X}}$, $t>0$ and
$\{n_k\}_{k=1}^{\infty}$ is a strictly increasing sequence of natural
numbers. Then for any separable closed linear space
$\Phi\subset\mmf{X}$ and $\phi\in\mmf{X}^*$ there exists a
subsequence
 $\{g_k\}_{k=1}^{\infty}$ of the sequence
$\{n_k\}_{k=1}^{\infty}$ such that there exists
\begin{eqnarray}
 \lim\limits_{k
\rightarrow\infty}{(\{\mathrm{F}({\textstyle\frac{t}{g_k}})\}^{[{g_k}s]}
g,\phi)}\label{eqn:l1}
\end{eqnarray}
for any $g\in \Phi$ and $s>0$. Furthermore, if such subsequence
 $\{g_k\}_{k=1}^{\infty}$ is chosen, then the family of the
 functions
 $\mathrm{T}_s: \Phi\times\{\phi\}\mapsto\mathbb{T}$,
 $s\geq 0$, defined by
\begin{equation}
\mathrm{T}_s (g,\,\phi)=\lim\limits_{k
\rightarrow\infty}{(\{\mathrm{F}({\textstyle\frac{t}{g_k}})\}^{[{g_k}\frac{s}{t}]}
g,\phi)},\label{eqn:l2}
\end{equation}
  satisfies the following conditions:

\noindent a) $\mathrm{T}_s (g, \phi)$ is linear w. r. t. $g$;

\noindent b) $\mathrm{T}_s (g, \phi)$ is continuous w. r. t. $s$ for
any $g\in\Phi$;

\noindent c) If $f\in\Phi\cap\mcl{D}(\mmm{F}_{mef}'(0))$ and
$h\in\Phi\cap\mmm{F}_{mef}'(0)f,$ then there exists
$$(\mathrm{T}_{s} (f, \phi))^{'}_{s}=\mathrm{T}_{s}(h,
\phi)$$ for each $s\geq 0$.

\end{proposition}

\begin{proof}[Proof of Proposition \ref{t:1}]
Choose a subsequence $\{g_k\}_{k=1}^{\infty}$ of the sequence
$\{n_k\}_{k=1}^{\infty}$ such that for any nonnegative rational
 $r$ there exists
$$\lim\limits_{k\to\infty}{(\{\mathrm{F}({\textstyle\frac{t}{g_k}})\}^{[{g_k}r]} g,\phi)},
\,\,g\in\mmf{C},$$ where $\mmf{C}$ is a dense countable subset of
$\Phi$. Indeed, it is possible to find such subsequence because the
family $\mmf{B}^{\mmm{F}}_s$ is equicontinuous for each $s>0$ and we
have a countable number of conditions imposed on sequence
$\{g_k\}_{k=1}^{\infty}$. Hence, it follows from the density of
$\mmf{C}$ in $\Phi$ that there exists
$$\lim\limits_{k\to\infty}{(\{\mathrm{F}({\textstyle\frac{t}{g_k}})\}^{[{g_k}r]} g,\phi)}$$
for any $g\in\Phi$ and $r\in[0,\,\infty)\cap\mmb{Q}$. Let us show
that the last limit exists for any real number $r\geq 0$. Fix $r>0$
and $\epsilon
> 0$. Put $\epsilon^{'}=\epsilon/3$. By Lemma~\ref{t:7} we can
choose $\epsilon^{''}(g,\epsilon^{'})$ such that the inequality
$$|(({\mathrm{F}}^{i}(t/k)-{\mathrm{F}}^{l}(t/k))g,\phi)|<\epsilon^{'}$$
follows from the inequalities $|(i-l)/k|<\epsilon^{''}$ and
$|t(i+l)/k|<2r$, where $i,\,\,l,\,\,k\in\mmb{N}$. Choose positive
$s\in\mmb{Q}$ such that $|r-s|<\epsilon^{''}/2$. There exists
$k_0\in\mmb{N}$ such that $2/g_{k_0}<\epsilon^{''}/2$. Thus, it
follows that
$$|[g_n r]-[g_n s]|/g_n \leq (|g_nr-g_n s|+2)/g_n
\leq\epsilon^{''}$$ for each natural $n>k_0$. Now it follows from
the existence of the limit
$$\lim\limits_{k\to\infty}\{(\mathrm{F}
({\textstyle\frac{t}{g_k}})\}^{[{g_k}s]}g,\phi)$$
that we can find  $n_0>k_0$
 such that for any natural $l,m>
n_0$ we have the following inequality
$$|((\{\mathrm{F}(t/{g_l})\}^{[{g_l}s]}-
\{\mathrm{F}(t/{g_m})\}^{[{g_m}s]})g, \phi)|<\epsilon/3.$$
Therefore, we get the following chain of the inequalities:

\begin{eqnarray}
 &\,&|((\{\mathrm{F}(t/{g_l})\}^{[{g_l}r]}-
\{\mathrm{F}(t/{g_m})\}^{[{g_m}r]})g, \phi)|\leq
|((\{\mathrm{F}(t/{g_l})\}^{[{g_l}r]}\nonumber\\&-&
\{\mathrm{F}(t/{g_l})\}^{[{g_l}s]})g, \phi)|
 +|((\{\mathrm{F}(t/{g_l})\}^{[{g_l}s]}-
\{\mathrm{F}(t/{g_m})\}^{[{g_m}s]})g, \phi)|\nonumber\\&+&
|((\{\mathrm{F}(t/{g_m})\}^{[{g_m}s]}-
\{\mathrm{F}(t/{g_m})\}^{[{g_m}r]})g, \phi)|
 \leq \epsilon/3 +\epsilon/3+\epsilon/3\leq\epsilon.\nonumber
\end{eqnarray}
As a consequence, we see that there exists
$$\lim\limits_{k\to\infty}{(\{\mathrm{F}({\textstyle\frac{t}{g_k}})\}^{[{g_k}s]}
g,\phi)},\,\, s\geq 0,\,\, g\in\Phi.$$ Thus, the first part of
Proposition~\ref{t:1} is proved. Now let us consider parts (a) --
(c) of Proposition~\ref{t:1}.

\begin{trivlist}
\item[(a)] Part $(a)$ immediately follows from the definition of
$\mmm{T}_s$. \item[(b)] Fix $g\in\Phi$. It follows from
Lemma~\ref{t:7} that for any  $\epsilon>0$ there exist $v_0>0$,
$k_0\in\mmb{N}$ such that
$$|(\{\mathrm{F}({\textstyle\frac{t}{g_k}})\}^{[{g_k}\frac{s+v}{t}]} g,\phi)-
(\{\mathrm{F}({\textstyle\frac{t}{g_k}})\}^{[{g_k}\frac{s}{t}]}
g,\phi)|<\epsilon,\,\,k\in\mmb{N},$$ for any
$v\in[-s,\,\infty)\cap(-v_0,\,v_0)$ and $k>k_0$. Tending $k$ to
infinity, we get the inequality
$$|\mmm{T}_{s+v}(g,\phi)-\mmm{T}_{s}(g,\phi)|<\epsilon.$$
From the arbitrariness of $\epsilon$ we easily infer part  $(b)$ of
Proposition~\ref{t:1}.

\item[(c)] Assume that $f\in\Phi\cap\mcl{D}(\mmm{F}_{mef}'(0))$
and $h\in\Phi\cap\mmm{F}_{mef}'(0)f$. Let $\{f_s\}_{s\geq 0}$ be a
sequence such that $\lim\limits_{s\to 0}{f_s}=f$ and
$\lim\limits_{s\to 0}{s^{-1}(\mmm{F}(s)-\mmm{F}(0))f_s}=h$. Fix arbitrary $s>0$. It
follows from Lemma~\ref{t:8} that for any  $\epsilon>0$ there exist
$v_0>0$, $k_0\in\mmb{N}$ such that
\begin{eqnarray}
\nonumber& &
|({\mathrm{F}}^{[{g_k}\frac{s}{t}]}({\textstyle\frac{t}{g_k}})f_{{\textstyle\frac{t}{g_k}}},\phi)-
 ({\mathrm{F}}^{
[{g_k}\frac{s+v}{t}]}({\textstyle\frac{t}{g_k}})f_{{\textstyle\frac{t}{g_k}}},\phi)-{\textstyle\frac{t}{g_k}}([{g_k}{\textstyle\frac{s}{t}}]-
[{g_k}{\textstyle\frac{s+v}{t}}])\\ &\times&
({\mathrm{F}}^{\min([{g_k}\frac{s}{t}],[{g_k}\frac{s+v}{t}])}({\textstyle\frac{t}{g_k}})h,\phi)|\nonumber
<\epsilon
|[{g_k}{\textstyle\frac{s}{t}}]-[{g_k}{\textstyle\frac{s+v}{t}}]|{\textstyle\frac{t}{g_k}},\,\,k\in\mmb{N},
\end{eqnarray}
for any $v\in[-s,\,\infty)\cap(-v_0,\,v_0)$ and $k>k_0$. Tending $k$
to infinity, by the equicontinuity of $\mmf{B}_{p}^{\mmm{F}}$,
$p\geq 0$,
 we get the inequality
\begin{equation}
|\mathrm{T}_{s}(f,\phi)-\mathrm{T}_{s+v}(f,\phi)-v\mathrm{T}_{\min{(s,
s+v)}}(h,\phi)|\leq \epsilon |v|. \nonumber
\end{equation}
From the arbitrariness of $\epsilon$ we deduce part $(c)$ of
Proposition~\ref{t:1}.
\end{trivlist}
  \end{proof}

\begin{proposition}\label{txt} Assume that
$\mathrm{F}\in\EuScript{F}_{\mmf{X}}$. Then the function
$\mmm{F}_{mef}'(0)$ is single-valued.
\end{proposition}
\begin{proof}[Proof of Proposition \ref{txt}] Fix $f\in\mcl{D}(\mmm{F}_{mef}'(0))$. Assume that
 for any $i\in\{0,1\}$
there exists the sequence
 $\{f^i_s\}_{s\geq 0}$
such that there exist  $\lim\limits_{s\to 0}{f^i_s}=f$ and
$\lim\limits_{s\to
0}{s^{-1}(\mmm{F}(s)-\mmm{I})f^i_s}=g_i,\,\,g_i\in\mmf{X}.$ It is
enough to show that $g_0=g_1$. We will argue by contradiction.
Assume that $g_0\neq g_1$. Then there exists $\phi\in\mmf{X}^{*}$
such that $(g_0,\,\phi)\neq (g_1,\,\phi)$. Let $\mmf{B}$ be
 the closure of the linear span of elements  $f,\,\,g_0,\,\,g_1.$
According to Proposition~\ref{t:1}, there exists a sequence
$\{g_k\}_{k=1}^{\infty}$ such that we can define the family of the
function $\mathrm{T}_s: \Phi\times\{\phi\}\mapsto\mathbb{T}$,
 $s\geq 0$, by formula \eqref{eqn:l2} with $\Phi=\mmf{B}$.
Thus, by part  (c) of Proposition~\ref{t:1}, we infer that
 $$(\mathrm{T}_{s} (f, \phi))^{'}_{s=0}=(g_0,
\phi)=(g_1, \phi).$$ Hence, we get contradiction with the assumption
$g_0\neq g_1$. \end{proof}

\begin{corollary} Assume that
$\mathrm{F}\in\EuScript{F}_{\mmf{X}}$. Then there exists the
(strong) effective derivative $\mmm{F}_{ef}'(0)$.
\end{corollary}

\begin{proposition}\label{n:1}Assume that
$\mathrm{F}\in\EuScript{F}_{\mmf{X}}$. Then the operator
$\mmm{F}_{ef}'(0)$ is closable.
\end{proposition}

\begin{proof}[Proof of Proposition \ref{n:1}] Assume that
 for any $i\in\{0,1\}$ and
 $\phi\in\mmf{X}^*$ there exists the sequence
  $\{f^{i,\alpha}_k\}_{k\in\mmb{N}},\,\,f^{i,\alpha}_k\in\mcl{D}(\mmm{F}_{ef}'(0)),$
   such that there exist
    $\lim\limits_{k\to \infty}{\|f^{i,\alpha}_k-f\|_{\alpha}}=0$ and
$\lim\limits_{k\to
\infty}{\|\mmm{F}_{ef}'(0)f^{i,\alpha}_k-g_i\|_{\alpha}}=0\,\,\mbox{for
some fixed}\,\,f,\, g_i\in\mmf{X}.$ It is enough to show that
$g_0=g_1$. We will argue by contradiction. Assume that $g_0\neq
g_1$. Then there exists $\phi\in\mmf{X}^{*}$ such that
$(g_0,\,\phi)\neq (g_1,\,\phi)$. Put
$f_k=f^{1,\alpha}_k-f^{0,\alpha}_k$, $k\in\mmb{N}$, and $h=g_1-g_0$.
Then $(h,\,\phi)\neq 0$. Let $\mmf{B}$ be
 the closure of the linear span of elements $f_i$ and $\mmm{F}_{ef}'(0)f_i$,
 where
$i\in\mmb{N}$.
 According to Proposition~\ref{t:1}, there exists a sequence
  $\{g_k\}_{k=1}^{\infty}$
 such that we can define the family of the
function $\mathrm{T}_s: \Phi\times\{\phi\}\mapsto\mathbb{T}$,
 $s\geq 0$, by formula \eqref{eqn:l2} with $\Phi=\mmf{B}$.
Thus, by part  (c) of Proposition~\ref{t:1}, we infer that
$$(\mathrm{T}_{s}
(f_i,\phi))^{'}_{s}=\mathrm{T}_{s}(\mathrm{F}_{ef}'(0)f_i,\phi)$$
for any $i\in\mmb{N}$. Therefore,
\begin{equation}
\mathrm{T}_{s}(f_i,\phi)-(f_i,\phi)=\int_0^{s}{\mathrm{T}_{t}(\mathrm{F}_{ef}'(0)f_i,\phi)\,
d t},i\in\mmb{N}.\label{eqn:eq-2}
\end{equation}
Putting $i\to\infty$  in expression \eqref{eqn:eq-2}, we get the
equality
$$
\mathrm{T}_{s}(0,\phi)-(0,\phi)=\int_0^{s}{\mathrm{T}_{s}(h,\phi)\,
d s}=0,\,\,s\geq 0.
$$
Indeed, it immediately follows from \eqref{eqn:eq-2} and the
equicontinuity of the family of the function
$\{\mathrm{T}_{l}(\cdot,\phi)|\, l\in [0,\, s]\}$. Hence, we get
\begin{equation}
(\int_0^{s}{\mathrm{T}_{t}(h,\phi)\,
dt})'_s=\mathrm{T}_{s}(h,\phi)=0,\,\,s\geq 0.\label{eqn:eq-3}
\end{equation}
Put $s=0$ in \eqref{eqn:eq-3}. Consequently, we have
$$
\mathrm{T}_{s}(h,\phi)=(h,\phi)=0.
$$
Thus, we get contradiction with the assumption $g_0\neq g_1$.
\end{proof}

As a consequence of Proposition~\ref{n:1}, we get
\begin{theorem}\label{cor:1}Let
$\mathrm{F}\in\EuScript{F}_{\mmf{X}}$. Let also a linear operator
$\mmm{Z}$ have the domain
$\mcl{D}(\mmm{Z})\subset\mcl{D}(\mmm{F}'_{ef}(0))$ and
$$\mmm{Z}f=\mmm{F}'_{ef}(0)f,\,\,f\in\mcl{D}(\mmm{Z}).$$ Then
the operator $\mmm{Z}$ possesses a closure.
\end{theorem}
From Theorem~\ref{cor:1} it  follows that any operator
$\mmm{Z}\in\mcl{Z}_{\mmf{X}}$ is closable. Below we will often use
the fact that an operator $\mmm{Z}\in\mcl{Z}_{\mmf{X}}$ possesses a
closure.

It turns out that if $\mmf{X}$ is a Hilbert space then any dissipative densely defined operator is contained in  $\mcl{Z}_{\mmf{X}}$.

\begin{proposition}\label{prop:11}
Let $\mmb{H}$ a Hilbert space. Then $\mcl{Z}_{\mmb{H}}$ contains all
 dissipative densely defined operators $\mmm{Z}$ in $\mmb{H}$.
\end{proposition}

\begin{proof}[ Proof of Proposition~\ref{prop:11}]
See Proposition 4.2 in \cite{Nekl2}.
\end{proof}

\section{Uniqueness of solution for the abstract Cauchy problem.}\label{1.3}
\begin{lemma}\label{lem:lem-4}
Suppose that the assumptions of Proposition~\ref{t:1} are satisfied. Assume also
that $r>0$,
$f_i\in\mcl{D}(\mmm{F}_{ef}'(0))\cap\Phi$,
$\mmm{F}_{ef}'(0)f_i\in\Phi$ for any $i\in\mmb{N}$. Furthermore, suppose that
$$
\lim\limits_{i\to\infty}{\|f_i-h\|_{\phi}^{\mmm{F},r}}=0,\,\,
\lim\limits_{i\to\infty}{\|\mmm{F}_{ef}'(0)f_i-g\|_{\phi}^{\mmm{F},r}}=0.
$$
Then the inequality
$$(\mathrm{T}_{s}(h,\phi))^{'}_{s}=\mathrm{T}_{s}(g,\phi)$$
holds for any $s\in[0,\,r).$
\end{lemma}
\begin{proof}[Proof of Lemma~\ref{lem:lem-4}]
Let $\mmm{Z}=\mmm{F}_{ef}'(0)$. Note that in the equalities
\begin{eqnarray}
&&\lim\limits_{l\to
0}{{{l}^{-1}(\mmm{T}_{s+l}(f_i,\phi)-\mmm{T}_{s}(f_i,\phi))}}=\lim\limits_{l\to
0}{{{l}^{-1}{\int_{s}^{s+l}{\mmm{T}_{m}(\mmm{Z}f_i,\phi)\, d
m}}}}\nonumber\\&=&\mmm{T}_{s}(\mmm{Z}f_i,\phi),\,\,s\in[0,\,r),\nonumber
\end{eqnarray}
we have the uniform convergence w.r.t. $i\in\mmb{N}$, because there exists
$$
\lim\limits_{i\to\infty}{\mmm{T}_{s}(\mmm{Z}f_i,\phi)}=\mmm{T}_{s}(g,\phi),\,\,s\in[0,\,r),
$$ and the family of the functions
$\{\mathrm{T}_{\cdot}(\mmm{Z}f_i,\phi)|\, i\in\mmb{N}\}$
is equicontinuous on $[0,\,l]$. Therefore, we get

\begin{eqnarray}
\mmm{T}_s (g,\phi)&=&\lim\limits_{i\to\infty}{\lim\limits_{l\to
0}{{{l}^{-1}(\mmm{T}_{s+l}(f_i,\phi)-\mmm{T}_{s}(f_i,\phi))}}}\nonumber\\&=&\lim\limits_{l\to
0}{\lim\limits_{i\to\infty}{{{l}^{-1}(\mmm{T}_{s+l}(f_i,\phi)-\mmm{T}_{s}(f_i,\phi))}}}\nonumber\\
&=&\lim\limits_{l\to
0}{{{l}^{-1}(\mmm{T}_{s+l}(h,\phi)-\mmm{T}_{s}(h,\phi))}}=(\mathrm{T}_{s}
(h,\phi))^{'}_{s},\,\,s\in[0,\,r).\nonumber
\end{eqnarray}
  \end{proof}
\begin{lemma}\label{n:2}Suppose that the assumptions of Proposition~\ref{t:1} are satisfied.
If $f\in\mcl{D}(\overline{\mmm{F}_{ef}'(0)})\cap\Phi$ and
$\overline{\mmm{F}_{ef}'(0)}f\in\Phi$, then
$$
(\mathrm{T}_{s}
(f,\phi))^{'}_{s}=\mathrm{T}_{s}(\overline{\mmm{F}_{ef}'(0)}f,\phi).
$$
\end{lemma}
\begin{proof}[Proof of Lemma~\ref{n:2}]
Fix $r>0$. It follows from
$f\in\mcl{D}(\overline{\mmm{F}_{ef}'(0)})$ that there exists a sequence
 $\{f_i\}_{i=1}^{\infty}$,
$f_i\in\mcl{D}(\mmm{F}_{ef}'(0)),$ such that
$$\lim\limits_{i\to\infty}{\|f_i-f\|_{\phi}^{\mmm{F},\,r}}=0$$ and
$$\lim\limits_{i\to\infty}{\|\mmm{F}_{ef}'(0)f_i-\overline{\mmm{F}_{ef}'(0)}f\|_{\phi}^{\mmm{F},\,r}}=0.$$
Let $\Phi'$ be the minimal closed space such that $\Phi$ and $f_i,\,\, \mmm{F}_{ef}'(0)f_i$ for all $i\in\mmb{N}$. Then, according to Proposition~\ref{t:1}, we can choose a subsequence $\{m_k\}_{k=1}^{\infty}$
of the sequence $\{g_k\}_{k=1}^{\infty}$ such that the family of the functions
$\mathrm{T}^1_{s}:\Phi'\times\{\phi\}\mapsto\mmb{T},\,\,s\geq 0,$
is defined by the equality
$$\mathrm{T}^1_{s}(g,\phi)=\lim\limits_{k
\rightarrow\infty}{(\{\mathrm{F}(t/{m_k})\}^{[{m_k}\frac{s}{t}]}
g,\phi)}, \,\,g\in\Phi'.$$ Thus, it follows from Lemma~\ref{lem:lem-4}
(applied with parameters $h=f$, $g=\overline{\mmm{F}_{ef}'(0)}f$)
that
$$(\mathrm{T}^1_{s}
(f,\phi))^{'}_{s}=\mathrm{T}^1_{s}(\overline{\mmm{F}_{ef}'(0)}f,\phi),\,s<r.$$
 Since the restriction of $\mathrm{T}^1_{s}$  to the set
$\Phi\times\{\phi\}$ coincides with $\mathrm{T}_{s}$ we see that
$$(\mathrm{T}_{s}
(f,\phi))^{'}_{s}=\mathrm{T}_{s}(\overline{\mmm{F}_{ef}'(0)}f,\phi),\,s<r.
$$
Since $r$ is arbitrary the Lemma is proved. \end{proof}

\begin{definition}\label{deft:21} Let $\mmm{Z}$ be closable operator in $\mmf{X}$. A function
 $g:[0,\,l)\mapsto\mmm{X},\,\,l>0,$ is called a local solution  of the system
\begin{eqnarray}
f'(t)&=&\overline{\mmm{Z}}f(t),\,\,t\in[0,\,l),\nonumber\\
f(0)&=&h\in\mcl{D}(\overline{\mmm{Z}}),\nonumber
\end{eqnarray}
(on semi-interval $[0,\,l)$), if the following conditions are satisfied:

\noindent 1) $g(s)\in\mcl{D}(\overline{\mmm{Z}})$ for any
$s\in[0,\,l).$

\noindent 2) $g(0) = h$ and $(g(s))'_s=\overline{\mmm{Z}}g(s),$ for any $s\in[0,\,l).$
\end{definition}

\begin{proposition}\label{t:6}Assume that $\mathrm{F}\in\EuScript{F}_{\mmf{X}}$
and $t,\,\,l>0$. Assume that there exists a local solution
$f:[0,\,l)\mapsto\mmm{X}$ of the system
\begin{eqnarray}
f'(s) &=& \overline{\mmm{F}_{ef}'(0)}f(s),\,\,s\in[0,\,l),\label{eqn:eq-4}\\
f(0) &=& h\in\mcl{D}(\overline{\mmm{F}_{ef}'(0)}). \label{eqn:eq-41}
\end{eqnarray}
then the solution is unique and the following equality is valid:
\begin{eqnarray}
f(s)=w\mbox{-}\lim\limits_{n\to\infty}{\{\mathrm{F}(t/{n})\}^{[n\frac{s}{t}]}h},\,
s\in [0,\,l).\label{eqn:l3}
\end{eqnarray}
\end{proposition}
\begin{proof}[Proof of Proposition~\ref{t:6}.] Let $f(s),\,\,s\in[0,\,l),$
be a local solution of system
\eqref{eqn:eq-4}-\eqref{eqn:eq-41}. Fix $\phi\in\mmf{X}^*$.
Let $\mmf{B}$ be the minimal closed linear subspace of $\mmm{X}$ such that $f(r)\in\mmf{B}$ for any
$r\in\mmb{Q}\cap [0,\, l)$. Hence, $f_s\in\mmf{B}$ for any $s\in [0,
l)$ and, consequently, $\overline{\mmm{F}'(0)}f_s\in\mmf{B}$ for any
$s\in [0, l)$. Let us choose a subsequence
$\{g_k\}_{k=1}^{\infty}$ of  arbitrary strictly increasing natural sequence
$\{n_k\}_{k=1}^{\infty}$ such
that the family of the functions  $\mathrm{T}_s: \Phi\times\{\phi\}\mapsto\mathbb{T}$, $s\geq
0$, is defined by the equality \eqref{eqn:l2} with $\Phi=\mmf{B}$.
    From Lemma~\ref{n:2} it follows that
$$(\mmm{T}_s
(f(v),\phi))'_s=\mmm{T}_s(\overline{\mmm{F}'(0)}f(v),\phi),\,\,s,\,v\in[0,\,l).$$
  So, by Proposition~\ref{t:1},
\begin{eqnarray}
&&(\mmm{T}_s (f(v-s),\phi))'_s \nonumber\\&=& \lim\limits_{a\to
0}{{a}^{-1}(\mmm{T}_{s+a}(f(v-a-s),\phi)-\mmm{T}_{s}(f(v-s),\phi))}\nonumber\\&=&
\lim\limits_{a\to
0}{{a}^{-1}(\mmm{T}_{s+a}((f(v-a-s)-f(v-s)),\phi))}\nonumber\\&+&
\lim\limits_{a\to
0}{{a}^{-1}(\mmm{T}_{s+a}(f(v-s),\phi)-\mmm{T}_{s}(f(v-s),\phi))}\nonumber\\&=&
(-\mmm{T}_{s}(\overline{\mmm{F}'(0)}f(v-s),\phi)))+\mmm{T}_{s}(\overline{\mmm{F}'(0)}f(v-s),\phi)=0,\nonumber
\end{eqnarray}
for any $v\in (0,\,l),\,\,s\in(0,\,v).$ From this we infer that the function $\mmm{T}_s (f(v-s),\phi)$ is a constant w.r.t.
 $s$. Thus,
$$\mmm{T}_v (g,\phi)=( f(v),\phi),\,\,v\in[0,\,l).$$
Now the existence of limit~\eqref{eqn:l3} follows from the arbitrariness of the sequence $\{n_k\}_{k=1}^{\infty}$
for any fixed $\phi\in\mmf{X}^*$. Therefore, $f$ is unique.
  \end{proof}

As a consequence of Proposition~\ref{t:6} we get
\begin{theorem}\label{cor:2}Assume that $\mathrm{F}\in\EuScript{F}_{\mmf{X}}$
and $t>0$. Assume also that a linear operator $\mmm{Z}$ has the domain
$\mcl{D}(\mmm{Z})\subset\mcl{D}(\mmm{F}_{ef}'(0))$ and
$$\mmm{Z}g=\mmm{F}_{ef}'(0)g,\,\,g\in\mcl{D}(\mmm{Z}).$$ If there exists a local solution
 $f:[0,\,l)\mapsto\mmm{X}$ of the system
\begin{eqnarray}
f'(s)&=&\overline{\mmm{Z}}f(s),\,\,s\in[0,\,l),\nonumber\\
f(0)&=&h\in\mcl{D}(\overline{\mmm{Z}}),\nonumber
\end{eqnarray}
then the solution is unique and the following equality is valid:
\begin{eqnarray}
f(s)=w\mbox{-}\lim\limits_{n\to\infty}{\{\mathrm{F}(t/{n})\}^{[n\frac{s}{t}]}h},\,\,s\in
[0,\,l).\nonumber
\end{eqnarray}
\end{theorem}
\begin{corollary}\label{cor_t:7}Assume that the conditions of Theorem~\ref{cor:2} are satisfied. Then the equality
$$f'(s)=w\mbox{-}\lim\limits_{n\to\infty}{\{\mathrm{F}(t/{n})\}^{[n\frac{s}{t}]}
\overline{\mmm{Z}}f(0)},$$
holds for any $s\in[0,\, l)$.
\end{corollary}
\begin{proof}[Proof of Corollary~\ref{cor_t:7}.]
Let $\mmf{B}$ be the same as in the proof of Proposition~\ref{t:6}. For any $\phi\in\mmf{X}^{*}$ choose a subsequence $\{g_k\}_{k=1}^{\infty}$ of  arbitrary strictly increasing sequence
$\{n_k\}_{k=1}^{\infty},\,\,n_k\in\mmb{N},$ such that the family of the functions  $\mathrm{T}_s:
\Phi\times\{\phi\}\mapsto\mathbb{T}$, $s\geq 0,$ is defined by the equality
\eqref{eqn:l2} with $\Phi=\mmf{B}$. From Lemma~\ref{n:2} and Theorem~\ref{cor:2} it follows that
$$(f'(s)
,\phi)=\lim\limits_{k\to\infty}{(\{\mathrm{F}(t/{g_k})\}^{[g_k\frac{s}{t}]}
\overline{\mmm{Z}}f(0),\phi)},\,\,s\in[0,\, l).$$ Since
 $\{n_k\}_{k=1}^{\infty}$ is arbitrary for any fixed $\phi\in\mmf{X}^{*}$, the Corollary is proved.   \end{proof}
\begin{theorem}\label{sec:1}
Assume that $\mathrm{F}\in\EuScript{F}_{\mmf{X}}$. Assume also that a linear operator  $\mmm{Z}$ has the
domain
$\mcl{D}(\mmm{Z})\subset\mcl{D}(\mmm{F}_{ef}'(0))$ and the following equality is valid:
$$\mmm{Z}g=\mmm{F}_{ef}'(0)g,\,\,g\in\mcl{D}(\mmm{Z}).$$ Let
$f:[0,\,l)\mapsto\mmm{X}$ be a local solution of the system
\begin{eqnarray}
f'(t)&=&\overline{\mmm{Z}}f(t),\,\,t\in[0,\,l),\nonumber\\
f(0)&=&h\in\mcl{D}(\overline{\mmm{Z}}).\nonumber
\end{eqnarray}
Then $f\in C_1([0,\,l),\mmm{X})$.
\end{theorem}
\begin{proof}[Proof of Theorem~\ref{sec:1}.]
It is sufficient to show that for $\alpha \in \Omega$,
$t\in[0,\,l)$  and $\epsilon>0$  there exists $\delta>0$ such that
 $\|f'(t)-f'(s)\|_{\alpha}\leq\epsilon$ for
$s\in[t-\delta,\,t+\delta]\cap[0,\,l).$ By Lemma~\ref{t:7}, we can choose $\delta_0>0$ and $n_0\in\mmb{N}$ such that
$$\|\{\mathrm{F}(t/{n})\}^{[n\frac{t}{t}]}
\overline{\mmm{Z}}f(0)-\{\mathrm{F}(t/{n})\}^{[n\frac{s}{t}]}
\overline{\mmm{Z}}f(0)\|_{\alpha}\leq\epsilon,\,\, n\in\mmb{N},$$
 for
$s\in[t-\delta,\,t+\delta]\cap[0,\,l)$, $t\in[0,\,l)$ and $n>n_0.$ Fix
$s\in[t-\delta,\,t+\delta]\cap[0,\,l)$.  By Hahn-Banach Theorem, there exists $\phi\in\mmm{X}^*$ such that $|(x,\,\phi)|\leq
\|\phi(x)\|_{\alpha}$ for all $x\in\mmm{X}$ and
$(f'(t)-f'(s),\phi)=\|f'(t)-f'(s)\|_{\alpha}$. By Corollary~\ref{cor_t:7}
$$\{\mathrm{F}(t/{n})\}^{[n\frac{t}{t}]}
\overline{\mmm{Z}}f(0)-\{\mathrm{F}(t/{n})\}^{[n\frac{s}{t}]}
\overline{\mmm{Z}}f(0)\to f'(t)-f'(s)$$ as $n\to 0$ in the weak topology. Thus,
\begin{eqnarray}\nonumber\lim\limits_{n\to\infty}{(\{\mathrm{F}(t/{n})\}^{[n\frac{t}{t}]}
\overline{\mmm{Z}}f(0)-\{\mathrm{F}(t/{n})\}^{[n\frac{s}{t}]}
\overline{\mmm{Z}}f(0),\phi)}\\=(f'(t)-f'(s),\phi)=\|f'(t)-f'(s)\|_{\alpha}
\leq\epsilon.\nonumber
\end{eqnarray} By the arbitrariness of
$s\in[t-\delta,\,t+\delta]\cap[0,\,l)$ the Theorem is proved.
  \end{proof}

\section{Chernoff type formula for local solution of the abstract Cauchy
problem.}\label{2.1}

\begin{theorem}\label{chern1}
Assume that $\mathrm{F}\in\EuScript{F}_{\mmf{X}}$.
Assume also that a linear operator $\mmm{Z}$ has the
domain $\mcl{D}(\mmm{Z})\subset\mcl{D}(\mmm{F}_{ef}'(0))$ and
$$\mmm{Z}g=\mmm{F}_{ef}'(0)g,\,\,g\in\mcl{D}(\mmm{Z}).$$ If there exists a local solution $f:[0,\,l)\mapsto\mmm{X}$ of the system
\begin{eqnarray}
f'(t)&=&\overline{\mmm{Z}}f(t),\,\,t\in[0,\,l),\nonumber\\
f(0)&=&h\in\mcl{D}(\overline{\mmm{Z}}),\nonumber
\end{eqnarray}
then $\mathrm{F}({\textstyle\frac{t}{n}})^{n}h$ tends to
$f(t)$ as $n\to\infty$ uniformly with respect to $t\in [0,\,t_0]$ for any
$t_0\in(0,\,l).$
\end{theorem}
\begin{remark}
Closability of the operator $\mmm{Z}$ of Theorem~\ref{chern1}
follows from Theorem~\ref{cor:1}.
\end{remark}
\begin{proof}[Proof of Theorem~\ref{chern1}.] \,
Fix arbitrary $t_0\in(0,\,l)$, $\alpha\in\Omega$
and $\epsilon>0$. It is sufficient to show
that there exists
$n_0\in\mmb{N}$ such that
$$\sup\limits_{t\in[0,\,t_0]}{\|\mmm{F}({\textstyle\frac{t}{n}})^n
h-f(t)\|_{\alpha}}<\epsilon$$ for any natural number
$n\geq n_0$. Denote $c_1=\epsilon/8t_0$, $c_2=\epsilon/4t_0$ and
$a_i=\min{(\epsilon/24i_0,\epsilon/32t_0)}$ for each
$i\in\mmb{N}$. Further, for each $i\in\mmb{N}$ choose a set
$\mmf{A}_i=\{y^i_0,\,y^i_1,\,\ldots \,
,y^i_i\}\subset\mcl{D}({\mmm{Z}})$ such that
\begin{eqnarray}\|\overline{\mmm{Z}}f({\textstyle\frac{t_0
k}{i}})-\mmm{Z}y^i_k\|_{\alpha}^{\mmm{F},t_0;\mmm{F},t_0}+\|f({\textstyle\frac{t_0
k}{i}})-y^i_k\|_{\alpha}^{\mmm{F},t_0;\mmm{F},t_0}<a_i,\label{e:12}
\end{eqnarray} where $k\in [0,\,i]\cap\mmb{Z}$. By Theorem \ref{sec:1},
we can choose $m_0\in\mmb{N}$
such that
\begin{eqnarray}
\|\overline{\mmm{Z}}f(r)-\overline{\mmm{Z}}f(s)\|_{\alpha}^{\mmm{F},t_0;\mmm{F},t_0}+
\|f(r)-f(s)\|_{\alpha}^{\mmm{F},t_0;\mmm{F},t_0}\leq c_1,
\label{e:10}
\end{eqnarray}
if $|r-s|<t_0/m_0,\,\,r,s\in [0,\,t_0]$. By Lemma~\ref{t:5}, we can choose $l_0>m_0$ such that from
$t\frac{i}{d}<t_0/l_0$, $i,\,d\in\mmb{N},\,\,t\geq 0$, follows the inequality
\begin{eqnarray}
{\textstyle\|(\mmm{F}(\frac{t}{d})^i-\mmm{I})\overline{\mmm{Z}}f(t_0\frac{
k}{m_0})\|_{\alpha}^{\mmm{F},t_0}\leq c_2} \label{e:9}
\end{eqnarray}
for each $k\in [0,\,m_0-1]\cap\mmb{Z}$. By inequalities \eqref{e:10}
and \eqref{e:9}, we get
\begin{eqnarray}
\|(\mmm{F}({\textstyle\frac{t}{d}})^i-\mmm{I})
\overline{\mmm{Z}}f(s)\|_{\alpha}^{\mmm{F},t_0}\leq
c_2+2c_1,\,\,s\in [0,\,t_0],\label{es:11}
\end{eqnarray}
if  $t\frac{i}{d}<t_0/l_0$, $i,\,d\in\mmb{N},\,\,t\geq 0$.
Fix natural $i_0\in(4 l_0,\,\infty)$.
 Let natural $d_0>2i_0$ such that
\begin{eqnarray}
 \|\overline{\mmm{Z}}f(r)-\overline{\mmm{Z}}f(s)\|_{\alpha}^{\mmm{F},t_0;\mmm{F},t_0}+
\|f(r)-f(s)\|_{\alpha}^{\mmm{F},t_0;\mmm{F},t_0}\leq a_{i_0},
\label{e:15}
\end{eqnarray}
if $|r-s|<t_0/d_0,$ $r,s\in [0,\,t_0]$. Then for any natural $n>d_0$ we can choose a set
$\{l_{n,d}\in\mmb{N}\cup\{0\}|d\in [0,\,i_0]\cap\mmb{Z}\}$ such that $l_{n,0}=0,\, l_{n,i_0}=n$ and
\begin{eqnarray}
|t_0l_{n,i}/n-t_0i /i_0|<t_0 /d_0.\label{dev:232}\end{eqnarray}
So, we have
\begin{eqnarray} 0<1/i_0 -2/d_0<l_{n,i}-i /i_0+1/i_0+(i-1)/i_0-l_{n,i-1}\label{dev:231}\\
=l_{n,i}-l_{n,i-1}<1/i_0 +2/d_0<2n/i_0,\,\,i\in
[1,\,i_0]\cap\mmb{N}.\label{dev:2}
\end{eqnarray} By inequalities \eqref{e:12}, \eqref{e:15} and
\eqref{dev:232} we get
\begin{eqnarray} \|\overline{\mmm{Z}}f({\textstyle\frac{t_0
l_{n,i}}{n}})-\mmm{Z}y^{i_0}_i\|_{\alpha}^{\mmm{F},t_0;\mmm{F},t_0}+\|f({\textstyle\frac{t_0
{l_{n,i}}}{n}})-y^{i_0}_i\|_{\alpha}^{\mmm{F},t_0;\mmm{F},t_0}<2a_{i_0}
\label{e:14}\end{eqnarray} for each $i\in [1,\,i_0]\cap\mmb{N}.$ Therefore,
by inequalities
 \eqref{es:11} and \eqref{e:14}, we get
\begin{eqnarray}
\|(\mmm{F}(t/d)^m-\mmm{I})\mmm{Z}y^{i_0}_i\|_{\alpha}^{\mmm{F},t_0}
\leq c_2+2c_1+4a_{i_0}\label{e:18}
\end{eqnarray}
for each $i\in [1,\,i_0]\cap\mmb{N}$, if $t\frac{m}{d}<t_0/l_0$,
$m,\,d\in\mmb{N},\,\,t\geq 0$. Choose $\{v^i_s\}_{s\geq 0}$,
$v^i_s\in\mmf{X}$, $i\in[0,\,i_0]\cap\mmb{Z},$ such that there exist
$\lim\limits_{s\to
0}{v^i_s}=y^{i_0}_i\,\,\mbox{and}\,\,\lim\limits_{s\to
0}{s^{-1}(\mmm{F}(s)-\mmm{I})v^i_s}=\mmm{Z}y^{i_0}_i.$ Choose $k_i$
such that from  $s\in[0,\,{\textstyle\frac{t_0}{k_i}})$ follows
\begin{eqnarray}\|v^i_s-y^{i_0}_i\|^{\mmm{F},t_0}_{\alpha}
+\|s^{-1}(\mmm{F}(s)-\mmm{I})v^i_s
-\mmm{Z}y^{i_0}_i\|^{\mmm{F},t_0;\mmm{F},t_0}_{\alpha}\leq
a_{i_0}.\label{dev:3}
\end{eqnarray}
Define $k=\sup\limits_{i\in[0,\,i_0]\cap\mmb{N}}{k_i}.$
Fix arbitrary
 $t\in(0,\,t_0]$ and natural $n>\max{(d_0, k)}.$ Let
$i_{t,n}=\sup\limits_{t_0 l_{n,i}/n<t}{(i+1)}$. Denote
$g_{i}=\min{(t_0l_{n,i}/n,t)}$ and $l_{n,i}'=[n g_{i}/t]$ for
$i\in[0,\,i_{t,n}]\cap\mmb{Z}.$ Then we have
\begin{eqnarray}
\label{dev:11}&&\|\mmm{F}({\textstyle\frac{t}{n}})^n
f-f(t)\|_{\alpha}
\\&\leq&
\sum_{i=0}^{i_{t,n}-1}{\|\mmm{F}({\textstyle\frac{t}{n}})^{n-l'(n,i)}f(g_i)}
-
{\mmm{F}({\textstyle\frac{t}{n}})^{n-l'(n,i+1)}f(g_{i+1})\|_{\alpha}}\,\,\,\,\,\,\\&\leq&
\sum_{i=0}^{i_{t,n}-1}{\|\mmm{F}({\textstyle\frac{t}{n}})^{l'(n,i+1)-l'(n,i)}f(g_i)}
- f(g_{i+1})\|^{\mmm{F},t_0}_{\alpha}\\&\leq&
\sum_{i=0}^{i_{t,n}-1}{\|\mmm{F}({\textstyle\frac{t}{n}})^{m_i}f(g_i)}
-
f(g_{i}+{\textstyle\frac{t}{n}}m_i)\|^{\mmm{F},t_0}_{\alpha},\label{in:1}
\end{eqnarray}
where $m_i=l'(n,i+1)-l'(n,i)$. For any
$i\in[0,\,i_{t,n})\cap\mmb{Z}$ we see that
\begin{eqnarray}
 \label{dev:10}&&{\|\mmm{F}({\textstyle\frac{t}{n}})^{m_i}f(g_i)} -
f(g_{i}+{\textstyle\frac{t}{n}}m_i)\|^{\mmm{F},t}_{\alpha}\\&\leq&
\|\mmm{F}({\textstyle\frac{t}{n}})^{m_i}f(g_i)-f(g_i)
-{\textstyle\frac{t}{n}}{m_i}\overline{\mmm{Z}}f(g_i)\|_{\alpha}^{\mmm{F},t_0}
\label{e:1}\\&+&\|f(g_{i}+{\textstyle\frac{t}{n}}m_i)-f(g_i)
-{\textstyle\frac{t}{n}}{m_i}\overline{\mmm{Z}}f(g_i)\|_{\alpha}^{\mmm{F},t_0}\\
&=&(A_i)+(B_i)\label{e:2}
\end{eqnarray}
  By inequalities \eqref{dev:3} and \eqref{e:14},
  we have
\begin{eqnarray}
&&(A_i)=\|\mmm{F}({\textstyle\frac{t}{n}})^{m_i}f(g_i)-f(g_i)
-{\textstyle\frac{t}{n}}{m_i}\overline{\mmm{Z}}f(g_i)\|_{\alpha}^{\mmm{F},t_0}\label{dev:4}\\&\leq&
\|\mmm{F}({\textstyle\frac{t}{n}})^{m_i}v^i_{t/n}-v^i_{t/n}
-{\textstyle\frac{t}{n}}{m_i}\mmm{Z}y^{i_0}_i\|_{\alpha}^{\mmm{F},t_0}
\\&+& 2\|v^i_{t/n}-f(g_i)\|_{\alpha}^{\mmm{F},t_0}+
{\textstyle\frac{t}{n}}{m_i}\|\overline{\mmm{Z}}f(g_i)-\mmm{Z}y^{i_0}_i
\|_{\alpha}^{\mmm{F},t_0}\\&\leq&
\|\mmm{F}({\textstyle\frac{t}{n}})^{m_i}v^i_{t/n}-v^i_{t/n}
-{\textstyle\frac{t}{n}}{m_i}\mmm{Z}y^{i_0}_i\|_{\alpha}^{\mmm{F},t_0}
\\&+&2\|v^i_{t/n}-y^{i_0}_i\|_{\alpha}^{\mmm{F},t_0}+2\|y^{i_0}_i-f(g_i)
\|_{\alpha}^{\mmm{F},t_0}\\&+&
{\textstyle\frac{t}{n}}{m_i}\|\overline{\mmm{Z}}f(g_i)-\mmm{Z}y^{i_0}_i
\|_{\alpha}^{\mmm{F},t_0}\,\,\,\,\\&\leq&
\|\mmm{F}({\textstyle\frac{t}{n}})^{m_i}v^i_{t/n}-v^i_{t/n}
-{\textstyle\frac{t}{n}}{m_i}\mmm{Z}y^{i_0}_i\|_{\alpha}^{\mmm{F},t_0}
\\&+&2a_{i_0}+4a_{i_0}+
{\textstyle\frac{t}{n}}{m_i}2a_{i_0}\label{dev:5}
\end{eqnarray}
for any $i\in[0,\,i_{t,n})\cap\mmb{Z}.$
 Note that
\begin{eqnarray}\nonumber{\textstyle t\frac{m_i}{n}}&=&
{\textstyle t\frac{l'(n,i+1)-l'(n,i)}{n}\leq
t\frac{(g_{i+1}-g_{i})n/t +2}{n}}
\\&\leq&\nonumber
g_{i+1}-g_{i}+2t_0/n<t_0(l_{n,i+1}-l_{n,i+1})/n+2t_0/i_0\\&<&
2t_0/i_0+t_0/2l_0<t_0/l_0,\,\,i\in[0,\,i_{t,n})\cap\mmb{Z}.\nonumber
\end{eqnarray}
Here we have used inequalities  \eqref{dev:231}-\eqref{dev:2}.
So, by inequalities \eqref{e:18} and \eqref{dev:3}, we have
\begin{eqnarray}
&&\|\mmm{F}({\textstyle\frac{t}{n}})^{m_i}v^i_{t/n}-v^i_{t/n}
-{\textstyle\frac{t}{n}}{m_i}\mmm{Z}y^{i_0}_i\|_{\alpha}^{\mmm{F},t_0}\label{e:3}\\&\leq&
{\textstyle\frac{t}{n}}\|m_i{\textstyle\frac{(\mmm{F}(\frac{t}{n})^{m_i-1}+\cdots+\mmm{I})}{m_i}}{\textstyle\frac{n}{t}}
(\mmm{F}({\textstyle\frac{t}{n}})-\mmm{I})v^i_{t/n}
-{m_i}\mmm{Z}y^{i_0}_i\|_{\alpha}^{\mmm{F},t_0}
\\
&\leq&{\textstyle\frac{t}{n}}m_i\|{\textstyle\frac{(\mmm{F}(\frac{t}{n})^{m_i-1}
+\cdots+\mmm{I})}{m_i}}{\textstyle\frac{n}{t}}
((\mmm{F}({\textstyle\frac{t}{n}})-\mmm{I})v^i_{t/n}-\mmm{Z}y^{i_0}_i)
\|_{\alpha}^{\mmm{F},t_0}\\&+&
{\textstyle\frac{t}{n}}m_i\|({\textstyle\frac{(\mmm{F}(\frac{t}{n})^{m_i-1}
+\cdots+\mmm{I})}{m_i}}-\mmm{I})\mmm{Z}y^{i_0}_i\|_{\alpha}^{\mmm{F},t_0}\\
&\leq&{\textstyle\frac{t}{n}}m_i\|{\textstyle\frac{n}{t}}
(\mmm{F}({\textstyle\frac{t}{n}})-\mmm{I})v^i_{t/n}-\mmm{Z}y^{i_0}_i
\|_{\alpha}^{\mmm{F},t_0;\mmm{F},t_0}\\&+&
{\textstyle\frac{t}{n}}m_i\|({\textstyle\frac{(\mmm{F}(\frac{t}{n})^{m_i-1}
+\cdots+\mmm{I})}{m_i}}-\mmm{I})\mmm{Z}y^{i_0}_i\|_{\alpha}^{\mmm{F},t_0}\\
&\leq&{\textstyle\frac{t}{n}}m_ia_{i_0}+{\textstyle\frac{t}{n}}
m_i(c_2+2c_1+4a_{i_0}),\,\,i\in[0,\,i_{t,n})\cap\mmb{Z}.\label{e:4}
\end{eqnarray}
 Combining inequalities
\eqref{dev:4}--\eqref{dev:5} and \eqref{e:3}--\eqref{e:4}, we get
\begin{eqnarray}
(A_i)\leq
6a_{i_0}+{\textstyle\frac{t}{n}}m_i(c_2+2c_1+7a_{i_0}),\,\,i\in[0,\,i_{t,n})\cap\mmb{Z}.\label{dev:8}
\end{eqnarray}
By inequality \eqref{e:15}, we have
\begin{eqnarray}
&&(B_i)=\|f(g_{i}+{\textstyle\frac{t}{n}}m_i)-f(g_i)
-{\textstyle\frac{t}{n}}{m_i}\overline{\mmm{Z}}f(g_i)\|_{\alpha}^{\mmm{F},t_0}\label{dev:6}
\\&\leq&\|\int_{0}^{{\textstyle\frac{t}{n}}m_i}{\ov{\mmm{Z}}f(g_i+s)\,d
s}-{\textstyle\frac{t}{n}}{m_i}\overline{\mmm{Z}}f(g_i)\|_{\alpha}^{\mmm{F},t_0}
\\&\leq&\|\int_{0}^{{\textstyle\frac{t}{n}}m_i}{\ov{\mmm{Z}}(f(g_i+s)-f(g_i))\,d
s}\|_{\alpha}^{\mmm{F},t_0}\leq {\textstyle\frac{t}{n}}m_ia_{i_0}, i\in[0,\,i_{t,n})\cap\mmb{Z}.
\label{dev:7}
\end{eqnarray}
 Therefore, by inequalities \eqref{dev:8} and \eqref{dev:6}--\eqref{dev:7} we infer that
\begin{eqnarray}(A_i)+(B_i)=
6a_{i_0}+{\textstyle\frac{t}{n}}m_i(c_2+2c_1+
8a_{i_0}),\,\,i\in[0,\,i_{t,n})\cap\mmb{Z}.
\end{eqnarray}
Thus, by inequalities \eqref{dev:11}--\eqref{in:1} and
\eqref{dev:10}--\eqref{e:2}, we get
\begin{eqnarray}
 \|\mmm{F}({\textstyle\frac{t}{n}})^n h-f(t)\|_{\alpha}&\leq&
\sum_{i=0}^{i_{t,n}-1}{6a_{i_0}+{\textstyle\frac{t}{n}}m_i(c_2+2c_1+
8a_{i_0})}\\&\leq& 6i_0a_{i_0}+t_0(c_2+2c_1+
8a_{i_0})\leq\epsilon\label{e:20}
\end{eqnarray}
for  $n>n_0=\max{(d_0, k)}.$  By the arbitrariness $\epsilon>0$,
$\alpha\in\Omega$ and $t_0\in(0,\,l)$ the Theorem is proved.
  \end{proof}

\begin{corollary}\label{chern1b}
Assume that $\mathrm{F}\in\EuScript{F}_{\mmf{X}}$. Assume also that a linear operator $\mmm{Z}$ has the domain
$\mcl{D}(\mmm{Z})\subset\mcl{D}(\mmm{F}_{ef}'(0))$ and
$$\mmm{Z}g=\mmm{F}_{ef}'(0)g,\,\,g\in\mcl{D}(\mmm{Z}).$$ If there exists a local solution $f:[0,\,l)\mapsto\mmm{X}$ of the system
\begin{eqnarray}
f'(t)&=&\overline{\mmm{Z}}f(t),\,\,t\in[0,\,l),\nonumber\\
f(0)&=&h\in\mcl{D}(\overline{\mmm{Z}}),\nonumber
\end{eqnarray}
then $\mathrm{F}({\textstyle\frac{t}{n}})^{[sn/t]}h$ tends to
$f(s)$ as $n\to\infty$ uniformly with respect to $s\in [0,\,t_0]$ for any
$t_0\in(0,\,l).$
\end{corollary}

\begin{proof}[Proof of Corollary~\ref{chern1b}.]
Note that
\begin{eqnarray}
\mathrm{F}({\textstyle\frac{t}{n}})^{[sn/t]}h - f(s) =
(\mathrm{F}({\textstyle\frac{t[sn/t]/n}{[sn/t]}})^{[sn/t]}h -
f(t[sn/t]/n))
\nonumber\\
+ (f(t[sn/t]/n) - f(s))=(A)+(B)\nonumber
\end{eqnarray}
From Theorem~\ref{chern1} it follows that $(A)$ tends to $0$
as $n\to\infty$ uniformly with respect to $s\in [0,\,t_0]$ for each
$t_0\in(0,\,l).$ From Theorem~\ref{sec:1} we infer that $(B)$
tends to $0$ as $n\to\infty$ uniformly with respect to $s\in
[0,\,t_0]$ for each $t_0\in(0,\,l).$ From this we get the Corollary.
\end{proof}

\begin{theorem}\label{chern1c}
Assume that $\mathrm{F}\in\EuScript{F}_{\mmf{X}}$.
Assume also that a linear operator $\mmm{Z}$ has the domain $\mcl{D}(\mmm{Z})\subset\mcl{D}(\mmm{F}_{ef}'(0))$ and
$$\mmm{Z}g=\mmm{F}_{ef}'(0)g,\,\,g\in\mcl{D}(\mmm{Z}).$$ If there exists a local solution $f:[0,\,l)\mapsto\mmm{X}$ of the system
\begin{eqnarray}
f'(t)&=&\overline{\mmm{Z}}f(t),\,\,t\in[0,\,l),\nonumber\\
f(0)&=&h\in\mcl{D}(\overline{\mmm{Z}}),\nonumber
\end{eqnarray}
then
$\mathrm{F}({\textstyle\frac{t}{n}})^{n}\overline{\mmm{Z}}h$
tends to $\overline{\mmm{Z}}f(t)$ as $n\to\infty$ uniformly with respect to
$t\in [0,\,t_0]$ for any $t_0\in(0,\,l).$
\end{theorem}
\begin{proof}[Proof of Theorem~\ref{chern1c}.] For arbitrary $t_0\in(0,\,l)$ fix $\delta>0$ such that $t_0+\delta<l$.
From Theorem~\ref{chern1} it follows that
\begin{eqnarray}
\mathrm{F}({\textstyle\frac{t}{n}})^{n}(f(m)-h)/m \to (f(t+m) -f(t))/m
\label{eqny:3}
\end{eqnarray}
as $n\to\infty$ uniformly with respect to $t\in [0,\,t_0]$ for any
$m\in (0,\,\delta ]$. Note that $\lim\limits_{m\to
0}{(f(m)-h)/m}=\overline{\mmm{Z}}h$ and the family of the operators
$\{\mathrm{F}({\textstyle\frac{t}{n}})^{n}|\,n\in\mmb{N},\,t\in
[0,\,t_0]\}$ is equicontinuous. From this it follows that there exists
$$\lim\limits_{n\to\infty}{\mathrm{F}({\textstyle\frac{t}{n}})^{n}\overline{\mmm{Z}}h}=\overline{\mmm{Z}}f(t)$$ for each $t\in [0,\,t_0]$ and the limit in expression~\eqref{eqny:3} is uniform with respect to  $m\in (0,\,\delta ]$. Thus, we infer that
$$\lim_{n\to\infty}{\mathrm{F}({\textstyle\frac{t}{n}})^{n}\overline{\mmm{Z}}h}=\overline{\mmm{Z}}f(t)$$
is uniform with respect to  $t\in [0,\,t_0]$. By the arbitrariness of
$t_0\in(0,\,l)$, the Theorem is proved.
\end{proof}

\begin{corollary}\label{chern1d}
Assume that  $\mathrm{F}\in\EuScript{F}_{\mmf{X}}$. Assume also that a linear operator $\mmm{Z}$ has the domain
$\mcl{D}(\mmm{Z})\subset\mcl{D}(\mmm{F}_{ef}'(0))$ and
$$\mmm{Z}g=\mmm{F}_{ef}'(0)g,\,\,g\in\mcl{D}(\mmm{Z}).$$ If there exists a local solution $f:[0,\,l)\mapsto\mmm{X}$ of the system
\begin{eqnarray}
f'(t)&=&\overline{\mmm{Z}}f(t),\,\,t\in[0,\,l),\nonumber\\
f(0)&=&h\in\mcl{D}(\overline{\mmm{Z}}),\nonumber
\end{eqnarray}
then
$\mathrm{F}({\textstyle\frac{t}{n}})^{[sn/t]}\overline{\mmm{Z}}h$
tends to $\overline{\mmm{Z}}f(s)$ as $n\to\infty$ uniformly with respect to
$s\in [0,\,t_0]$ for any $t_0\in(0,\,l).$
\end{corollary}

\begin{proof}[Proof of Corollary~\ref{chern1d}.]
Note that
\begin{eqnarray}
\mathrm{F}({\textstyle\frac{t}{n}})^{[sn/t]}\overline{\mmm{Z}}h -
\overline{\mmm{Z}}f(s) =
(\mathrm{F}({\textstyle\frac{t[sn/t]/n}{[sn/t]}})^{[sn/t]}\overline{\mmm{Z}}h
- \overline{\mmm{Z}}f(t[sn/t]/n))
\nonumber\\
+ (\overline{\mmm{Z}}f(t[sn/t]/n) -
\overline{\mmm{Z}}f(s))=(A_n)+(B_n)\nonumber
\end{eqnarray}
From Theorem~\ref{chern1c} we infer that $(A_n)$ tends to $0$
as $n\to\infty$ uniformly with respect to $s\in [0,\,t_0]$ for each
$t_0\in(0,\,l).$ From Theorem~\ref{sec:1}  it follows that $(B_n)$
tends to $0$ as $n\to\infty$ uniformly with respect to $s\in
[0,\,t_0]$ for each $t_0\in(0,\,l).$
\end{proof}

\section{Chernoff, Lie-Trotter and Trotter-Kato theorems for sequentially complete locally convex spaces.}\label{2.2}

The following theorem is the extension of Chernoff theorem~\cite{Chernoff}.

\begin{theorem}\label{chern}
Let  $\mmm{Z}$ be a generator of the $lC_0$-semigroup in sequentially complete locally convex space $\mmf{X}$ and
$\EuScript{D}$ be a core of the generator $\mmm{Z}$. Assume that there exists a function
$\mmm{F}\in\EuScript{F}_{\mmf{X}}$ such that
$$\mmm{F}_{ef}'(0)f =\mmm{Z}f$$
for each $f\in\EuScript{D}.$ Then $\mmm{F}(t/n)^n f$ converges to
$\exp{(t\mmm{Z})}f$ as $n\to\infty$ for all $f\in\mmf{X}$
uniformly with respect to $t\in[0,\,t_0]$ for any $t_0>0$.
\end{theorem}

\begin{proof}[Proof of Theorem~\ref{chern}.]
From Theorem~\ref{chern1} it follows that for each
$f\in\EuScript{D}$ $\mmm{F}(t/n)^n f$ tends to $\exp{(t\mmm{Z})}f$
as $n\to\infty$ uniformly with respect to $t\in[0,\,t_0]$ for any
$t_0>0$. From equicontinuity of the set of the functions
$\{\mathrm{F}({\textstyle\frac{t}{n}})^{n}|\,n\in\mmb{N},\,t\in
[0,\,t_0]\}$ we deduce that $\mmm{F}(t/n)^n f$ tends to
$\exp{(t\mmm{Z})}f$ as $n\to\infty$ for each $f\in\mmf{X}$ uniformly
with respect to $t\in[0,\,t_0]$ for any $t_0>0$.
\end{proof}

The following theorem can be considered as a extension of Trotter-Kato theorem for $lC_0$-semigroups in sequentially complete locally convex spaces. It states the equivalence between convergences of $lC_0$-semigroups and their generators.

\begin{theorem}\label{t-k}
For any $s\geq 0$ let $\mmm{Z}_s$ be a generator of
the $lC_0$-semigroup in sequentially complete locally convex space $\mmf{X}$. Let also the set
$$\mmf{A}_{l_0}=\{\exp{(l\mmm{Z}_s)}|l\in[0,\,l_0],\,\,s\geq 0\}$$
be equicontinuous for any $l_0\geq 0$. Let
$\EuScript{D}$ consist of all $f\in\mcl{D}(\mmm{Z}_0)$ for which
there exists a sequence $\{f_s\}_{s>0}$,
$f_s\in\mcl{D}(\mmm{Z}_s),$ such that
$$\lim\limits_{s\to 0}{f_s}=f$$ and $$\lim\limits_{s\to 0}{\mmm{Z}_s
f_s}=\mmm{Z}_0 f. $$ Then the following conditions are equivalent:

\noindent i) $\exp{(l\mmm{Z}_s)}f$ converges to $\exp{(l\mmm{Z}_0)}f$
as $s \to 0$ uniformly with respect to $l\in[0,\,t_0]$ for any $t_0>0$ and
$f\in\mmf{X}.$

\noindent ii) $\EuScript{D}$  is a core of the generator $\mmm{Z}_0$.

\end{theorem}
\begin{proof}[Proof of Theorem~\ref{t-k}.]
 Put $\mmm{T}_s(l)=\exp{(l\mmm{Z}_s)}$. For any
$\alpha\in\Omega$ and $d>0$ define the function
$\|\cdot\|_{\alpha}^{\mmf{A}_{d}}:\mmf{X}\to\mmb{R}$ by the equality
 $$\|x\|_{\alpha}^{\mmf{A}_{d}}=\sup\limits_{g\in\mmf{A}_{d}}{\|g(x)\|_{\alpha}},\,\,
 x\in\mmf{X}.$$ Firstly, let us show implication $(\Rightarrow)$.
It is easy to show that the set of elements of the type
$$\int_{s_1}^{s_2}{\mmm{T}_s(l)f\,dl},$$ $s_1, s_2\geq 0$,
$f\in\mmf{X},$ is a core of $\mmm{Z}_s$ for each
$s\geq 0$. Particularly, it follows from Theorem~\ref{cog:1}. Note that for any element of the type $$\int_{s_1}^{s_2}{\mmm{T}_0(l)f\,dl},$$
 $s_1, s_2\geq 0$, $f\in\mmf{X},$ the
following conditions are satisfied:
$$\lim\limits_{s\to 0}{\int_{s_1}^{s_2}{\mmm{T}_s(l)f\,dl}}=
\int_{s_1}^{s_2}{\mmm{T}_0(l)f\,dl}$$ and
\begin{eqnarray}\lim\limits_{s\to
0}{\mmm{Z}_s\int_{s_1}^{s_2}{\mmm{T}_s(l)f\,dl}} = \lim\limits_{s\to
0}{(\mmm{T}_s(s_2)f-\mmm{T}_s(s_1)f)}\\=
(\mmm{T}_0(s_2)-\mmm{T}_0(s_1))f=\mmm{Z}_0\int_{s_1}^{s_2}{\mmm{T}_0(l)f\,dl}
\end{eqnarray}
Therefore, $$\int_{s_1}^{s_2}{\mmm{T}_0(l)f\,dl}\in\EuScript{D}$$
and we get implication $(\Rightarrow)$.

\noindent Now let us show implication $(\Leftarrow)$. We will argue by contradiction. Assume that implication $(\Leftarrow)$
is not true. Then there exist $\alpha_0\in\Omega$, $t_0,\delta>0$,
$g\in\mmf{X}$, decreasing sequence
$\{s_k\}_{k\in\mmb{N}},\,\,s_k>0,$ such that
$\lim\limits_{k\to\infty}{s_k}=0$ and
\begin{eqnarray}\liminf\limits_{k\to\infty}{ \sup\limits_{s\in[0,\,t_0]}
{\|\mmm{T}_{s_k}(s)g-\mmm{T}_{0}(s)g\|_{\alpha_0}}}>\delta.\label{fg:8}
\end{eqnarray}
Firstly, let us show that
\begin{eqnarray}\lim\limits_{l\to
0}{\sup\limits_{k\in\mmb{N}}{\|\exp{(l\mmm{Z}_{s_k})}f-f\|_{\alpha}}}=0\label{fg:1}
\end{eqnarray}
for each $f\in\mmf{X}$ and $\alpha\in\Omega$. Fix
$\alpha\in\Omega$ and $\epsilon>0$. If $f\in\EuScript{D}$,then we can choose a sequence $\{f_s\}_{s>0}$,
$f_s\in\mcl{D}(\mmm{Z}_s),$ such that
$$\lim\limits_{s\to 0}{\mmm{Z}_s f_s}=\mmm{Z}_0 f.$$
So, we have
\begin{eqnarray}
&&\label{fg:2}\|\mmm{T}_{s_k}(l)f-f\|_{\alpha}\leq
\|\mmm{T}_{s_k}(l)f_{s_k}-f_{s_k}\|_{\alpha}+
\|f-f_{s_k}\|_{\alpha}\\&+&\|\mmm{T}_{s_k}(l)f_{s_k}-\mmm{T}_{s_k}(l)f\|_{\alpha}\leq
\|\int_0^l{\mmm{T}_{s_k}(r)\mmm{Z}_{s_k}f_{s_k}\,dr}\|_{\alpha}\\&+&
\|f-f_{s_k}\|_{\alpha}+\|\mmm{T}_{s_k}(l)f_{s_k}-\mmm{T}_{s_k}(l)f\|_{\alpha}\\&\leq&
l\|\mmm{Z}_{s_k}f_{s_k}\|_{\alpha}^{\mmf{A}_{l}}+
2\|f-f_{s_k}\|_{\alpha}^{\mmf{A}_{l}}\label{fg:3}
\end{eqnarray}
Choose $k_0\in\mmb{N}$ such that
$$\|\mmm{Z}_{s_k}f_{s_k}-\mmm{Z}_{0}f_{0}\|_{\alpha}^{\mmf{A}_{l}}\leq\epsilon$$
and
$$\|f-f_{s_k}\|_{\alpha}^{\mmf{A}_{l}}<{\textstyle \frac{\epsilon}{4}}$$
for any natural $k\leq k_0$. Denote
$$m=\sup\limits_{k\leq k_0}
{(\|\mmm{Z}_{s_k}f_{s_k}\|_{\alpha}^{\mmf{A}_{l}}
+\|\mmm{Z}_{0}f_{0}\|_{\alpha}^{\mmf{A}_{l}})}+\epsilon,\,\,k\in\mmb{N},$$
and put $l_0=\frac{\epsilon}{2m}$. Then, by inequalities
\eqref{fg:2}-\eqref{fg:3}, we deduce that
$$\|\mmm{T}_{s_k}(l)f-f\|_{\alpha}<\epsilon$$ for any
$l\in[0,\,l_0]$. By the arbitrariness $\epsilon$ we get inequality
\eqref{fg:1} for each $f\in\EuScript{D}.$ Therefore, $\exp{(l\mmm{Z}_{s_k})}f$ converges to $f$ as $l\to 0$
uniformly with respect to $k\in\mmb{N}$ for each $f\in\EuScript{D}.$ By the density of $\EuScript{D}$ in $\mmf{X}$ we easily infer that the latter fact and equality \eqref{fg:1} are true for each
$f\in\mmf{X}$. Therefore, we can choose $l_1>0$ such that
\begin{eqnarray}\sup\limits_{k\in\mmb{N}}{\|\exp{(l\mmm{Z}_{s_k})}g-g\|
_{\alpha_0}^{\mmf{A}_{t_0}}}<\delta/3\label{fg:5}
\end{eqnarray}
for each $l\in[0,\,l_1].$ Fix arbitrary
$t\in(0,\,t_0].$ Define the function $\sigma:[0,\,\infty)\to\mmb{R}$
by the
following conditions:

\noindent i) $\sigma(s)=s_1,\,\,s\in[t,\,\infty).$

\noindent ii)
$\sigma(s)=s_{k+1},\,\,s\in[\frac{t}{k+1},\,\frac{t}{k}),\,\,k\in\mmb{N}.$

\noindent iii) $\sigma(0)=0.$

\noindent Further, define the function
$\mmm{F}:[0,\,\infty)\to\mcl{L}(\mmf{X})$ by the equality
$$\mmm{F}(s)f=\mmm{T}_{\sigma(s)}(s)f,\,\,s\geq 0,\,\,f\in\mmf{X}.$$
It is evidently that $\mmm{F}\in\EuScript{F}_{\mmf{X}}.$ Now let us show that $\mmm{F}\in\EuScript{F}_{\mmf{X}}$ and
$$\mmm{F}_{ef}'(0)f =\mmm{Z}f$$
for each $f\in\EuScript{D}.$ Fix arbitrary $f\in\EuScript{D}$.
Then there exists a sequence $\{f_s\}_{s>0}$,
$f_s\in\mcl{D}(\mmm{Z}_s),$ such that $\lim\limits_{s\to 0}{f_s}=f$
and $$\lim\limits_{s\to 0}{\mmm{Z}_s f_s}=\mmm{Z}_0 f.$$ Note that
\begin{eqnarray}
&&\lim\limits_{s\to
0}{\|s^{-1}(\mmm{F}(s)-\mmm{I})f_{\sigma(s)}-\mmm{Z}_0f\|_{\alpha}}
\nonumber\\&=&\lim\limits_{s\to
0}{\|s^{-1}(\mmm{T}_{\sigma(s)}(s)-\mmm{I})f_{\sigma(s)}-\mmm{Z}_0f\|_{\alpha}}\nonumber\\&=&\lim\limits_{s\to
0}{\|s^{-1}\int_{0}^{s}{\mmm{T}_{\sigma(s)}(l)\mmm{Z}_{\sigma(s)}f_{\sigma(s)}\,dl}
-\mmm{Z}_{\sigma(s)}f_{\sigma(s)}}\nonumber\\&+&{\mmm{Z}_{\sigma(s)}f_{\sigma(s)}-\mmm{Z}_0f\|_{\alpha}}\nonumber
\end{eqnarray}
\begin{eqnarray}
&=&
\lim\limits_{s\to
0}{\|s^{-1}\int_{0}^{s}{(\mmm{T}_{\sigma(s)}(l)-\mmm{I})
\mmm{Z}_{\sigma(s)}f_{\sigma(s)}\,dl}\|_{\alpha}}
\nonumber\\&+&\lim\limits_{s\to
0}{\|\mmm{Z}_{\sigma(s)}f_{\sigma(s)}-\mmm{Z}_0f\|_{\alpha}}\nonumber\\&
=&\lim\limits_{s\to
0}{\|s^{-1}\int_{0}^{s}{(\mmm{T}_{\sigma(s)}(l)-\mmm{I})
(\mmm{Z}_{\sigma(s)}f_{\sigma(s)}-\mmm{Z}_0
f)\,dl}\|_{\alpha}}\nonumber\\&+&\lim\limits_{s\to
0}{\|s^{-1}\int_{0}^{s}{(\mmm{T}_{\sigma(s)}(l)-\mmm{I}) \mmm{Z}_0
f\,dl}\|_{\alpha}}\nonumber\\&=&\lim\limits_{s\to
0}{2\|\mmm{Z}_{\sigma(s)}f_{\sigma(s)}-\mmm{Z}_0 f
\|_{\alpha}^{\mmf{A}_{s}}}=0\nonumber
\end{eqnarray}
Here we have used equality \eqref{fg:1}. Thus,
$\mmm{F}\in\EuScript{F}_{\mmf{X}}$ and $\mmm{F}_{ef}'(0)f =\mmm{Z}f$
for each $f\in\EuScript{D}.$ Therefore, by Theorem~\ref{chern}, we infer that
$$\lim\limits_{n\to\infty}{\mmm{F}(t/n)^n
g}=\lim\limits_{n\to\infty}{\mmm{T}_{s_n}(t)g}=\mmm{T}_0(t)g$$ for arbitrary  $t\in(0,\,t_0].$ Therefore, by inequality
\eqref{fg:5} we have
\begin{eqnarray}\sup\limits_{k\in\mmb{N}}{\|\exp{(l\mmm{Z}_{0})}g-g\|
_{\alpha_0}^{\mmf{A}_{t_0}}}\leq\delta/3  \label{fg:7}
\end{eqnarray}
for each $l\in[0,\,l_1].$ Put $m_0=[t/l_1]$ and define
$t_i=i l_1$ for each natural $i\leq m_0$. Choose
$k_0\in\mmb{N}$ such that
$$\|\mmm{T}_{s_n}(t_i)g-\mmm{T}_0(t_i)g\|_{\alpha}<\delta/3$$
for each natural numbers $i\leq m_0$ and $n\geq k_0$. Fix arbitrary $s\in[0,\,t_0]$. Choose natural $k\leq m_0$
such that $|s-t_k|\leq l_1$. Then, by inequalities \eqref{fg:5} and \eqref{fg:7} we infer that
\begin{eqnarray}&&\|\mmm{T}_{s_n}(s)g-\mmm{T}_0(s)g\|_{\alpha}
\leq\|\mmm{T}_{s_n}(t_i)g-\mmm{T}_0(t_i)g\|_{\alpha}\nonumber\\&+&
\|\mmm{T}_{s_n}(t_i)g-\mmm{T}_{s_n}(s)g\|_{\alpha}+
\|\mmm{T}_{0}(t_i)g-\mmm{T}_{0}(s)g\|_{\alpha}\nonumber\\&\leq&
\delta/3+ \|\mmm{T}_{s_n}(|t_i-s|)g-g\|_{\alpha}^{\mmf{A}_{t_0}}+
\|\mmm{T}_{0}(|t_i-s|)g-g\|_{\alpha}^{\mmf{A}_{t_0}}<\delta\nonumber
\end{eqnarray}
for each natural $n\geq k_0$. Then, by arbitrariness
$s\in[0,\,t_0]$, it follows that
\begin{eqnarray}\liminf\limits_{k\to\infty}{ \sup\limits_{s\in[0,\,t_0]}
{\|\mmm{T}_{s_k}(s)g-\mmm{T}_{0}(s)g\|_{\alpha_0}}}\leq\delta.\nonumber
\end{eqnarray}
Thus, we get contradiction with the initial assumption.
  \end{proof}

The following theorem is the extension of Lie-Trotter theorem.

\begin{theorem}\label{t:22}
Let $\mmm{A}$, $\mmm{B}$ and $\mmm{Z}$ be generators of the
$lC_0$-semigroups in sequentially complete
locally convex space  $\mmf{X}$ and $\EuScript{D}$ be a core of the generator $\mmm{Z}$.
Suppose that $\EuScript{D}\subset\mcl{D}(\mmm{A})\cap\mcl{D}(\mmm{B})$,
$$ \mmm{Z}f=\mmm{A}f+\mmm{B}f$$
for all $f\in\EuScript{D}$ and the function
$$ \{[0,\,\infty)\ni s\to\exp{(s\mmm{A})}\exp{(s\mmm{B})}\}\in\EuScript{E}_{\mmf{X}}.$$
Then
$(\exp{({\textstyle\frac{t}{n}}\mmm{A})}\exp{({\textstyle\frac{t}{n}}\mmm{B})})^n
f$ tends to $\exp{(t\mmm{Z})}f$ as $n\to\infty$ for all
$f\in\mmf{X}$ uniformly with respect to $t\in[0,\,t_0]$ for any $t_0>0$.
\end{theorem}
\begin{proof}[Proof of Theorem~\ref{t:22}.] Define the function
$\mmm{F}(\cdot)$ by the equality
$$\mmm{F}(t)=\exp{(t\mmm{A})}\exp{(t\mmm{B})},\,t\in [0,\,\infty).$$
Then $\mmm{F}(0)=\mmm{I}$. Moreover, if $f\in \EuScript{D}$, then
\begin{eqnarray}
\lim\limits_{t\to 0}{t^{-1}(\mmm{F}(t)f-f)}&=&\lim\limits_{t\to
0}{(\exp{(\mmm{A}t)}t^{-1}(\exp{(\mmm{B}t)}f-f))}\nonumber\\&+&
\lim\limits_{t\to
0}{t^{-1}(\exp{(\mmm{A}t)}f-f)}=\mmm{B}f+\mmm{A}f.\nonumber
\end{eqnarray}
Thus, we can apply Theorem~\ref{chern}.
  \end{proof}

In conclusion, let us show that results of \cite{Albanese0},
\cite{Albanese} are not applicable to arbitrary equicontinuous
semigroups in sequentially complete locally convex spaces. Indeed,
it suffices to construct a generator $\mmm{Z}$ of the
$lC_0$-semigroups $\mmm{T}$, such that there doesn't exist
 $\lambda\in\mmb{C}$, for which the image of the operator
 $\lambda\mmm{I}-\mmm{Z}$ is dense in $\mmm{X}$. The following example presents such generator.
 \begin{example}\label{exa}
 Let $\mmf{X}=C(\mmb{C},\mmb{C})$ be the local convex space with topology generated by semi-norms $\|\cdot\|_r, \,\,r>0,$, defined by the equality $$\|f\|_r=\sup\limits_{|x|\leq r}{|f(x)|},\,\,x\in\mmb{C},$$ for each $f\in\mmf{X}$. Define the multiplication operator $\mmm{Z}:\mmf{X}\to\mmf{X}$ by the equality $\mmm{Z}(f)(x)=xf(x),\,\,f\in\mmf{X},\,\,x\in\mmb{C}.$
 Then $\mmm{Z}$ generates $lC_0$-semigroup $\mmm{T}$, defined by
 $$\mmm{T}(s)(f)(x)=e^{sx}f(x),\,\,s\geq 0,\,\,f\in\mmf{X},\,\,x\in\mmb{C}.$$ Note that for any $\lambda\in\mmb{C}$ the closure of the image of the operator
 $\lambda\mmm{I}-\mmm{Z}$ doesn't contain any $f\in\mmf{X}$ for which
 $f(\lambda)\neq 0$. Thus, the image of the operator
 $\lambda\mmm{I}-\mmm{Z}$ isn't dense in $\mmm{X}$ for any $\lambda\in\mmb{C}$.
 \end{example}

\section{Necessary and sufficient conditions of existence of local solution
for equation $\dot{x}=\mmm{A}x$ where
$\mmm{A}\in\mcl{Z}_{\mmf{X}}$.}\label{3.1}

\begin{lemma}\label{t:10b} Let $\mathrm{F}\in\EuScript{F}_{\mmf{X}}$ and $\mmf{X}$  be a sequentially complete locally convex space.  Assume that for some $t,\,\,l>0$ and
subspace $\Phi\subset\mmf{X}$ there exists a strictly increasing natural sequence $\{f_k\}_{k=1}^{\infty}$ such that there exists the limit
\begin{equation}
w\mbox{-}\lim\limits_{k
\rightarrow\infty}{\{\mathrm{F}({\textstyle\frac{t}{g_k}})\}^{[{g_k}\frac{s}{t}]}
g,} \label{eqn:m1}
\end{equation}
for any $g\in \Phi$ and rational $s\in[0,\,l)$.
  Then limit~\eqref{eqn:m1} exists for any  $g\in \Phi$,
   $s\in[0,\,l)$.
\end{lemma}
\begin{proof}[Proof of Lemma~\ref{t:10b}.]
The fact that the sequence
$\{\{\mathrm{F}({\textstyle\frac{t}{g_k}})\}^{[{g_k}\frac{s}{t}]}g\}_{n\in\mmb{N}}$
is a Cauchy sequence in the topology $\sigma (\mmf{X},\mmf{X}^*)$ for any $g\in\Phi$ and  $s\in[0,\,l)$ follows from the condition $\mathrm{F}\in\EuScript{F}_{\mmf{X}}$.
Thus, by Lemma~\ref{t:7}, sequential completeness of the space
$\mmf{X}$ and the condition $\mathrm{F}\in\EuScript{F}_{\mmf{X}}$,
we get the Lemma.
\end{proof}
\begin{lemma}\label{t:10}
Let $\mathrm{F}\in\EuScript{F}_{\mmf{X}}$ and $\mmf{X}$ be a locally convex space.  Assume that for some $t,\,\,l>0$ and separable closed linear subspace
 $\Phi\subset\mmf{X}$ there exists a strictly increasing natural sequence $\{f_k\}_{k=1}^{\infty}$ such that there exists the limit
\begin{equation}
w\mbox{-}\lim\limits_{k
\rightarrow\infty}{\{\mathrm{F}({\textstyle\frac{t}{g_k}})\}^{[{g_k}\frac{s}{t}]}
g,}
\end{equation}
for any $g\in \Phi$ and  $s\in[0,\,l)$. Then
  the family of the operators $\Phi_s: \Phi\mapsto
  \mmf{X}$, $s\in[0,\,l)$, defined by the equality
\begin{equation}
\Psi_s g=w\mbox{-}\lim\limits_{k
\rightarrow\infty}{\{\mathrm{F}({\textstyle\frac{t}{g_k}})\}^{[{g_k}
\frac{s}{t}]} g}
\end{equation}
  satisfies the following conditions:

\noindent a) $\Psi_s g$ is linear with respect to $g$ for each
$s\in[0,\,l)$;

\noindent b) $\Psi_s g$ is continuous with respect to
$s\in[0,\,l)$ for each $g\in \Phi$;

\noindent c) If $f\in \Phi\cap\mcl{D}(\ov{\mmm{F}'_{ef}(0)})$, then the sequence
$\{\mathrm{F}({\textstyle\frac{t}{g_k}})\}^{[{g_k}
\frac{s}{t}]}\ov{\mmm{F}'_{ef}(0)}f$, $s\in[0,\,l)$, is a Cauchy sequence in the topology $\sigma(\mmf{X},\mmf{X}^*)$ and there exists
$$(\Psi_{s} f, \phi)^{'}_{s}= \lim\limits_{k
\rightarrow\infty}{(\{\mathrm{F}({\textstyle\frac{t}{g_k}})\}^{[{g_k}
\frac{s}{t}]}\ov{\mathrm{F}'_{ef}(0)}f,\phi)},\,\,s\in[0,\,l),\,\,\phi\in\mmf{X}^*;$$

\noindent d) If $\phi\in\mcl{D}({\mmm{F}^*}'(0))$, then
$$(\Psi_{m}f, \phi)- (\Psi_{p}f,
\phi)=\int_{p}^{m}{(\Psi_{s}f, {\mmm{F}^*}'(0)\phi)\,d s},\,\,
m,\,\,p\in[0,\,l),$$ for any $f\in \Phi$;

\noindent e) If $\phi\in\mcl{D}({\mmm{F}^*}'(0))$, then $$
(\Psi_{s}f, \phi)^{'}_{s}=(\Psi_{s}f, {\mmm{F}^*}'(0)\phi)$$ for any $f\in \Phi;$

\noindent f) Let $f\in \Phi\cap\mcl{D}(\ov{\mmm{F}'_{ef}(0)})$,
$\Psi_{s}f\in\mcl{D}(\ov{\mmm{F}'_{ef}(0)})$ for any
$s\in[0,\,l)$ and
 $\mcl{D}({\mmm{F}^*}'(0))$ be *-dense in $\mmf{X}^*$. If there
 exists
$$w\mbox{-}\lim\limits_{k
\rightarrow\infty}{\{\mathrm{F}({\textstyle\frac{t}{g_k}})\}^{[{g_k}
\frac{s}{t}]}\ov{\mmm{F}'_{ef}(0)}f}=q(s)\in\tilde{\mmf{X}},\,\,s\in[0,\,l),$$ where
$\tilde{\mmf{X}}$ is the completion of the space $\mmf{X}$, then
the following equality holds
$$(\Psi_{s}f)^{'}_{s}=\ov{\mmm{F}'_{ef}(0)}\Psi_{s}f=q(s)$$
for any $s\in[0,\,l)$;

\noindent g) Let $f\in \Phi\cap\mcl{D}(\mmm{F}'_{ef}(0))$,
$\Psi_{s}f\in\mcl{D}(\ov{\mmm{F}'_{ef}(0)})$ for each
$s\in[0,\,l)$ and
 $\mcl{D}({\mmm{F}^*}'(0))$ be *-dense in $\mmf{X}^*$. Then
the following equality holds
$$(\Psi_{s}f)^{'}_{s}=\ov{\mmm{F}'_{ef}(0)}\Psi_{s}f$$
for each $s\in[0,\,l)$.
\end{lemma}
\begin{proof}[Proof of Lemma~\ref{t:10}.]
 Part (a)
is trivial. Part (b) follows from part  (b)  of Proposition~\ref{t:1} and Lemma~\ref{t:7}.
Let us show part (c). Fix arbitrary subsequence
$\{g'_k\}_{k=1}^{\infty}$ of the sequence
$\{g_k\}_{k=1}^{\infty}$. By Proposition~\ref{t:1}, choose subsequence $\{f_k\}_{k=1}^{\infty}$ of the sequence
$\{g'_k\}_{k=1}^{\infty}$  such that we can define the set of the functions
$\mathrm{T}_{s}:\Phi'\times\{\phi\}\mapsto\mmb{T},\,\,s\geq 0,$
 by the inequality
$$\mathrm{T}_{s}(g,\phi)=\lim\limits_{k
\rightarrow\infty}{(\{\mathrm{F}(t/{f_k})\}^{[{f_k}\frac{s}{t}]}
y,\phi)}, \,\,y\in\Phi',$$ where
$\Phi'=\span1{\{f,\ov{\mmm{F}'_{ef}(0)}f\}}.$ Then, by Lemma~
\ref{n:2}, we infer that $$(\mathrm{T}_{s}
(f,\phi))^{'}_{s}=\mathrm{T}_{s}(\overline{\mmm{F}_{ef}'(0)}f,\phi),\,\,s\geq
0.$$ Since $$\mathrm{T}_{s} (f,\phi)=(\Psi_s
f,\phi),\,\,s\in[0,\,l),$$  and $\{g'_k\}_{k=1}^{\infty}$ is
arbitrary we infer that the subsequence
$\{\mathrm{F}({\textstyle\frac{t}{g_k}})\}^{[{g_k}
\frac{s}{t}]}\ov{\mmm{F}'_{ef}(0)}f$, $s\in[0,\,l)$, is a Cauchy
sequence in the topology $\sigma(\mmf{X},\mmf{X}^*)$ and
$$(\Psi_{s} f, \phi)^{'}_{s}= \lim\limits_{k
\rightarrow\infty}{(\{\mathrm{F}({\textstyle\frac{t}{g_k}})\}^{[{g_k}
\frac{s}{t}]}\ov{\mathrm{F}'_{ef}(0)}f,\phi)},\,\,s\in[0,\,l),\,\,\phi\in\mmf{X}^*.$$
 Let us show part (d). Assume that
$\phi\in\mcl{D}({\mmm{F}^*}'(0))$. Fix $\epsilon>0$ and
$m,\,\,l\in\mmb{T}$ such that $m>l\geq 0$. Choose $r_0$, $j_0$
such that
\begin{eqnarray}(({\mathrm{F}({\textstyle\frac{t}{g_k}})}^{[{g_k}\frac{g}{t}]}-
{\mathrm{F}({\textstyle\frac{t}{g_k}})}^{[{g_k}\frac{h}{t}]})f,{\mathrm{F}^{*}}'(0)\phi)
<\epsilon, \,\,k\geq j_0, \,\,|g-h|<r_0,\\
g,h\in[0,\,l).\end{eqnarray} The existence of such $r_0$, $j_0$
follows from Lemma~\ref{t:7}. Hence, we have
\begin{eqnarray}|\Psi_h(f,{\mathrm{F}^{*}}'(0)\phi)-\Psi_g(f,{\mathrm{F}^{*}}'(0)\phi)|\leq\epsilon
,\label{q:1}
\end{eqnarray} if $|g-h|<r_0,$ $g,h\in[0,\,l)$.
Therefore, we get
\begin{eqnarray}
&&|\sum_{v=0}^{v_0-1}{\Psi_{l+v(m-l)/v_0}(f,
{\mmm{F}^*}'(0)\phi)((m-l)/v_0)}\nonumber\\&-&
\int_{l}^{m}{\Psi_{s}(f, {\mmm{F}^*}'(0)\phi)\,d s}|\leq \epsilon
|m-l|,\label{eqn:d0}
\end{eqnarray}
$v_0=[(m-l)/r_0]+1$.
Choose $j_1>j_0$ such that
\begin{eqnarray}
|\Psi_{l+v(m-l)/v_0}(f, {\mmm{F}^*}'(0)\phi)-
({\mathrm{F}({\textstyle\frac{t}{g_k}})}^{[{g_k}(l+v(m-l)/v_0)/t]}f,{\mathrm{F}^{*}}'(0)\phi)|<\epsilon\nonumber
\end{eqnarray}
for any $v\in\{0,1,\dots ,v_0\}$, $k>j_1$. Thus,
\begin{eqnarray}
&\,&|((m-l)/v_0)\sum_{v=0}^{v_0-1}{({\mathrm{F}({\textstyle\frac{t}{g_k}})}
^{[{g_k}(l+v(m-l)/v_0)/t]}f,{\mathrm{F}^{*}}'(0)\phi)}\label{eqn:d1}\\
&-&\sum_{v=0}^{v_0-1}{\Psi_{l+v(m-l)/v_0}(f,
{\mmm{F}^*}'(0)\phi)((m-l)/v_0)}|<\epsilon |m-l|.\label{eqn:d2}
\end{eqnarray}
Note that
\begin{eqnarray}
&\,&\label{q:2}\Psi_{m}(f,\phi)-\Psi_{l}(f,\phi)\\&=&\lim\limits_{k
\rightarrow\infty}{({\mathrm{F}({\textstyle\frac{t}{g_k}})}^{[{g_k}\frac{m}{t}]}
f-{\mathrm{F}({\textstyle\frac{t}{g_k}})}^{[{g_k}\frac{l}{t}]} f,\phi)}\\
&=&\lim\limits_{k
\rightarrow\infty}{((\mathrm{F}({\textstyle\frac{t}{g_k}})
-\mathrm{I})
\sum_{i=[{g_k}\frac{l}{t}]}^{[{g_k}\frac{m}{t}]-1}{{\mathrm{F}({\textstyle\frac{t}{g_k}})}^{i}f},\phi)}\label{q:3}\\
&=&\lim\limits_{k \rightarrow\infty}{
({\textstyle\frac{t}{g_k}}\sum_{i=[{g_k}\frac{l}{t}]}^{[{g_k}\frac{m}{t}]-1}{{\mathrm{F}({\textstyle\frac{t}{g_k}})}^{i}f},
{\textstyle\frac{g_k}{t}}(\mathrm{F}^{*}({\textstyle\frac{t}{g_k}}) -\mathrm{I})\phi)}\label{q:4}\\
&=& \lim\limits_{k \rightarrow\infty}{
({\textstyle\frac{t}{g_k}}\sum_{i=[{g_k}\frac{l}{t}]}^{[{g_k}\frac{m}{t}]-1}{{\mathrm{F}({\textstyle\frac{t}{g_k}})}^{i}f},
{\textstyle\frac{g_k}{t}}(\mathrm{F}^{*}({\textstyle\frac{t}{g_k}})
-\mathrm{I})\phi-{\mathrm{F}^{*}}'(0)\phi)}\,\,\,\,\,\,\,\,\,\label{q:5}\\
&+&\lim\limits_{k \rightarrow\infty}{
({\textstyle\frac{t}{g_k}}\sum_{i=[{g_k}\frac{l}{t}]}^{[{g_k}\frac{m}{t}]-1}{{\mathrm{F}({\textstyle\frac{t}{g_k}})}^{i}f},
{\mathrm{F}^{*}}'(0)\phi)}\label{q:6}\\
&=&\lim\limits_{k \rightarrow\infty}{
({\textstyle\frac{t}{g_k}}\sum_{i=[{g_k}\frac{l}{t}]}^{[{g_k}\frac{m}{t}]-1}{{\mathrm{F}({\textstyle\frac{t}{g_k}})}^{i}f},
{\mathrm{F}^{*}}'(0)\phi)}\label{q:7}\\
&=&\lim\limits_{k
\rightarrow\infty}{{\textstyle\frac{t}{g_k}}\sum_{s=0}^{d_k}{({\mathrm{F}({\textstyle\frac{t}{g_k}})}^{[{g_k}(l+s({\textstyle\frac{t}{g_k}}))/t]}f,{\mathrm{F}^{*}}'(0)\phi)}}\label{q:8}\\
&=&\lim\limits_{k
\rightarrow\infty}{\sum_{v=0}^{v_0-1}{\sum_{s\in\mmf{B}_{k,v,v_0}}{{\textstyle\frac{t}{g_k}}({\mathrm{F}({\textstyle\frac{t}{g_k}})}^{[{g_k}(l+s({\textstyle\frac{t}{g_k}}))/t]}f,{\mathrm{F}^{*}}'(0)\phi)}}},\label{q:9}
\end{eqnarray}
where $d_k=[(m-l){\textstyle\frac{g_k}{t}}]$ and
$$\mmf{B}_{k,v,v_0}=\{s\in\mmb{N}|\,v(m-l)/v_0\leq s {\textstyle\frac{t}{g_k}} <
(v+1)(m-l)/v_0\}.$$

\noindent So, by inequality~\eqref{q:1} and inequalities
\eqref{q:2}
--- \eqref{q:9}, we get
\begin{eqnarray}
&&|\label{eqn:d3}
\sum_{v=0}^{v_0-1}{{\textstyle({\mathrm{F}({\textstyle\frac{t}{g_k}})}
^{[{g_k}(l+v
\frac{m-l}{v_0})/t]}f,{\mathrm{F}^{*}}'(0)\phi)}({\textstyle
\frac{m-l}{v_0}})}\\&+& \Psi_{l}(f,\phi)-\Psi_{m}(f,\phi)|
\\&=&|\lim\limits_{k
\rightarrow\infty}{\sum_{v=0}^{v_0-1}{\sum_{s\in\mmf{B}_
{k,v,v_0}}{{\textstyle\frac{t}{g_k}}({\mathrm{F}({\textstyle\frac{t}{g_k}})}^{[{g_k}(l+s({\textstyle\frac{t}{g_k}}))/t]}f,{\mathrm{F}^{*}}'(0)\phi)}}}\\
&-&((m-l)/v_0)\sum_{v=0}^{v_0-1}{({\mathrm{F}({\textstyle\frac{t}{g_k}})}
^{[{g_k}(l+v(m-l)/v_0)/t]}f,{\mathrm{F}^{*}}'(0)\phi)}|\,\,\,\,\,\,\\
&\leq&\sum_{v=0}^{v_0-1}|\lim\limits_{k
\rightarrow\infty}{\sum_{s\in\mmf{B}_{k,v,v_0}}{{\textstyle\frac{t}{g_k}}({\mathrm{F}({\textstyle\frac{t}{g_k}})}^{[g_k(l+s({\textstyle\frac{t}{g_k}}))/t]}f,{\mathrm{F}^{*}}'(0)\phi)}}\\
&-&({\textstyle
\frac{m-l}{v_0}}(\mathrm{F}({\textstyle\frac{t}{g_k}})^{[g_k(l+v(m-l)/v_0)/t]}f,{\mathrm{F}^*}'(0)\phi)|<\epsilon
v_0 {\textstyle \frac{|m-l|}{v_0}}.\label{eqn:d4}
\end{eqnarray}
Combining inequalities \eqref{eqn:d0}, \eqref{eqn:d1} ---
\eqref{eqn:d2} and \eqref{eqn:d3}
--- \eqref{eqn:d4}, we infer that
\begin{eqnarray}
|\Psi_{m}(f,\phi)-\Psi_{l}(f,\phi)-\int_{l}^{m}{\Psi_{s}(f,
{\mmm{F}^*}'(0)\phi)\,d s}|<3\epsilon |m-l|\nonumber
\end{eqnarray}
So, part $(d)$ of Lemma~\ref{t:10} follows from the
arbitrariness of
 $\epsilon$. Part $(e)$ is a direct
consequence of parts $(b)$ and
$(d)$. Let us prove part $(f)$. Firstly, let us verify the equality
$q(s)=\ov{\mmm{F}'_{ef}(0)}\Psi_{s}f$, $s\in[0,\,l).$
 Note that for
$h\in\mcl{D}(\mmm{F}_{ef}'(0))$ there exists a sequence
$\{h_s\}_{s\in[0,\,\infty)}$ such that $\lim\limits_{s\to
0}{h_s}=h$ and
$$\lim\limits_{s\to 0}{s^{-1}(\mmm{F}(s)-\mmm{I})h_s}=\mmm{F}_{ef}'(0)h.$$
Therefore, we get
\begin{eqnarray}
(\mmm{F}_{ef}'(0)h,g)=\lim\limits_{s\to
0}{(s^{-1}(\mmm{F}(s)-\mmm{I})h_s,g)}\label{hg:1}\\=\lim\limits_{s\to
0}{(h_s,s^{-1}(\mmm{F}^{*}(s)-\mmm{I})g)}=(h,{\mmm{F}^*}'(0)g)\label{hg:2}
\end{eqnarray}
for each $g\in\mcl{D}({\mmm{F}^*}'(0)).$ Fix
arbitrary $w\in\mcl{D}(\ov{\mmm{F}_{ef}'(0)})$ and
$g\in\mcl{D}({\mmm{F}^*}'(0)).$ Choose the sequence
$\{w_k\}_{k=1}^{\infty}$, $w_k\in\mcl{D}(\mmm{F}_{ef}'(0))$, such that there exist
$$\lim\limits_{k\to\infty}{(w_k,u)}=(w,u),$$
$$\lim\limits_{k\to\infty}{(\mmm{F}_{ef}'(0)w_k,u)}=(\ov{\mmm{F}_{ef}'(0)}w,u)$$
for any $u\in\{g,\,{\mmm{F}^*}'(0)g\}$. Then, by equalities
\eqref{hg:1}-\eqref{hg:2}, we get
\begin{eqnarray}(\ov{\mmm{F}_{ef}'(0)}w,g)=\lim\limits_{k\to\infty}{(\mmm{F}_{ef}'(0)w_k,g)}\nonumber\\=
\lim\limits_{k\to\infty}{(w_k,{\mmm{F}^*}'(0)g)}=(w,{\mmm{F}^*}'(0)g).\nonumber
\end{eqnarray}
So, we infer that
\begin{eqnarray}(\ov{\mmm{F}_{ef}'(0)}w,g)=(w,{\mmm{F}^*}'(0)g)\nonumber
\end{eqnarray}
for any $w\in\mcl{D}(\ov{\mmm{F}_{ef}'(0)})$ and
$g\in\mcl{D}({\mmm{F}^*}'(0)).$ Therefore, by parts (c) and (e)
we get \begin{eqnarray}&&\nonumber(\Psi_{s} f,
\phi)^{'}_{s}= \lim\limits_{k
\rightarrow\infty}{(\{\mathrm{F}({\textstyle\frac{t}{g_k}})\}^{[{g_k}
\frac{s}{t}]}\ov{\mathrm{F}'_{ef}(0)}f,\phi)}\\&=&(q(s),\phi)=(\Psi_{s}f,{\mmm{F}^*}'(0)\phi)\nonumber\\&=&(\ov{\mmm{F}_{ef}'(0)}\Psi_{s}f,\phi),
\,\,s\in[0,\,l),\,\,\phi\in\mcl{D}({\mmm{F}^*}'(0)).\nonumber
\end{eqnarray} For arbitrary $\phi\in\mmf{X}^*$ choose
the sequence $\{\phi_k\}_{k=1}^{\infty}$,  $\phi_k\in\mcl{D}({\mmm{F}^*}'(0)),$ such that there exists
\begin{eqnarray}\lim\limits_{k\to\infty}{(x,\phi_k)}=(x,\phi)\label{fg:10}
\end{eqnarray} for each $x\in\tilde{\mmf{X}}$.
 Therefore, we get
$$(\ov{\mmm{F}_{ef}'(0)}\Psi_{s}f,\phi)=\lim\limits_{k\to\infty}{(\ov{\mmm{F}_{ef}'(0)}\Psi_{s}f,\phi_k)}=\lim\limits_{k\to\infty}{(q(s),\phi_k)}=(q(s),\phi),$$
for each $\phi\in\mmf{X}^*$. Thus,
$$\ov{\mmm{F}_{ef}'(0)}\Psi_{s}f=q(s)=w\mbox{-}\lim\limits_{k
\rightarrow\infty}{\{\mathrm{F}({\textstyle\frac{t}{g_k}})\}^{[{g_k}
\frac{s}{t}]}\ov{\mmm{F}'_{ef}(0)}f}.$$ So, by Lemma~\ref{t:7},
we deduce that the function $q(\cdot)$ is continuous in $\mmf{X}$.
Hence, we have
\begin{eqnarray}(\Psi_{m}f-\Psi_{p}f,\phi)=(\Psi_{m}f,
\phi)- (\Psi_{p}f, \phi)\\=\int_{p}^{m}{(q(s), \phi)\,d
s}=(\int_{p}^{m}{q(s) \,d s}, \phi),\,\, m,\,\,p\in[0,\,l),
\end{eqnarray}
for any $\phi\in\mmf{X}^*$, where $$\int_{p}^{m}{q(s) \,d
s}\in\tilde{\mmf{X}}.$$ By the arbitrariness of $\phi$, we get the equality
$$\Psi_{m}f-\Psi_{p}f=\int_{p}^{m}{q(s) \,d s},\,\, m,\,\,p\in[0,\,l).$$
Thus, by continuity of the function $q(\cdot)$, we infer that
$$(\Psi_{s}f)^{'}_{s}=q(s)=\ov{\mmm{F}'_{ef}(0)}\Psi_{s}f$$ for any $s\in[0,\,l).$ Let us show part (g).
Note that it follows from part (c)  that the sequence $\{{\mathrm{F}({\textstyle\frac{t}{g_k}})}^{[{g_k}
\frac{s}{t}]}\mmm{F}'_{ef}(0)f\}_{k\in\mmb{N}}$, $s\in[0,\,l)$, is
a Cauchy sequence in the topology $\sigma(\mmf{X},\mmf{X}^*)$. It suffices to prove that there exists
$$w\mbox{-}\lim\limits_{k
\rightarrow\infty}{\{\mathrm{F}({\textstyle\frac{t}{g_k}})\}^{[{g_k}
\frac{s}{t}]}\mmm{F}'_{ef}(0)f}=q(s)\in\tilde{\mmf{X}},$$ where
$\tilde{\mmf{X}}$ is the completion of $\mmf{X}$. Indeed, we can easily deduce part (g) from the latter fact and part (f).
Below we will figure that the set
$\{\|\cdot\|_{\alpha}|\alpha\in\Omega\}$ consists of all continuous semi-norm of the space $\mmf{X}$.
Let us show the fact that if a sequence $\{x_k\}_{k\in\mmb{N}}$
is a Cauchy sequence in the topology $\sigma(\mmf{X},\mmf{X}^*)$ and for
each $\alpha\in\Omega$ there exists $\{y_k\}_{k\in\mmb{N}}$ such that there exists $w\mbox{-}\lim_{n\to\infty}{y_n}$ and $\limsup_{n\to\infty}{\|x_n -y_n\|_{\alpha}}\leq 1$, then there exists $w\mbox{-}\lim_{n\to\infty}{x_n}\in\tilde{\mmf{X}}$. Define
$\mcl{B}=\{\mmf{B}_{\alpha}|\,
\alpha\in\Omega\}$, where $\mmf{B}_{\alpha}=\{y\in\mmf{X}|\ex\{y_k\}_{k\in\mmb{N}}\,\, \mbox{such that}\,\, \limsup_{n\to\infty}{\|x_n -y_n\|_{\alpha}}\leq 1
\,\,\mbox{and}\,\, y=w\mbox{-}\lim_{n\to\infty}{y_n}\}$. Then $\mcl{B}$ is a filter base. Indeed, $\mcl{B}$ is not empty, doesn't contain
empty set and for each $\alpha_1,\alpha_2\in\Omega$ there exists $\alpha_3\in\Omega$ such that $\mmf{B}_{\alpha_3}\subset\mmf{B}_{\alpha_1}\cap\mmf{B}_{\alpha_2}$, where
$\alpha_3$ can be defined by the equality $\|\cdot\|_{\alpha_3}=\|\cdot\|_{\alpha_1}+\|\cdot\|_{\alpha_2}$. Moreover,
 $\mcl{B}$ is a Cauchy filter by construction. Thus, there exists $q\in\tilde{\mmf{X}}$, which the filter $\mcl{B}$ converges to. Hence, we easily deduce that there exists $w\mbox{-}\lim_{n\to\infty}{x_n}=q$. Now it suffices to prove that for each   $\alpha\in\Omega$ there exists a sequence $\{y_k\}_{k\in\mmb{N}}$ such that there exist $w\mbox{-}\lim_{n\to\infty}{y_n}$ and
\begin{eqnarray}
\limsup_{n\to\infty}{\|\{\mathrm{F}({\textstyle\frac{t}{g_k}})\}^{[{g_k}
\frac{s}{t}]}\mmm{F}'_{ef}(0)f -y_n\|_{\alpha}}\leq 1.\label{eqn41}
\end{eqnarray}
Fix arbitrary $\alpha\in\Omega$. Let $\{f_s\}_{s\geq 0}$ be a sequence
such that $\lim\limits_{s\to 0}{f_s}=f$ and
$\lim\limits_{s\to 0}{s^{-1}(\mmm{F}(s)-\mmm{F}(0))f_s}=\mmm{F}'_{ef}(0)f$.
By Lemma~\ref{t:8}, we infer that there exist
$v_0>0$, $k_0\in\mmb{N}$ such that
\begin{eqnarray}
\nonumber& &
\|{\mathrm{F}}^{[{g_k}\frac{s}{t}]}({\textstyle\frac{t}{g_k}})f_{{\textstyle\frac{t}{g_k}}}-
 {\mathrm{F}}^{
[{g_k}\frac{s+v}{t}]}({\textstyle\frac{t}{g_k}})f_{{\textstyle\frac{t}{g_k}}}-{\textstyle\frac{t}{g_k}}([{g_k}{\textstyle\frac{s}{t}}]-
[{g_k}{\textstyle\frac{s+v}{t}}])\\ &\times&
{\mathrm{F}}^{[{g_k}\frac{s}{t}]}({\textstyle\frac{t}{g_k}})\mmm{F}'_{ef}(0)f\|_{\alpha}\nonumber
< |[{g_k}{\textstyle\frac{s}{t}}]-[{g_k}{\textstyle\frac{s+v}{t}}]|{\textstyle\frac{t}{g_k}},\,\,k\in\mmb{N},
\end{eqnarray}
for any $v\in[0,\,2v_0)$ and $k>k_0$. Therefore, the sequence $\{y_k\}_{k\in\mmb{N}}$ defined by the inequalities
$$y_k=([{g_k}{\textstyle\frac{s}{t}}]-[{g_k}{\textstyle\frac{s+v}{t}}]|{\textstyle\frac{t}{g_k}})^{-1}({\mathrm{F}}^{[{g_k}\frac{s}{t}]}({\textstyle\frac{t}{g_k}})f_{{\textstyle\frac{t}{g_k}}}-
 {\mathrm{F}}^{
[{g_k}\frac{s+v}{t}]}({\textstyle\frac{t}{g_k}})f_{{\textstyle\frac{t}{g_k}}}),\,\,k>k_0,$$
$$y_k=0,\,\,k\leq k_0, $$
holds inequality~\eqref{eqn41} and there exists $w\mbox{-}\lim_{n\to\infty}{y_n}=(\Psi_{s+v}-\Psi_{s})/v.$
By the arbitrariness of $\alpha\in\Omega$ we get part (g).
\end{proof}

\begin{theorem}\label{vvmain:10a} \label{vvmain:10b} Let $\mmm{Z}\in\mcl{Z}_{\mmf{X}}$ and $t>0$. Assume that the function
$\mathrm{F}\in\EuScript{F}_{\mmf{X}}$ satisfies the following conditions:

\noindent i) $\mcl{D}(\mmm{F}_{ef}'(0))\supset\mcl{D}(\mmm{Z})$ and
$\mmm{F}_{ef}'(0)g=\mmm{Z}g,\,\,g\in\mcl{D}(\mmm{Z})$;

\noindent ii) $\mcl{D}({\mmm{F}^*}'(0))$  is *-dense in $\mathbf{X}^*$.

\noindent Then the local solution $f:[0,\,l)\mapsto\mmm{X},\,\,l>0,$
of the system
\begin{eqnarray}
f'(v)&=&\overline{\mmm{Z}}f(v),\,v\in[0,\,l)\nonumber\\
f(0)&=&h\in\mcl{D}(\mmm{Z}),\nonumber
\end{eqnarray}
exists on the semi-interval $[0,\,l)$ iff the following conditions
are satisfied:

\noindent a) It is possible to choose a weakly convergent
sub-subsequence for any subsequence of the sequence
$\{\mathrm{F}^{[n s/t]}(\frac{t}{n})h\}_{n\in \Nat}$,
 for any $s\in[0,\,l)$;

\noindent b) For any $s\in[0,\,l)$  there exists
a sequence $\{f_n^{s}\}_{n=1}^{\infty}$,
$f_n^{s}\in\mcl{D}(\ov{\mmm{Z}}),$ such that there exists
$$w\mbox{-}\lim\limits_{n\to\infty}{(\mathrm{F}^{[n s/t]}({\textstyle
\frac{t}{n}})h-f_n^{s})}=0$$ and it is possible to choose a weakly
convergent sub-subsequence for any subsequence of the sequence
$\{\ov{\mmm{Z}}f_n^{s}\}_{n\in \Nat}$.

\noindent Moreover, if conditions (a)-(b) are satisfied, then
$$f(v)=\lim\limits_{n\to\infty}{\{{\mmm{F}(v/n)}\}^n h},\,v\in[0,\,l).$$
\end{theorem}
\begin{proof}[Proof of Theorem~\ref{vvmain:10a}.] ($\Rightarrow$) It follows from Corollaries~\ref{chern1b},~\ref{chern1d} and $f(v)\in\mcl{D}(\ov{\mmm{Z}})$ for any $v\in[0,\,l)$ (the last property is true because
the weak closure of
 $\mmm{Z}$ coincides with the closure in the strong topology).

\noindent ($\Leftarrow$) Let $\{s_k\}_{k=1}^{\infty}$ be a sequence of all nonnegative rational numbers contained in $[0,\,l)$.
  Fix a strictly increasing sequence  $\{n_i\}_{i=1}^{\infty}$,
of natural numbers. Choose a subsequence
$\{n_i^1\}_{i=1}^{\infty}$ of the sequence
$\{n_i\}_{i=1}^{\infty}$ such that there exists the limit
 $w\mbox{-}\lim\limits_{i\to\infty}{\mmm{F}^{[{n_i^1}
s_1/t]}(t/{n_i^1})h}$. Further, choose a subsequence
$\{n_i^2\}$ of the sequence $\{n_i^1\}$ such that there exists the limit
 $w\mbox{-}\lim\limits_{i\to\infty}{\mmm{F}^{[{n_i^2}
  s_2/t]}(t/{n_i^2})h}$. In the same way, for any natural $k\geq 3$, we choose the subsequence $\{n_i^k\}_{i=1}^{\infty}$
of the subsequence $\{n_i^{k-1}\}_{i=1}^{\infty}$ such that there exists the limit $w\mbox{-}\lim\limits_{i\to\infty}{\mmm{F}^{[{n_i^k}
s_k/t]}(t/{n_i^k})h}$. Note that for diagonal sequence $\{n_i^i\}_{i=1}^{\infty}$ there exists the limit $w\mbox{-}\lim\limits_{i\to\infty}{\mmm{F}^{[n_i^i
s_k/t]}(t/{n_i^i})h}$ for each $k\in\mmb{N}$. Thus, by Lemma~\ref{t:7}, we infer that
$\{\mmm{F}^{[n_i^i
  s/t]}(t/{n_i^i})h\}_{i=1}^{\infty}$ is a Cauchy sequence in the topology $\sigma(\mmf{X},\mmf{X}^{*})$
  for any  $s\in[0,\,l)$. Therefore, by condition
  (a), there exists
$$w\mbox{-}\lim\limits_{i\to\infty}{\mmm{F}^{[n_i^i
  s/t]}(t/{n_i^i})h}$$ for any  $s\in[0,\,l).$
  Denote $d_i=n_i^i,\,\,i\in\mmb{N},$ and
  $$\mmm{T}_s g=w\mbox{-}\lim\limits_{i\to\infty}{\mmm{F}^{[d_i
  s/t]}(t/{d_i})g}, \,\,g\in\span1{(\{h\})},\,\,s\in[0,\,l).$$ Fix arbitrary $m\in[0,\,l)$. Choose a subsequence
  $\{d'_k\}_{k\in\Nat}$ of the sequence $\{d_k\}_{k\in\Nat}$ such that there exists the limit $w\mbox{-}\lim\limits_{i\to\infty}{\ov{\mmm{Z}}f_{d'_i}^{m}}$. Thus, there exist
   $w\mbox{-}\lim\limits_{i\to\infty}{f_{d'_i}^{m}}$,
  $w\mbox{-}\lim\limits_{i\to\infty}{\ov{\mmm{Z}}f_{d'_i}^{m}}$ and, therefore,  from the coincidence of the weak closure of
  $\mmm{Z}$ with the closure in the strong topology it follows that
   $\mmm{T}_{m}h\in\mcl{D}(\overline{\mmm{Z}})$ for any
  $m\in[0,\,l)$. Then, from
  part (g) of Lemma~\ref{t:10}
  we infer that
  $$(\mathrm{T}_{s}h)^{'}_{s}=\overline{\mmm{Z}}\mathrm{T}_{s}h$$ for each $s\in[0,\,l)$.
  Consequently, $\mathrm{T}_{s}h$
is the solution of the equation $$f'(s)=\overline{\mmm{Z}}f(s)$$ with the initial condition $f(0)=h$ on the semi-interval $[0,\,l)$. Thus, by Theorem~\ref{chern1}, we have
$$f(s)=\lim\limits_{n\to\infty}{{\mmm{F}(s/n)}^n
h},\,\,s\in[0,\,l).$$
\end{proof}

\begin{lemma}\label{t:44} Let
$\mathrm{F}\in\EuScript{F}_{\mmf{X}}$. Assume that for some $t,l>0$ and separable closed linear subspace
$\Phi\subset\mmf{X}$ there exists a strictly increasing natural sequence $\{g_k\}_{k=1}^{\infty}$, such that there exists the limit
\begin{equation}
\lim\limits_{k
\rightarrow\infty}{\{\mathrm{F}({\textstyle\frac{t}{g_k}})\}^{[{g_k}\frac{s}{t}]}
g}
\end{equation}
 for any $g\in \Phi$ and  $s\in [0,\,l)$.
  Then the family of the operators $\Psi_s: \Phi\mapsto\mmf{X}$,
 $s\in [0,\,l)$, defined by the equality
\begin{equation}
\Psi_s g=\lim\limits_{k
\rightarrow\infty}{\{\mathrm{F}({\textstyle\frac{t}{g_k}})\}^{[{g_k}\frac{s}{t}]}
g},\nonumber
\end{equation}
  satisfies the following condition:
 if $ f\in \Phi\cap\mcl{D}(\ov{\mmm{F}_{ef}'(0)})$ and
$\Psi_{s} f\in \mcl{D}(\ov{\mmm{F}_{ef}'(0)})$ for any $s\in
[0,\,l)$, then
\begin{eqnarray} (\Psi_{s}f)'_s
=\ov{\mmm{F}_{ef}'(0)}\Psi_{s}f=w\mbox{-}\lim\limits_{k
\rightarrow\infty}{\{\mathrm{F}({\textstyle\frac{t}{g_k}})\}^{[{g_k}
\frac{s}{t}]}\ov{\mmm{F}'_{ef}(0)}f}\nonumber
\end{eqnarray}
for each $s\in [0,\,l).$
\end{lemma}
\begin{proof}[Proof of Lemma~\ref{t:44}.]
Fix arbitrary $\phi\in\mmm{X}^*$ and $s\in[0,\,l)$.
Denote
$$\Phi'=\ov{\span1{\{\ov{\mmm{F}'_{ef}(0)}\Psi_s f, \Psi_s f,
\ov{\mmm{F}'_{ef}(0)} f,  f, \}}}.$$ By Proposition~\ref{t:1}, we can choose a subsequence
$\{f_k\}_{k=1}^{\infty}$ of the sequence $\{g_k\}_{k=1}^{\infty}$
such that the family of the functions
$\mathrm{T}_{v}:\Phi'\times\{\phi\}\mapsto\mmb{T},\,\,v\geq 0,$
is defined by the equality
$$\mathrm{T}_{v}(g,\phi)=\lim\limits_{k
\rightarrow\infty}{(\{\mathrm{F}(t/{f_k})\}^{[{f_k}\frac{v}{t}]}
g,\phi)}, \,\,g\in\Phi'.$$ From Lemma~\ref{n:2} it follows that
$$(\mathrm{T}_{v}
(\Psi_s f,\phi))^{'}_{v}=\mathrm{T}_{v}(\ov{\mmm{F}'_{ef}(0)}\Psi_s
f,\phi).$$ Let us show that
$$(\mathrm{T}_{v} (\Psi_s f,\phi))= (\Psi_{s+v} f,\phi),\,\, v\in[0,\,l-s).$$
Note that
\begin{eqnarray}
&&(\mathrm{T}_{v} (\Psi_s f,\phi))=
\lim\limits_{k\to\infty}{(\{\mathrm{F}({\textstyle\frac{t}{f_k}})\}^{[{f_k}
\frac{v}{t}]}\Psi_s f,\phi)}\nonumber \\&=&
\lim\limits_{k\to\infty}{(\{\mathrm{F}({\textstyle\frac{t}{f_k}})\}^{[{f_k}
\frac{v}{t}]}(\Psi_s
f-\{\mathrm{F}({\textstyle\frac{t}{f_k}})\}^{[{f_k}
\frac{s}{t}]}f),\phi)}\nonumber\\&+&
\lim\limits_{k\to\infty}{(\{\mathrm{F}({\textstyle\frac{t}{f_k}})\}^{[{f_k}
\frac{v}{t}]}\{\mathrm{F}({\textstyle\frac{t}{f_k}})\}^{[{f_k}
\frac{s}{t}]}f,\phi)}\nonumber\\&=&
\lim\limits_{l\to\infty}{(\{\mathrm{F}({\textstyle\frac{t}{f_k}})\}^{[{f_k}
\frac{v}{t}]}\{\mathrm{F}({\textstyle\frac{t}{f_k}})\}^{[{f_k}
\frac{s}{t}]}f,\phi)}\nonumber\\&=&
\lim\limits_{k\to\infty}{(\{\mathrm{F}({\textstyle\frac{t}{f_k}})\}^{[{f_k}
\frac{v+s}{t}]}f,\phi)}= ((\Psi_{s+l} f,\phi).\nonumber
\end{eqnarray}
Thus, we have
$$(\mathrm{T}_{v} (\Psi_s f,\phi))= (\Psi_{s+v} f,\phi)=(\mathrm{T}_{s+v} ( f,\phi)),\,\, v\in[0,\,l-s).$$
From this, by part (c) of Proposition~\ref{t:1}, we deduce that
\begin{eqnarray}(\Psi_{s+v} f,\phi)^{'}_{v}&=&(\mathrm{T}_{v}
(\Psi_s f,\phi))^{'}_{v}=\mathrm{T}_{v}(\ov{\mmm{F}'_{ef}(0)}\Psi_s
f,\phi) \nonumber\\&=&\mathrm{T}_{v+s}(\ov{\mmm{F}'_{ef}(0)}
f,\phi),\,\,s\in[0,\,l).\nonumber\end{eqnarray} So, by part (c) of Lemma~\ref{t:10}, we get
\begin{eqnarray}&&(\Psi_{s}
f,\phi)^{'}_{s}=(\ov{\mmm{F}'_{ef}(0)}\Psi_s f,\phi)
\nonumber\\&=&\lim\limits_{k
\rightarrow\infty}{(\{\mathrm{F}({\textstyle\frac{t}{g_k}})\}^{[{g_k}
\frac{s}{t}]}\ov{\mathrm{F}'_{ef}(0)}f,\phi)},
\,\,s\in[0,\,l).\nonumber\end{eqnarray} By the arbitrariness of $\phi$,
we have
$$\ov{\mmm{F}'_{ef}(0)}\Psi_s f=w\mbox{-}\lim\limits_{k
\rightarrow\infty}{\{\mathrm{F}({\textstyle\frac{t}{g_k}})\}^{[{g_k}
\frac{s}{t}]}\ov{\mmm{F}'_{ef}(0)}f},\,\,s\in[0,\,l).$$
Therefore, by part (b) of Lemma~\ref{t:10}, we deduce that the function
$\ov{\mmm{F}'_{ef}(0)}\Psi_s f$, $s\in[0,\,l)$, is continuous with respect to $s$. Thus, we get
\begin{eqnarray}&&(\Psi_{m}f-\Psi_{p}f,\phi)=(\Psi_{m}f,
\phi)- (\Psi_{p}f,
\phi)\nonumber\\&=&\int_{p}^{m}{(\ov{\mmm{F}'_{ef}(0)}\Psi_s f,
\phi)\,d s}=(\int_{p}^{m}{\ov{\mmm{F}'_{ef}(0)}\Psi_s f \,d s},
\phi),\nonumber\\ && \,\, m,\,\,p\in[0,\,l),\nonumber
\end{eqnarray}
for arbitrary $\phi\in\mmf{X}^*$; moreover,
$$\int_{p}^{m}{\ov{\mmm{F}'_{ef}(0)}\Psi_s f \,d
s}\in\tilde{\mmf{X}},$$ where $\tilde{\mmf{X}}$ is the completion of the space $\mmf{X}$. By the arbitrariness of $\phi$, we infer that
$$\Psi_{m}f-\Psi_{p}f=\int_{p}^{m}{\ov{\mmm{F}'_{ef}(0)}\Psi_s f \,d s},\,\, m,\,\,p\in[0,\,l).$$
Consequently, by the continuity of the function $\ov{\mmm{F}'_{ef}(0)}\Psi_s
f$ with respect to $s$, we get
$$(\Psi_{s}f)^{'}_{s}=\ov{\mmm{F}'_{ef}(0)}\Psi_{s}f$$ for
each $s\in[0,\,l).$   \end{proof}

\begin{theorem}\label{vvt:54}
Let  $\mmm{Z}$ be a densely defined linear operator in $\mmf{X}$ and $t>0$.
Assume that there exists the function
$\mathrm{F}\in\EuScript{F}_{\mmf{X}}$ such that
$\mcl{D}(\mmm{F}_{ef}'(0))\supset\mcl{D}(\mmm{Z})$ and
$\mmm{F}_{ef}'(0)g=\mmm{Z}g,\,\,g\in\mcl{D}(\mmm{Z})$.

\noindent Then the local solution $f:[0,\,l)\mapsto\mmm{X},\,\,l>0,$
of the system
\begin{eqnarray}
f'(v)&=&\overline{\mmm{Z}}f(v),\nonumber\\
f(0)&=&h\in\mcl{D}(\overline{\mmm{Z}}),\nonumber
\end{eqnarray}
exists on the semi-interval $[0,\,l)$ iff the following conditions are satisfied:

\noindent a)  It is possible to choose a  convergent sub-subsequence
for any subsequence of the sequence $\{\mathrm{F}^{[n
s/t]}(\frac{t}{n})h\}_{n\in \Nat}$,
 for any $s\in[0,\,l)$;

\noindent b) For any $s\in[0,\,l)$ there exists a sequence
$\{f_n^{s}\}_{n=1}^{\infty}$, $f_n^{s}\in\mcl{D}(\ov{\mmm{Z}}),$
such that there exists
$w\mbox{-}\lim\limits_{n\to\infty}{(\mathrm{F}^{[n
s/t]}(\frac{t}{n})h-f_n^{s})}=0$ and it is possible to choose a
weakly convergent sub-subsequence for any subsequence of the
sequence $\{\ov{\mmm{Z}}f_n^{s}\}_{n\in \Nat}$.

\noindent Moreover, if conditions (a)-(b) are satisfied, then
$$f(v)=\lim\limits_{n\to\infty}{\{{\mmm{F}(v/n)}\}^n h},\,v\in[0,\,l).$$
\end{theorem}

\begin{proof}[Proof of Theorem~\ref{vvt:54}.]
($\Rightarrow$) It follows from Theorem~\ref{chern1} and
$f(v)\in\mcl{D}(\ov{\mmm{Z}})$ for any $v\in[0,\,l)$.

\noindent($\Leftarrow$) Let $\{s_k\}_{k=1}^{\infty}$ be a sequence of all nonnegative rational numbers contained in $[0,\,l)$.
  Fix a strictly increasing natural sequence $\{n_i\}_{i=1}^{\infty}$. Choose a subsequence
$\{n_i^1\}_{i=1}^{\infty}$ of the sequence
$\{n_i\}_{i=1}^{\infty}$ such that there exists the limit
 $\lim\limits_{i\to\infty}{\mmm{F}^{[{n_i^1}
s_1/t]}(t/{n_i^1})h}$. Further, choose a subsequence
$\{n_i^2\}$ of the sequence $\{n_i^1\}$ such that there exists the
limit $\lim\limits_{i\to\infty}{\mmm{F}^{[{n_i^2}
  s_2/t]}(t/{n_i^2})h}$. In the same way, for any natural $k\geq 3$, we choose
the subsequence $\{n_i^k\}_{i=1}^{\infty}$
of the subsequence $\{n_i^{k-1}\}_{i=1}^{\infty}$ such that there exists the
limit $\lim\limits_{i\to\infty}{\mmm{F}^{[{n_i^k}
s_k/t]}(t/{n_i^k})h}$. Note that for diagonal sequence $\{n_i^i\}_{i=1}^{\infty}$ there exists the
limit $\lim\limits_{i\to\infty}{\mmm{F}^{[n_i^i
s_k/t]}(t/{n_i^i})h}$ for each $k\in\mmb{N}$. Thus, by Lemma~\ref{t:7} we infer that
$\{\mmm{F}^{[n_i^i
  s/t]}(t/{n_i^i})h\}_{i=1}^{\infty}$ is a Cauchy sequence in $\mmf{X}$ for any
   $s\in[0,\,l)$. Therefore, by condition~(a), there exists
$$\lim\limits_{i\to\infty}{\mmm{F}^{[n_i^i
  s/t]}(t/{n_i^i})h}$$ for any  $s\in[0,\,l).$
  Denote $d_i=n_i^i,\,\,i\in\mmb{N},$ and
  $$\mmm{T}_s g=\lim\limits_{i\to\infty}{\mmm{F}^{[d_i
  s/t]}(t/{d_i})g}, \,\,g\in\span1{(\{h\})},\,\,s\in[0,\,l).$$ Fix arbitrary $m\in[0,\,l)$. Choose a subsequence
  $\{d'_k\}_{k\in\Nat}$ of the sequence $\{d_k\}_{k\in\Nat}$ such
that there exists the limit $w\mbox{-}\lim\limits_{i\to\infty}{\ov{\mmm{Z}}f_{d'_i}^{m}}$. Thus, there exist
   $w\mbox{-}\lim\limits_{i\to\infty}{f_{d'_i}^{m}}$,
  $w\mbox{-}\lim\limits_{i\to\infty}{\ov{\mmm{Z}}f_{d'_i}^{m}}$ and, therefore, from the coincidence of the weak closure of
  $\mmm{Z}$ with the closure in the strong
topology it follows that
   $\mmm{T}_{m}h\in\mcl{D}(\overline{\mmm{Z}})$ for any
  $m\in[0,\,l)$. Then, from Lemma~\ref{t:44}
  we infer that
  $$(\mathrm{T}_{s}h)^{'}_{s}=\overline{\mmm{Z}}\mathrm{T}_{s}h$$ for each $s\in[0,\,l)$.
  Consequently,  $\mathrm{T}_{s}h$
is the solution of the equation $$f'(s)=\overline{\mmm{Z}}f(s)$$ with the initial condition $f(0)=h$ on the semi-interval $[0,\,l)$. Thus, by Theorem~\ref{chern1} we get
$$f(s)=\lim\limits_{n\to\infty}{{\mmm{F}(s/n)}^n
h},\,\,s\in[0,\,l).$$
\end{proof}
From Theorem~\ref{vvt:54} we get the following consequence.
\begin{corollary} \label{t:55} Let $t>0$. Assume that $\mmm{C}$
and $\mmm{D}$ are generators of the $lC_0$-semigroups
$\exp{(s\mmm{C})}$ and $\exp{(s\mmm{D})}$ in sequentially
complete locally convex space $\mmf{X}$. Assume also, that the following conditions are satisfied:

\noindent i) $\mcl{D}(\mmm{C})\cap\mcl{D}(\mmm{D})$ is dense in
$\mmf{X}$;

\noindent ii) $\{[0,\,\infty)\ni
s\mapsto\exp{{\textstyle(s\mmm{C})}}\exp{{\textstyle(s\mmm{D})}}\}
\in\EuScript{E}_{\mmf{X}}$.

\noindent Then the local solution $f:[0,\,l)\mapsto\mmm{X},\,\,l>0,$
of the system
\begin{eqnarray}
f'(v)&=&(\overline{\mmm{C}+\mmm{D}})f(v),\nonumber\\
f(0)&=&h\in\mcl{D}(\overline{\mmm{C}+\mmm{D}}),\nonumber
\end{eqnarray}
exists on the semi-interval $[0,\,l)$ iff the following conditions are satisfied:

\noindent a) It is possible to choose a convergent sub-subsequence
for any subsequence of the sequence
$\{\{\exp{{\textstyle(\frac{t}{n}\mmm{C})}}\exp{{\textstyle(\frac{t}{n}\mmm{D})}}\}^{[n
s]}h\}_{n\in \Nat}$, for any $s\in[0,\,l/t)$;

\noindent b) For any $s\in[0,\,l/t)$ there exists a sequence
$\{f_n^{s}\}_{n=1}^{\infty}$,
$f_n^{s}\in\mcl{D}(\ov{\mmm{C}+\mmm{D}}),$ such that there exists
$w\mbox{-}\lim\limits_{n\to\infty}{\{\exp{{\textstyle(\frac{t}{n}\mmm{C})}}\exp{{\textstyle(\frac{t}{n}\mmm{D})}}\}^{[n
s]}f-f_n^{s})}=0$ and it is possible to choose a weakly convergent
sub-subsequence for any subsequence of the sequence
$\{({\ov{\mmm{C}+\mmm{D}}})f_n^{s}\}_{n=1}^{\infty}$.

\noindent Moreover, if conditions (a)-(b) are satisfied, then
$$f(s) =\lim\limits_{n\to\infty}{\{\exp{{\textstyle (\frac{s}{n}\mmm{C})}}\exp{{\textstyle (\frac{s}{n}\mmm{D})}}\}^n
f},\, s\in[0,\,l),$$ for each $f\in\mmm{X}.$
\end{corollary}

\begin{proof}[Proof of Corollary~\ref{t:55}.] Put $\mathrm{F}(s)=\exp(s \mmm{C})\exp(s \mmm{D})$, $s\in[0,\,\infty)$.
From
\begin{eqnarray}
 & &\lim\limits_{s \rightarrow 0}{s^{-1}(\mmm{F}(s)f-\mmm{F}(0)f)}\nonumber\\&=&\lim\limits_{s \rightarrow
0}{s^{-1}(\exp(s \mmm{C})\exp(s \mmm{D})f-\exp(s
\mmm{C})f)}\nonumber\\&+& \lim\limits_{s \rightarrow
0}{s^{-1}(\exp(s
\mmm{C})f-f)}=\mmm{C}f+\mmm{D}f,\,\,f\in\mcl{D}(\mmm{C})\cap\mcl{D}(\mmm{D}),\nonumber
\end{eqnarray}
it follows that $\mcl{D}(\mmm{F}_{ef}'(0))$ is dense in $\mmm{X}$.
Thus, from condition~(ii), we deduce that
$\mmm{F}\in\EuScript{F}_{\mmf{X}}$. From Theorem~\ref{cor:1} we infer
that the operator $\mmm{C}+\mmm{D}$ possesses a closure. Applying Theorem~\ref{vvt:54} for the operator $\mmm{Z}=\mmm{C}+\mmm{D}$ we get the Corollary.
\end{proof}

\section{Criteria for closure of operator to be the generator of the $lC_0$-semigroup}\label{3.2}

\begin{theorem}\label{cog:1} Let $\mmf{X}$ be a  complete locally convex space and $\mathrm{F}\in\EuScript{F}_{\mmf{X}}$. Assume that a linear operator $\mmm{Z}$ has the domain
$\mcl{D}(\mmm{Z})\subset\mcl{D}(\mmm{F}_{ef}'(0))$ and
$$\mmm{Z}g=\mmm{F}_{ef}'(0)g,\,\,g\in\mcl{D}(\mmm{Z}).$$  Assume also that $\mmf{A}\subset\mcl{D}(\mmm{Z})$ is a dense linear subset of  $\mmf{X}$ and there exists a fixed
$l>0$ such
that there exists a local solution
$f:[0,\,l)\mapsto\mmm{X}$ of the system
\begin{eqnarray}
f'(s)&=&\overline{\mmm{Z}}f(s),\,\,s\in[0,\,l),\label{tt:21}\\
f(0)&=&f_0\label{tt:21a}
\end{eqnarray}
for each $f_0\in\mmf{A}$. Then the operator $\mmm{Z}$ possesses a closure and
its closure is the generator of the $lC_0$-semigroup
$\mmm{S}$. Furthermore,
the following equality is satisfied:
\begin{eqnarray} \mmm{S}(t) g=\lim\limits_{n\to\infty}{{\mmm{F}(t/n)}^n
g},\, t\geq 0,\label{st:14}
\end{eqnarray} for all $g\in\mmm{X}.$
\end{theorem}
\begin{proof}[Proof of Theorem~\ref{cog:1}.]
By Theorem~\ref{cor:2}, there exists
\begin{eqnarray} w\mbox{-}\lim\limits_{n\to\infty}{{\mmm{F}(t/n)}^{[ns/t]}
f},\,\, s\in[0,\,l),\label {st:12}
\end{eqnarray} for each $f\in\mmf{A}$. From the density of
$\mmf{A}$ in $\mmf{X}$ and Lemma~\ref{t:10} we easily infer that there
exists limit~\eqref{st:12} for each $f\in\mmf{X}$. Put
$\mmm{G}(s)
f=w\mbox{-}\lim\limits_{n\to\infty}{{\mmm{F}(t/n)}^{[ns/t]} f},\,
s\in[0,\,l),\,f\in\mmf{X}$. From the uniqueness of the local solution of system \eqref{tt:21}--\eqref{tt:21a} we infer that
\begin{eqnarray}
\mmm{G}(s_1)\mmm{G}(s_2)f=\mmm{G}(s_1+s_2)f,\,s_1,\,s_2,\,s_1+s_2\in[0,\,l),\label{st:13}
\end{eqnarray}
for each $f\in\mmf{A}$. From this and part (a) of Lemma~\ref{t:10}
it follows that equality (\ref{st:13}) holds for any
$f\in\mmf{X}$.
 Define the function
$\mmm{S}:[0,\,\infty)\mapsto\mcl{L}(\mmf{X})$ by the equality
$$\mmm{S}(s)=\{\mmm{G}(l/2)\}^{[2s/l]}\mmm{G}(s-[2s/l]l/2),\,\,s\geq 0.$$ Then,
from parts (a), (b) of Lemma~\ref{t:10} and the definition of $\mmm{S}$
we easily deduce that the function  $\mmm{S}$ is the $lC_0$-semigroup. Furthermore, from part (c) of Lemma~\ref{t:10} it follows that
$\mcl{D}(\mmm{Z})\subset\mcl{D}(\mmm{F}_{ef}'(0))\subset\mcl{D}(\mmm{S}'(0))$
and, by the closedness of $\mmm{S}'(0)$,
$\overline{\mmm{Z}}f=\mmm{S}'(0)f$ for each
$f\in\mcl{D}(\overline{\mmm{Z}})$. Let us show that
$\mcl{D}(\overline{\mmm{Z}})=\mcl{D}(\mmm{S}'(0))$. Fix
$g\in\mmf{A}$. Then
$$\sum_{k=1}^{n}{{\textstyle\frac{s}{n}\mmm{S}(k\frac{s}{n})g}}\in\mcl{D}(\overline{\mmm{Z}}),\,\,s\in[0,\,l/2],$$ and, by properties of $lC_0$-semigroups, we infer that there exist the limits
\begin{eqnarray}\lim\limits_{n\to\infty}{\sum_{k=1}^{n}{{\textstyle\frac{s}{n}\mmm{S}(k\frac{s}{n})g}}}&=&\int_{0}^{s}{\mmm{S}(t)g\,d
t},\nonumber\\
\lim\limits_{n\to\infty}{\ov{\mmm{Z}}\sum_{k=1}^{n}{{\textstyle\frac{s}{n}\mmm{S}(k\frac{s}{n})g}}}&=&\int_{0}^{s}{\mmm{S}(t)\ov{\mmm{Z}}g\,d
t}=\mmm{S}(s)g-g,\,\,s\in[0,\,l/2].\nonumber
\end{eqnarray}  Thus, by the closedness of $\overline{\mmm{Z}}$, we deduce that
 $\int_{0}^{s}{\mmm{S}(t)g\,d
t}\in\mcl{D}(\overline{\mmm{Z}})$ and
$$\ov{\mmm{Z}}\int_{0}^{s}{\mmm{S}(t)g\,d t}=\mmm{S}(s)g-g,\,\,s\in[0,\,l/2].$$
Let $f\in\mmf{X}$. Fix $\alpha\in\Omega$. Choose a sequence $\{f_n\}_{n=1}^{\infty},\,\,f_n\in\mmf{A}$, such that
$\lim\limits_{n\to\infty}{\|f_n-f\|_{\alpha}^{\mmm{F},l}}=0$. Then we have
$$\|\lim\limits_{n\to\infty}{\int_{0}^{s}{\mmm{S}(t)f_n\,d
t}}-\int_{0}^{s}{\mmm{S}(t)f\,d t}\|_{\alpha}=0,$$
$$\|\lim\limits_{n\to\infty}{\overline{\mmm{Z}}\int_{0}^{s}{\mmm{S}(t)f_n\,d
t}}-(\mmm{S}(s)f-f)\|_{\alpha}$$
$$=\|\lim\limits_{n\to\infty}{\mmm{S}(s)f_n-f_n}-(\mmm{S}(s)f-f)\|_{\alpha}=0,\,s\in[0,\,l/2].$$
 Therefore, by the arbitrariness of $\alpha\in\Omega$
 and the closedness of $\overline{\mmm{Z}}$, we infer that
$$\int_{0}^{s}{\mmm{S}(t)f\,d t}\in\mcl{D}(\overline{\mmm{Z}})$$ and
$$\overline{\mmm{Z}}\int_{0}^{s}{\mmm{S}(t)f\,d t}=\mmm{S}(s)f-f$$
for $s\in[0,\,l/2]$.  Note that for $f\in\mcl{D}(\mmm{S}'(0))$
there exist
$$\lim\limits_{s\to 0}{s^{-1}\int_{0}^{s}{\mmm{S}(t)f\,d t}}=f,$$
$$\lim\limits_{s\to
0}{s^{-1}\overline{\mmm{Z}}\int_{0}^{s}{\mmm{S}(t)f\,d
t}}=\lim\limits_{s\to 0}{s^{-1}(\mmm{S}(s)f-f)}=\mmm{S}'(0)f.$$
So, by the closedness of $\overline{\mmm{Z}}$ it follows that
$f\in\mcl{D}(\overline{\mmm{Z}})$. Therefore,
$\overline{\mmm{Z}}$ is the generator of the $lC_0$-semigroup
$\mmm{S}$ and, by Theorem~\ref{chern}, we get equality
\eqref{st:14}.
\end{proof}
\begin{remark}\label{rem:1.1}
If $\mmf{A}$ is sequentially dense in $\mmf{X}$, then the condition that $\mmf{X}$ is a complete locally convex space in Theorem~\ref{cog:1}
can be changed to the condition that $\mmf{X}$ is a sequentially complete locally convex space.
\end{remark}

\begin{lemma}\label{st:31} Let $\mmf{X}$ be a sequentially complete locally convex space and
a function $\mmm{T}:[0,\,\infty)\to\mcl{L}(\mmf{X})$
is a $lC_0$-semigroup. Then $\mcl{D}({\mmm{T}^{*}}'(0))$
is *-dense in $\mmf{X}^{*}$, where ${\mmm{T}^{*}}'(0)$ is the (strong) derivative at the point $0$ of the function $\mmm{T}^{*}$.
\end{lemma}
\begin{proof}[Proof of Lemma~\ref{st:31}.]
Firstly, consider the case in which $\mmf{X}$ is a complete locally convex space.
Denote $\EuScript{D}((0,\,\infty))$ the set of all smooth functions  $\phi:(0,\,\infty) \to \mmb{R}$ with compact support
$\supp{\phi}\subset(0,\,\infty)$. Let $\mmf{Y}$ consist of all
$f\in\mmf{X}^*$ for which there
exist $g\in\mmf{X}^*$,
$\phi\in\EuScript{D}((0,\,\infty))$ such that
$$f=\int_0^{\infty}{\phi(s)\mmm{T}^*(s)g\,ds},$$
where the integral is in the sense of Pettis with respect to the topology $\sigma(\mmf{X}^*,\mmf{X})$ (\cite{Pettis}).
Note that $\mmf{Y}$  is *-dense in $\mmf{X}^{*}$. Indeed, we can choose a sequence of non-negative
functions $\phi_n,\,\,n\in\mmb{N}$,
$\phi_n\in\EuScript{D}((0,\,\infty))$, such that
$\supp{\phi_n}\in(0,\,1/n]$ and
$$\int_0^{\infty}{\phi_n(s)\,ds}=1.$$
Then $\phi_n$ tends to $\delta-$function as
$n\to\infty.$ Therefore, for any $f\in\mmf{X}^*$ and
$x\in\mmf{X}$ we see that
$$\lim\limits_{n\to\infty}{(x,\int_0^{\infty}{\phi_n(s)\mmm{T}^*(s)f\,ds})}=(x,f).$$
Now let us show that $\mmf{Y}\subset\mcl{D}({\mmm{T}^{*}}'(0)).$
Fix
$g\in\mmf{X}^*$, $\phi\in\EuScript{D}((0,\,\infty))$ and a bounded set $\mmf{A}\subset\mmf{X}$. Let
$\supp{\phi}\in(0,\,a),\,\,a>0,$ and
$$f=\int_0^{\infty}{\phi(s)\mmm{T}^*(s)g\,ds}.$$ Choose $t_0>0$
such that $[0,\,t_0]\cap\supp{\phi}=\emptyset.$ Then for
$t\in(0,\,t_0)$ we have
\begin{eqnarray}
&&\|\mmm{T}^*(t)f-f+t\int_{0}^{\infty}{\phi'(s)\mmm{T}^*(s)g\,ds}\|_{\mmf{A}^{\circ}}
\nonumber\\&=&\|\int_0^{\infty}{\phi(s)\mmm{T}^*(s+t)g\,ds}-\int_0^{\infty}{\phi(s)\mmm{T}^*(s)g\,ds}\nonumber\\&+&
\int_{0}^{\infty}{t\phi'(s)\mmm{T}^*(s)g\,ds}\|_{\mmf{A}^{\circ}}
=\|\int_t^{\infty}{\phi(s-t)\mmm{T}^*(s)g\,ds}\nonumber\\&-&\int_0^{\infty}{\phi(s)\mmm{T}^*(s)g\,ds}+
\int_{0}^{\infty}{t\phi'(s)\mmm{T}^*(s)g\,ds}\|_{\mmf{A}^{\circ}}\nonumber\\
&=&\|\int_t^{a+t_0}{(t\phi'(s)+\phi(s-t)-\phi(s))\mmm{T}^*(s)g\,ds}\|_{\mmf{A}^{\circ}}\nonumber\\
&=&\|\int_t^{a+t_0}{\int_{s-t}^s{(\phi'(s)-\phi'(l))\,dl}\mmm{T}^*(s)g\,ds}\|_{\mmf{A}^{\circ}}=(A)\nonumber
\end{eqnarray}
Denote $b(s)=\sup\limits_{|l-k|\leq
s}{|\phi'(l)-\phi'(k)|},\,\,l,k\in(0,\,a).$ Then
\begin{eqnarray}
(A)&\leq&(t_0+a)tb(t)\sup\limits_{s\in[0,\,a+t_0]}{\|\mmm{T}^*(s)g\|_{\mmf{A}^{\circ}}}\nonumber\\
&=&(t_0+a)tb(t)\sup\limits_{s\in[0,\,a+t_0]}{\sup\limits_{x\in\mmf{A}}{|(x,\mmm{T}^*(s)g)|}}\nonumber\\&=&
(t_0+a)tb(t)\sup\limits_{s\in[0,\,a+t_0]}{\sup\limits_{x\in\mmf{A}}{|(\mmm{T}(s)x,g)|}}\nonumber\\&\leq&
(t_0+a)tb(t)\sup\limits_{s\in[0,\,a+t_0]}{\sup\limits_{x\in\mmf{A}}{\|\mmm{T}(s)x\|_{g}}}\nonumber\\&\leq&
(t_0+a)tb(t){\sup\limits_{x\in\mmf{A}}{\|x\|^{\mmm{T},a+t_0}_{g}}}=(B_t).\nonumber
\end{eqnarray}
Since the semi-norm
$\|\cdot\|^{\mmm{T},a+t_0}_{g}$ is continuous  in $\mmf{X}$ we infer that ${\sup\limits_{x\in\mmf{A}}{\|x\|^{\mmm{T},a+t_0}_{g}}}<\infty$
and, therefore, $\lim\limits_{t\to 0}{t^{-1}B_{t}}=0.$ By the arbitrariness of $\mmf{A}$, we infer that $f\in\mcl{D}({\mmm{T}^{*}}'(0)).$  So, we've shown  the Lemma in the case under consideration.

Consider the general case in which $\mmf{X}$ is a sequentially complete locally convex space. Let $\tilde{\mmf{X}}$  be the completion of the space $\mmf{X}$ and the function $\tilde{\mmm{T}}:[0,\,\infty)\to\mcl{L}(\tilde{\mmf{X}})$ be defined by the condition that $\tilde{\mmm{T}}(s)$ is the continuous extension of $\mmm{T}(s)$ to $\tilde{\mmf{X}}$ for each $s\in[0,\,\infty)$. Then $\tilde{\mmm{T}}$ is the $lC_0$-semigroup, $\tilde{\mmf{X}}^*=\mmf{X}^*$ and $\tilde{\mmm{T}}^*=\mmm{T}^*$.
Now, the Lemma follows from previous case.
\end{proof}

\begin{theorem}\label{vt:4} Assume that  $\mmm{Z}$ is a densely defined linear operator in
a complete locally convex space  $\mmf{X}$ and
$t>0$. Assume also that there
exists a function
$\mathrm{F}:[0,\,\infty)\to\mcl{L}({\mmf{X}})$ such that:

\noindent i) $\mmm{F}\in\EuScript{F}_{\mmf{X}}$;

\noindent ii) $\mcl{D}(\mmm{F}_{ef}'(0))\supset\mcl{D}(\mmm{Z})$  and
$\mmm{F}_{ef}'(0)g=\mmm{Z}g,\,\,g\in\mcl{D}(\mmm{Z})$;

\noindent iii) $\mcl{D}({\mmm{F}^*}'(0))$ is  *-dense in
$\mathbf{X}^*$.

\noindent Then the operator $\mmm{Z}$ possesses a closure and its closure
is the generator of the $lC_0$-semigroup iff
there exists a dense
linear subspace
$\mmf{A}\subseteq\mcl{D}(\mmm{Z})$ such that for any
$f\in\mmf{A}$ and $s\geq 0$ there exists a sequence
$\{f_n^{s}\}_{n=1}^{\infty}$, $f_n^{s}\in\mcl{D}({\ov{\mmm{Z}}}),$
satisfying the following conditions:

\noindent a) $w\mbox{-}\lim\limits_{n\to\infty}{(\mathrm{F}^{[n
s/t]}(\frac{t}{n})f-f_n^{s})}=0$;

\noindent b) It is possible to choose a weakly convergent
sub-subsequence for any subsequence of the sequence
$\{{\ov{\mmm{Z}}}f_n^{s}\}_{n=1}^{\infty}$;

\noindent c) It is possible to choose a weakly convergent
sub-subsequence for any subsequence of the sequence
$\{f_n^{s}\}_{n=1}^{\infty}$.
\end{theorem}
\begin{proof}[Proof of Theorem~\ref{vt:4}.]
($\Rightarrow$) It follows from Theorem~\ref{vvmain:10a}.

\noindent ($\Leftarrow$)  By Theorem~\ref{vvmain:10a} there exists
the solution of the equation $$f'(s)=\overline{\mmm{Z}}f(s)$$ on
$[0,\,\infty)$ for any initial condition  $f(0)=g$, where
$g\in\mmf{A}$. Thus, by Theorem~\ref{cog:1},
$\overline{\mmm{Z}}$ is the generator of the $lC_0$-semigroup
$\mmm{S}$ and
$$\mmm{S}(t) f=\lim\limits_{n\to\infty}{{\mmm{F}(t/n)}^n
f},\,\,f\in\mmm{X}.$$
\end{proof}

\begin{theorem}\label{vmain:1} Assume that $\mmm{Z}$ is a densely defined linear operator in a complete locally
convex space $\mmf{X}$ and $t>0$. Then the operator $\mmm{Z}$
 possesses a closure and its closure is the generator of the $lC_0$-semigroup
iff there exists a function
$\mathrm{F}\in\EuScript{F}_{\mmf{X}}$ such that:

\noindent i) $\mcl{D}(\mmm{F}_{ef}'(0))\supset\mcl{D}(\mmm{Z})$ and
$\mmm{F}_{ef}'(0)g=\mmm{Z}g,\,\,g\in\mcl{D}(\mmm{Z})$;

\noindent ii) $\mcl{D}({\mmm{F}^*}'(0))$ is *-dense in $\mathbf{X}^*$;

\noindent iii)  There exists a dense linear subspace
$\mmf{A}\subseteq\mcl{D}(\mmm{Z})$ such that for any
$f\in\mmf{A}$, $s\geq 0$ there
exists a sequence
$\{f_n^{s}\}_{n=1}^{\infty}$, $f_n^{s}\in\mcl{D}({\ov{\mmm{Z}}}),$
satisfying the following conditions:

\noindent a) $w\mbox{-}\lim\limits_{n\to\infty}{(\mathrm{F}^{[n
s/t]}(\frac{t}{n})f-f_n^{s})}=0$;

\noindent b) It is possible to choose a weakly convergent
sub-subsequence for any subsequence of the sequence
$\{{\ov{\mmm{Z}}}f_n^{s}\}_{n\in\Nat}$;

\noindent c) It is possible to choose a weakly convergent
sub-subsequence for any subsequence of the sequence
$\{f_n^{s}\}_{n\in\Nat}$.

\noindent Furthermore, if conditions  (i)-(iii) are satisfied, then
$\exp{(t\overline{\mmm{Z}})}
f=\lim\limits_{n\to\infty}{{\mmm{F}(t/n)}^n f}$ for each $f\in\mmm{X}.$
\end{theorem}

\begin{proof}[Proof of Theorem~\ref{vmain:1}.]
($\Rightarrow$) It follows from Theorems~\ref{vt:4} and~\ref{chern}.

\noindent ($\Leftarrow$) It follows from Lemma~\ref{st:31} and Theorem~\ref{vt:4}.
\end{proof}

\begin{corollary}\label{vmain:10} Assume that $\mmm{Z}$ is a densely defined linear operator in a complete locally
convex space $\mmf{X}$ and $t>0$. Then the operator $\mmm{Z}$
 possesses a closure and its closure is the generator of the $lC_0$-semigroup
iff there exists a function
$\mathrm{F}\in\EuScript{F}_{\mmf{X}}$ such that:

\noindent i) $\mcl{D}(\mmm{F}_{ef}'(0))\supset\mcl{D}(\mmm{Z})$ and
$\mmm{F}_{ef}'(0)g=\mmm{Z}g,\,\,g\in\mcl{D}(\mmm{Z})$;

\noindent ii) $\mcl{D}({\mmm{F}^*}'(0))$ is *-dense in $\mathbf{X}^*$;

\noindent iii) It is possible to choose a weakly convergent
sub-subsequence for any subsequence of the sequence
$\{\mathrm{F}^{[n s]}(\frac{t}{n})f\}_{n\in \Nat}$, for any $s>0$
and $f\in\mmf{X}$;

\noindent iv) There exists a dense linear subspace
$\mmf{A}\subseteq\mcl{D}(\mmm{Z})$ such that for any
$f\in\mmf{A}$ and $s\geq 0$  there
exists a sequence $\{f_n^{s}\}_{n=1}^{\infty}$,
$f_n^{s}\in\mcl{D}(\ov{\mmm{Z}}),$ such that
$$w\mbox{-}\lim\limits_{n\to\infty}{(\mathrm{F}^{[n s]}({\textstyle
\frac{t}{n}})f-f_n^{s})}=0$$ and
$$w\mbox{-}\lim\limits_{n\to\infty}{(\mathrm{F}^{[n
s]}({\textstyle\frac{t}{n}})\mmm{Z}f-\ov{\mmm{Z}}f_n^{s})}=0.$$

Furthermore, if conditions  (i)-(iv) are satisfied, then
$\exp{(t\overline{\mmm{Z}})}
f=\lim\limits_{n\to\infty}{{\mmm{F}(t/n)}^n f}$ for each $f\in\mmm{X}.$
\end{corollary}

\begin{proof}[Proof of Corollary~\ref{vmain:10}.]
 It follows from Theorem~\ref{vt:4}.
\end{proof}

\begin{corollary}\label{vt:3} Let $\mmm{C}$
and $\mmm{D}$ be generators of the $lC_0$-semigroups
$\exp{(s\mmm{C})}$ and $\exp{(s\mmm{D})}$ in a
complete locally convex space $\mmf{X}$ and $t>0$ be fixed. Assume that there exists a set
$\mmf{B}\subset\mcl{D}((\exp(s\mmm{C})^*)'_{s=0})\cap\mcl{D}((\exp(s\mmm{D})^*)'_{s=0})$
and the following conditions are satisfied:

\noindent i) $\mcl{D}(\mmm{C})\cap\mcl{D}(\mmm{D})$ is dense in
$\mmf{X}$;

\noindent ii) $\mmf{B}$  is *-dense in $\mmf{X^*}$;

\noindent iii) $\{[0,\,\infty)\ni
s\mapsto\exp{{\textstyle(s\mmm{C})}}\exp{{\textstyle(s\mmm{D})}}\}
\in\EuScript{E}_{\mmf{X}}$;

\noindent iv) The function $g(s,x)=\exp{(s\mmm{D}^*)}x,\,\, s\geq 0,\,\,
x\in \mmf{X},$ is continuous at the point $s=0$ for each $x\in
(\exp(s\mmm{C})^*)'_{s=0}(\mmf{B})$.

\noindent Then the sum of $\mmm{C}$ and $\mmm{D}$ possesses a closure and its closure is the generator of the $lC_0$-semigroup
iff there exists a   dense linear subspace
$\mmf{A}\subseteq\mcl{D}(\mmm{C})\cap\mcl{D}(\mmm{D})$ such that for any  $f\in\mmf{A}$ and $s\geq 0$ there
exists a sequence
$\{f_n^{s}\}_{n=1}^{\infty}$,
$f_n^{s}\in\mcl{D}(\ov{\mmm{C}+\mmm{D}}),$ satisfying the following conditions:

\noindent a) $w\mbox{-}\lim\limits_{n\to\infty}{((\exp{{(\textstyle
\frac{t}{n}\mmm{C})}}\exp{{(\textstyle\frac{t}{n}\mmm{D})}})^{[n
s]}f-f_n^{s})}=0$;

\noindent b) It is possible to choose a weakly convergent
sub-subsequence for any subsequence of the sequence
$\{f_n^{s}\}_{n=1}^{\infty}$;

\noindent c) It is possible to choose a weakly convergent
sub-subsequence for any subsequence of the sequence
$\{(\ov{\mmm{C}+\mmm{D}})f_n^{s}\}_{n=1}^{\infty}$.

\noindent Furthermore, if conditions  (a)-(c) are satisfied, then
$$\exp{(s(\overline{\mmm{C}+\mmm{D}}))} f=\lim\limits_{n\to\infty}{\{\exp{{\textstyle (\frac{s}{n}\mmm{C})}}\exp{{\textstyle (\frac{s}{n}\mmm{D})}}\}^n
f},\, s\geq 0,$$ for each $f\in\mmm{X}.$
\end{corollary}

\begin{proof}[Proof of Corollary~\ref{vt:3}.] Put $\mathrm{F}(s)=\exp(s \mmm{C})\exp(s \mmm{D})$ for any $s\in[0,\,\infty)$.
By
\begin{eqnarray}
 & &\lim\limits_{s \rightarrow 0}{s^{-1}(\mmm{F}(s)f-\mmm{F}(0)f)}\nonumber\\&=&\lim\limits_{s \rightarrow
0}{s^{-1}(\exp(s \mmm{C})\exp(s \mmm{D})f-\exp(s
\mmm{C})f)}\nonumber\\&+& \lim\limits_{s \rightarrow
0}{s^{-1}(\exp(s
\mmm{C})f-f)}=\mmm{C}f+\mmm{D}f,\,\,f\in\mcl{D}(\mmm{C})\cap\mcl{D}(\mmm{D}),\nonumber
\end{eqnarray}
and
\begin{eqnarray} &&\lim\limits_{s \rightarrow
0}{s^{-1}(\mmm{F}^*(s)-\mmm{F}(0))f}\nonumber\\&=&\lim\limits_{s
\rightarrow 0}{s^{-1}\{(\exp(s \mmm{D}))^*(\exp(s
\mmm{C}))^*-(\exp(s \mmm{D}))^*\}f}\nonumber\\&+& \lim\limits_{s
\rightarrow 0}{s^{-1}\{(\exp(s \mmm{D}))^*f-f\}}=\lim\limits_{l
\rightarrow 0}{(\exp(l
\mmm{D}))^*(\exp(s\mmm{C})^*)'|_{s=0}f}\nonumber\\&+&(\exp(s\mmm{D})^*)'|_{s=0}f=
(\exp(s\mmm{C})^*)'|_{s=0}f+(\exp(s\mmm{D})^*)'|_{s=0}f,\nonumber\\&&f\in
\mmf{B},\,\,\,\,\,\,\,\,\,\,\,\,\,\,\,\,\,\,\,\,\,\,\,\,\nonumber
\end{eqnarray}
we deduce that $\mcl{D}(\mmm{F}_{ef}'(0))$ is dense in $\mmm{X}$ and
$\mcl{D}({\mmm{F}^{*}}'(0))$ is *-dense in $\mmm{X}^*$. Thus,
by condition (iii), we infer that
$\mmm{F}\in\EuScript{F}_{\mmf{X}}$. From Theorem~\ref{cor:1} it follows that operator $\mmm{C}+\mmm{D}$ possesses a closure. Therefore, applying
Theorem~\ref{vt:4} for the operator $\mmm{Z}=\mmm{C}+\mmm{D}$, we get the Corollary. \end{proof}

\begin{theorem}\label{vt:54} Assume that  $\mmm{Z}$ is a densely defined linear operator in a complete locally
convex space $\mmf{X}$ and
$t>0$. Assume also that there exists a function
$\mathrm{F}\in\EuScript{F}_{\mmf{X}}$ such that
$\mcl{D}(\mmm{F}_{ef}'(0))\supset\mcl{D}(\mmm{Z})$  and
$\mmm{F}_{ef}'(0)g=\mmm{Z}g,\,\,g\in\mcl{D}(\mmm{Z})$. Then the operator
$\mmm{Z}$ possesses a closure and its closure is the generator of the
$lC_0$-semigroup iff there exists a dense linear subspace $\mmf{A}\subseteq\mcl{D}(\mmm{Z})$ such that for any $f\in\mmf{A}$ and $s\geq 0$ there
exists a sequence $\{f_n^{s}\}_{n=1}^{\infty}$,
$f_n^{s}\in\mcl{D}({\ov{\mmm{Z}}}),$ satisfying the following conditions:

\noindent a) $w\mbox{-}\lim\limits_{n\to\infty}{(\mathrm{F}^{[n
s]}(\frac{t}{n})f-f_n^{s})}=0$;

\noindent b)  It is possible to choose a convergent sub-subsequence
for any subsequence of the sequence $\{\mathrm{F}^{[n
s/t]}(\frac{t}{n})f\}_{n\in \Nat}$;

\noindent c) It is possible to choose a weakly convergent
sub-subsequence for any subsequence of the sequence
$\{{\ov{\mmm{Z}}}f_n^{s}\}_{n=1}^{\infty}$.
\end{theorem}
\begin{proof}[Proof of Theorem~\ref{vt:54}.]
($\Rightarrow$) It follows from Theorem~\ref{vvt:54}.

\noindent ($\Leftarrow$) By Theorem~\ref{vvt:54}  there exists the solution of the equation $$f'(s)=\overline{\mmm{Z}}f(s)$$ on
$[0,\,\infty)$ for any initial condition $f(0)=g$, where
$g\in\mmf{A}$. Thus, by Theorem~\ref{cog:1},
$\overline{\mmm{Z}}$ is the generator of the $ lC_0$-semigroups
$\mmm{S}$ and
$$\mmm{S}(t) f=\lim\limits_{n\to\infty}{{\mmm{F}(t/n)}^n
f},\,\,f\in\mmm{X}.$$
\end{proof}

\begin{corollary} \label{vt:55} Let $\mmm{C}$
and $\mmm{D}$ be generators of the $lC_0$-semigroups
$\exp{(s\mmm{C})}$ and $\exp{(s\mmm{D})}$ in a
complete locally convex  space $\mmf{X}$ and $t>0$ be fixed. Assume that the following conditions are satisfied:

\noindent i) $\mcl{D}(\mmm{C})\cap\mcl{D}(\mmm{D})$ is dense in
$\mmf{X}$.

\noindent ii) $\{[0,\,\infty)\ni
s\mapsto\exp{{\textstyle(s\mmm{C})}}\exp{{\textstyle(s\mmm{D})}}\}
\in\EuScript{E}_{\mmf{X}}$.

\noindent Then the sum of $\mmm{C}$ and $\mmm{D}$ possesses a closure and its closure is the generator of the $lC_0$-semigroup
iff there exists a   dense linear subspace
$\mmf{A}\subseteq\mcl{D}(\mmm{C})\cap\mcl{D}(\mmm{D})$ such that for any  $f\in\mmf{A}$ and $s\geq 0$ there
exists a sequence
$\{f_n^{s}\}_{n=1}^{\infty}$,
$f_n^{s}\in\mcl{D}(\ov{\mmm{C}+\mmm{D}}),$ satisfying the following conditions:

\noindent a) $w\mbox{-}\lim\limits_{n\to\infty}{((\exp{{(\textstyle
\frac{t}{n}\mmm{C})}}\exp{{(\textstyle\frac{t}{n}\mmm{D})}})^{[n
s]}f-f_n^{s})}=0$;

\noindent b) It is possible to choose a  convergent sub-subsequence
for any subsequence of the sequence $\{\mathrm{F}^{[n
s/t]}(\frac{t}{n})f\}_{n\in \Nat}$;

\noindent c) It is possible to choose a weakly convergent
sub-subsequence for any subsequence of the sequence
$\{(\ov{\mmm{C}+\mmm{D}})f_n^{s}\}_{n=1}^{\infty}$.

\noindent Furthermore, if conditions (a)-(c) are satisfied, then
$$\exp{(s(\overline{\mmm{C}+\mmm{D}}))} f=\lim\limits_{n\to\infty}{\{\exp{{\textstyle (\frac{s}{n}\mmm{C})}}\exp{{\textstyle (\frac{s}{n}\mmm{D})}}\}^n
f},\, s>0,$$ for each $f\in\mmm{X}.$
\end{corollary}

\begin{proof}[Proof of Corollary~\ref{vt:55}.] Put $\mathrm{F}(s)=\exp(s \mmm{C})\exp(s \mmm{D})$ for any $s\in[0,\,\infty)$.
By
\begin{eqnarray}
 & &\lim\limits_{s \rightarrow 0}{s^{-1}(\mmm{F}(s)f-\mmm{F}(0)f)}\nonumber\\&=&\lim\limits_{s \rightarrow
0}{s^{-1}(\exp(s \mmm{C})\exp(s \mmm{D})f-\exp(s
\mmm{C})f)}\nonumber\\&+& \lim\limits_{s \rightarrow
0}{s^{-1}(\exp(s
\mmm{C})f-f)}=\mmm{C}f+\mmm{D}f,\,\,f\in\mcl{D}(\mmm{C})\cap\mcl{D}(\mmm{D}),\nonumber
\end{eqnarray}
we infer that $\mcl{D}(\mmm{F}_{ef}'(0))$ is dense in $\mmm{X}$.
Thus, by condition (ii), we deduce that
$\mmm{F}\in\EuScript{F}_{\mmf{X}}$. From Theorem~\ref{cor:1} it follows that operator $\mmm{C}+\mmm{D}$ possesses
a closure. Therefore, applying Theorem~\ref{vt:55} for the operator $\mmm{Z}=\mmm{C}+\mmm{D}$, we get the Corollary. \end{proof}

\begin{remark}\label{ref:9.1}
If $\mmf{A}$ is sequentially dense in $\mmf{X}$, then the condition that $\mmf{X}$ is a complete locally convex
space in Theorems~\ref{vt:4}-\ref{vt:54} and Corollaries~\ref{vmain:10}-\ref{vt:55} can be changed to the condition that $\mmf{X}$ is a sequentially complete locally convex space. Indeed, it is enough to apply Theorem~\ref{cog:1} and Remark~\ref{rem:1.1}.
\end{remark}

\section{Consequences for semigroups defined on Banach spaces.}\label{3.3}

If $\mmf{X}$ is a reflexive separable Banach space or Hilbert space, then the formulations of some previous theorems can be considerably simplified
 and we get the following results:
\begin{corollary}\label{newc:1}Let  $\mmm{Z}:\mcl{B}\supset \mcl{D}(\mmm{Z})\to\mcl{B}$
be a densely defined linear operator and $t>0$. Then the operator $\mmm{Z}$ possesses a closure and its closure is the
generator of the $lC_0$-semigroup iff there exists a function
$\mathrm{F}:[0,\infty)\mapsto\mathcal{L}(\mcl{B})$ such that:

\noindent i) $\mathrm{F}(0)=\mathrm{I}$ and there exist
$a\in\mathbb{R}$ and $M\geq 1$ such that
$\|\mathrm{F}^{m}(\frac{s}{n})\|\leq M\exp(a s\frac{m}{n})$ for any $n,\, m\in\mmb{N}$ and $s\geq 0$;

\noindent ii) $\mcl{D}(\mmm{F}_{ef}'(0))\supset\mcl{D}(\mmm{Z})$ and
$\mmm{F}_{ef}'(0)g=\mmm{Z}g,\,\mbox{для}\,g\in\mcl{D}(\mmm{Z})$;

\noindent iii) $\mcl{D}({\mmm{F}^*}'(0))$  is *-dense in $\mcl{B}^*$;

\noindent iv) There exists a dense linear subspace
$\mmf{A}\subseteq\mcl{D}(\mmm{Z})$ such that for any
$f\in\mmf{A}$ and $s\geq 0$ there
exists a sequence
$\{f_n\}_{n=1}^{\infty}$, $f_n\in\mcl{D}(\mmm{Z}),$ such that there exists
$$w\mbox{-}\lim\limits_{n\to\infty}{(\mathrm{F}^{[n
s]}({\textstyle\frac{t}{n}})f-f_n)}=0,$$ and the sequences $\{\mmm{Z}f_n\}_{n=1}^{\infty}$ is bounded.

\noindent Furthermore, if conditions (i)-(iv) are satisfied, then
 $\exp{(s\overline{\mmm{Z}})}
f=\lim\limits_{n\to\infty}{{\mmm{F}(s/n)}^n f}, \,s>0,$ for each
$f\in\mcl{B}.$
\end{corollary}

\begin{proof}[Proof of Corollary~\ref{newc:1}.] Since  $\{(\mathrm{F}^{[n
s]}({\textstyle\frac{t}{n}})f-f_n)\}_{n=1}^{\infty}$ is weak convergent and $\{(\mathrm{F}^{[n
s]}({\textstyle\frac{t}{n}})f)\}_{n=1}^{\infty}$ is bounded we infer that $\{(\mathrm{F}^{[n
s]}({\textstyle\frac{t}{n}})f-f_n)\}_{n=1}^{\infty}$ and $\{f_n\}_{n=1}^{\infty}$ are bounded. Note that for any bounded sequences
$\{g_n\}_{n=1}^{\infty},\,\,g_n\in\mcl{B},$ the following condition is satisfied:
 It is possible to choose a weakly convergent sub-subsequence for any subsequence of the
sequence
$\{g_n\}_{n\in\Nat}$. Therefore, by Theorem~\ref{vmain:1}, we get the Corollary.
\end{proof}

\begin{corollary}\label{newc:2} Let  $\mmm{Z}:\mcl{B}\supset \mcl{D}(\mmm{Z})\to\mcl{B}$ be
 a densely defined linear operator in $\mcl{B}$
and $t>0$. Assume  that there exists a function
$\mathrm{F}:[0,\infty)\mapsto\mathcal{L}(\mcl{B})$ such that:

\noindent i) $\mathrm{F}(0)=\mathrm{I}$ and there exist
$a\in\mathbb{R}$ and $M\geq 1$ such that
$$\|\mathrm{F}^{m}({\textstyle\frac{s}{n}})\|\leq M\exp{({\textstyle a s\frac{m}{n}})}$$
for any $n,\, m\in\mmb{N}$ and $s\geq 0$;

\noindent ii) $\mcl{D}(\mmm{F}_{ef}'(0))\supset\mcl{D}(\mmm{Z})$ and
$\mmm{F}_{ef}'(0)g=\mmm{Z}g,\,\,g\in\mcl{D}(\mmm{Z})$;

\noindent iii) $\mcl{D}({\mmm{F}^*}'(0))$ is *-dense in $\mathbf{B}^*$.

\noindent Then the operator $\mmm{Z}$ possesses a closure and its closure is the generator of the $lC_0$-semigroup iff there exists a dense linear subspace
$\mmf{A}\subseteq\mcl{D}(\mmm{Z})$ such that for any
$f\in\mmf{A}$ and $s\geq 0$ there
exists a sequence
$\{f_n\}_{n=1}^{\infty}$, $f_n\in\mcl{D}(\mmm{Z}),$ such that $w\mbox{-}\lim\limits_{n\to\infty}{(\mathrm{F}^{[n
s]}({\textstyle\frac{t}{n}})f-f_n)}=0$, and  $\{\mmm{F}_{ef}'(0)f_n\}_{n=1}^{\infty}$ is bounded.
\end{corollary}

\begin{proof}[Proof of Corollary~\ref{newc:2}.] Since  $\{(\mathrm{F}^{[n
s]}({\textstyle\frac{t}{n}})f-f_n)\}_{n=1}^{\infty}$ is weak convergent and $\{(\mathrm{F}^{[n
s]}({\textstyle\frac{t}{n}})f)\}_{n=1}^{\infty}$ is bounded we infer that $\{(\mathrm{F}^{[n
s]}({\textstyle\frac{t}{n}})f-f_n)\}_{n=1}^{\infty}$ and $\{f_n\}_{n=1}^{\infty}$ are bounded. Note that for any bounded sequences
$\{g_n\}_{n=1}^{\infty},\,\,g_n\in\mcl{B},$ the following condition is satisfied:
 It is possible to choose a weakly convergent sub-subsequence for any subsequence of the
sequence
$\{g_n\}_{n\in\Nat}$. Therefore, by Theorem~\ref{vt:4}, we get the Corollary.
\end{proof}

\begin{corollary}\label{newc:3} Let  $\mmm{C}$
and $\mmm{D}$ be  generators of the $lC_0$-semigroups
$\exp(s\mmm{C})$ and $\exp(s\mmm{D})$ in $\mcl{B}$ and $t>0$. Assume that the
following conditions are satisfied:

\noindent i) $\mcl{D}(\mmm{C})\cap\mcl{D}(\mmm{D})$ is dense in
$\mcl{B}$;

\noindent ii) $\mcl{D}(\mmm{C}^*)\cap\mcl{D}(\mmm{D}^*)$ is *-dense in
$\mcl{B^*}$;

\noindent iii) There exist $a\in\mathbb{R}$ and $M\geq 1$ such that
$$\|\{\exp{{\textstyle (\frac{s}{n}\mmm{C})}}\exp{{\textstyle (\frac{s}{n}\mmm{D})}}\}^m\|\leq M\exp{({\textstyle a
s\frac{m}{n}})}$$ for any $n,\, m\in\mmb{N}$ and $s>0$.

\noindent Then the sum of $\mmm{C}$ and $\mmm{D}$ possesses a closure and its closure is the generator of the $lC_0$-semigroup iff  there exists a dense linear subspace
$\mmf{A}\subseteq\mcl{D}(\mmm{C})\cap\mcl{D}(\mmm{D})$ such that for any $f\in\mmf{A}$ and $s\geq 0$ there exists a sequence
$\{f_n\}_{n=1}^{\infty}$,
$f_n\in\mcl{D}(\mmm{C})\cap\mcl{D}(\mmm{D}),$ satisfying the following
conditions:

\noindent a) $w\mbox{-}\lim\limits_{n\to\infty}{((\exp{({\textstyle
\frac{t}{n}\mmm{C}})}\exp{({\textstyle\frac{t}{n}\mmm{D}})})^{[n
s]}f-f_n)}=0$;

\noindent b) $\{(\mmm{C}+\mmm{D})f_n\}_{n=0}^{\infty}$ is bounded.

\noindent Furthermore, if conditions (a)-(b) are satisfied, then
$$\exp{(s(\overline{\mmm{C}+\mmm{D}}))} f=\lim\limits_{n\to\infty}{\{\exp{{\textstyle (\frac{s}{n}\mmm{C})}}\exp{{\textstyle
(\frac{s}{n}\mmm{D})}}\}^n f},\, \, s>0,$$ for each
$f\in\mcl{B}.$
\end{corollary}

\begin{proof}[Proof of Corollary~\ref{newc:3}.] Since  $\{((\exp{({\textstyle
\frac{t}{n}\mmm{C}})}\exp{({\textstyle\frac{t}{n}\mmm{D}})})^{[n
s]}f-f_n)\}_{n=1}^{\infty}$ is weak convergent and $\{(\exp{({\textstyle
\frac{t}{n}\mmm{C}})}\exp{({\textstyle\frac{t}{n}\mmm{D}})})^{[n
s]}f\}_{n=1}^{\infty}$ is bounded we infer that $\{((\exp{({\textstyle
\frac{t}{n}\mmm{C}})}\exp{({\textstyle\frac{t}{n}\mmm{D}})})^{[n
s]}f-f_n)\}_{n=1}^{\infty}$ and $\{f_n\}_{n=1}^{\infty}$ are bounded. Note that for any bounded sequences
$\{g_n\}_{n=1}^{\infty},\,\,g_n\in\mcl{B},$ the following condition is satisfied:
 It is possible to choose a weakly convergent sub-subsequence for any subsequence of the
sequence
$\{g_n\}_{n\in\Nat}$. From Theorem~1.34 in~\cite{Dav} it follows that
$\exp(s\mmm{C}^*)$ and  $\exp(s\mmm{D}^*)$ are $lC_0$-semigroups and, therefore, condition (iv) of Corollary~\ref{vt:3} is satisfied for
$\mmf{B}=\mcl{D}(\mmm{C}^*)\cap\mcl{D}(\mmm{D}^*)$. Thus, by Corollary~\ref{vt:3}, we get the Corollary.
\end{proof}

\begin{remark} \label{rem:5} If the operators $i\mmm{C}$ and $i\mmm{D}$ are self-adjoint, then conditions (ii), (iii) of Corollary~\ref{newc:3} can be omitted.
\end{remark}

\begin{corollary}\label{chern1hilbert}
Assume that $\mcl{B}$ is a Hilbert space and  $\mmm{Z}$ is
 a densely defined dissipative linear operator in $\mcl{B}$. Then there exists a function $\mmm{F}\in\EuScript{F}_{\mcl{B}}$ such that
$$F'(0)|_{\mcl{D}(\mmm{Z})}=\mmm{Z}.$$ Moreover, If there exist a local solution $f:[0,\,l)\mapsto\mcl{B}$ of the system
\begin{eqnarray}
f'(t)&=&\overline{\mmm{Z}}f(t),\,\,t\in[0,\,l),\nonumber\\
f(0)&=&h\in\mcl{D}(\overline{\mmm{Z}}),\nonumber
\end{eqnarray}
then $\mathrm{F}({\textstyle\frac{t}{n}})^{n}h$ tends to
$f(t)$ as $n\to\infty$ uniformly with respect to $t\in [0,\,t_0]$ for any
$t_0\in(0,\,l)$.
\end{corollary}
\begin{proof}[Proof of Corollary~\ref{chern1hilbert}.]
It follows from Theorem~\ref{chern1} and Proposition~\ref{prop:11}.
\end{proof}


\begin{thebibliography}{99}

\bibitem{Albanese0}
A. A. Albanese, F. K\"{u}hnemund, Trotter-Kato approximation
theorems for locally equicontinuous semi-groups, Riv. Mat. Univ.
Parma 1 (2002) 19-53.
\bibitem{Albanese}
A. A. Albanese, E. Mangino, Trotter-Kato theorems for bi-continuous
semigroups and applications to Feller semigroups, J. Math. Anal.
Appl. 289 (2004) 477-492.
\bibitem{Chernoff}
P. R. Chernoff, Note on product formulas for operator semigroups, J.
Funct. Anal. 2 (1968), 238-242.

\bibitem{Dav}
E. B. Davies, One-Parameter Semigroups, St. John's College, Oxford,
England (1980).
\bibitem{Engel}
K.-J. Engel and R. Nagel, One-Parameter Semigroups for Linear
Evolution Equations, Graduate Texts. in Math., vol. 194,
Springer-Verlag, 2000.

\bibitem{Grosu}
C. Grosu, The Trotter product formula in locally convex spaces, I.
Rev. Roumaine Math. Pures Appl.- 31.- 1986, no. 1, p.29-42.

\bibitem{Komatsu}
H.  Komatsu,  Semi-groups  of  operators  in  locally  convex
spaces,  J.  Math.  Soc. Japan,  16  (1964),  230-262.
\bibitem{Komura}
T.  Komura,  Semigroups  of  operators  in  locally  convex  spaces,
J.  Functional Analysis,  2  (1968),  258-296.
\bibitem{Kuhnemund}
F. K\"{u}hnemund, Bi-continuous semigroups on spaces with two
topologies: theory and applications, Ph.D. thesis,
T\"{u}bingen, 2001.
\bibitem{Kuhnemund2}
F. K\"{u}hnemund, A Hille–Yosida theorem for bi-continuous
semigroups, Semigroup Forum, submitted for publication.
\bibitem{Kuhnemund3}
F. K\"{u}hnemund, Approximation of bi-continuous semigroups, J.
Approx. Theory, submitted for publication.
\bibitem{Nekl2} A. Y. Neklyudov,
Chernoff and Trotter type product formulas, Mat. Sb., 200:10 (2009), 81-106


\bibitem{Nelson}
E. Nelson, Feynman integrals and the Schrodinger equation, J. Math.
Phis., 1964, Vol. 5, no. 3, p. 332-343.
\bibitem{Ouchi}
S.  Ouchi, Semi-groups  of  operators  in  locally  convex spaces,
Math.  Soc. Japan, Vol.  25, No.  2,  1973.

\bibitem{Pettis}
B.J. Pettis,   "On integration in vector spaces"  Trans. Amer. Math. Soc., 44:2 (1938) p. 277-304.
\bibitem{Schwartz}
L.  Schwartz,  Lectures  on  mixed  problems  in  partial
differential  equations  and the  representation  of  semi-groups,
Tata  Institute  of  Fundamental  Research,  1958.
\bibitem{Hamil}
O. G. Smolyanov, A. G. Tokarev, A. Truman, Hamiltonian Feynman path
integral via the Chernoff formula, J. Math. Phys., Vol. 43, No. 10,
October 2002.
\bibitem{Witt}
O. G. Smolyanov, H. von Weizsacker, O. Wittich, Brownian Motion on a
Manifold as Limit of Stepwise Conditioned Standard Brownian Motions,
Canadian Math. Society Conference Proceedings, 589-602, 29 (2000).


\bibitem{Trotter}
H. Trotter, On the product of semigroups of operators, Proc. Amer.
Math. Soc. 10 (1959), 545-551.
\bibitem{Yosida}
K.  Yosida,  Functional  analysis,  Springer,  Berlin,  1965.




\end{thebibliography}
\end{document}